\documentclass{lmcs}
\pdfoutput=1

\usepackage{lastpage}
\lmcsdoi{21}{2}{12}
\lmcsheading{}{\pageref{LastPage}}{}{}%
{May~24,~2024}{May~16,~2025}{}


\usepackage{hyperref}
\usepackage{tikz}
\usepackage{tikz-cd}
\usetikzlibrary{matrix} 
\usetikzlibrary{arrows}
\usetikzlibrary{arrows.meta}
\usetikzlibrary{decorations.pathmorphing,shapes}
\usetikzlibrary{decorations.markings}
\usetikzlibrary{automata}

\tikzset{
  ->, 
  >=stealth', 
  node distance=2cm, 
  inner sep=1mm,
  every edge/.append style={thick},
  every node/.style={thick,font=\scriptsize},
  every state/.style={font=\normalsize,minimum size=2em}, 
  machine/.style={rectangle,draw,minimum size=2em},
  initial text=$ $, 
}

\usepackage[T1]{fontenc}
\usepackage[utf8]{inputenc}
\usepackage{amssymb}
\usepackage{stmaryrd}
\usepackage{mathpartir}
\usepackage{etoolbox}
\usepackage{wrapfig}

\tikzset{spanmap/.style={
            decoration={markings,
            mark= at position 0.5 with {
                  \node[transform shape] (tempnode) {$|$};
                  }
              },
              postaction={decorate}
}
}

\tikzset{doublespanmap/.style={
            decoration={markings,
            mark= at position 0.5 with {
                  \node[transform shape] (tempnode) {$||$};
                  }
              },
              postaction={decorate}
}
}

\DeclareMathAlphabet\mathbfcal{OMS}{cmsy}{b}{n}

\newcommand\ccat[1]{\mathcal{#1}}
\newcommand\cgraph[1]{\mathbb{#1}}

\newcommand\coper[1]{\mathcal{#1}}
\newcommand\cspec[1]{\mathbb{#1}}

\newcommand{\Dependency}[1]{\mathit{D}}
\newcommand{\Independency}[1]{\mathit{I}}
\newcommand{\permeq}[1]{\sim}

\newcommand\cC{\ccat{C}}
\newcommand\cD{\ccat{D}}
\newcommand\cE{\ccat{E}}
\newcommand\cQ{\ccat{Q}}

\newcommand\oD{\coper{D}}

\newcommand\oO{\coper{O}}
\newcommand\oP{\coper{P}}
\newcommand\oQ{\coper{Q}}

\newcommand\sS{\cspec{S}}
\newcommand\sR{\cspec{R}}

\newcommand\N{{\mathbb N}}
\newcommand\TBase{{\mathcal T}}

\newcommand\bigcat[1]{\mathrm{#1}}
\newcommand\Set{\bigcat{Set}}
\newcommand\FinSet{\bigcat{FinSet}}
\newcommand\Span[1]{\bigcat{Span}(#1)}
\renewcommand\Subset{\bigcat{Subset}}

\newcommand\Rel{\bigcat{Rel}}
\newcommand\Cat{\bigcat{Cat}}
\newcommand\Graph{\bigcat{Grph}}

\newcommand\Dist{\bigcat{Dist}}

\newcommand\DiGr{\mathrm{DiGr}}
\newcommand\series{\mathbin{;}}

\newcommand\finULF{\mathrm{finULF}}

\newcommand\TranSet{\mathit{Tran}}

\newcommand\Monoid[1][]{\ccat{B}_{#1}}
\newcommand\Bouquet[1][]{\cgraph{B}_{#1}}
\newcommand\Bracket[1]{#1^{\bow\eow}}

\newcommand\eow{\$}
\newcommand\eowOb{\top}
\newcommand\bow{\textrm{\^{}}}
\newcommand\bowOb{\bot}

\newcommand\set[1]{\{\,{#1}\,\}}

\newcommand\Words[1]{\coper{W}\,{#1}}

\newcommand\Wordsonly{\coper{W}}
\newcommand\Wordsbulletonly{\coper{W}_\bullet}

\newcommand\defin[1]{\textbf{#1}}
\newcommand\defeq{:=}

\newcommand\Contour[1]{\ccat{C}\,{#1}}
\newcommand\Contouronly{\ccat{C}}
\newcommand\UnivGrammar[1]{\mathsf{Univ}_{#1}}

\newcommand{\Ccategory}{\ccat{C}}

\newcommand{\Ooperad}{\coper{O}}

\newcommand{\Qoperad}{\coper{Q}}
\newcommand{\Qcategory}{\ccat{Q}}

\newcommand{\Sspecies}{\cspec{S}}
\newcommand{\Rspecies}{\cspec{R}}
\newcommand{\Tspecies}{\cspec{T}}

\newcommand\AffMon[1]{!_{\mathsf{aff}}{#1}}
\newcommand\CartMon[1]{!_{\mathsf{cart}}{#1}}

\newcommand\Free{\mathcal{F}}
\newcommand\FreeOper[2][]{\Free_{#1}\,#2}
\newcommand\FreeOperFunctor{\Free}
\newcommand\Forget{\mathbb{U}}
\newcommand\ForgetOper[1]{{#1}}
\newcommand\ForgetOperFunctor{\Forget}

\newcommand\FreeCat[1]{\Free\,{#1}}

\makeatletter
\newcommand\FactorOperad[2][\@nil]{\def\tmp{#1}\ifx\tmp\@nnil\mathsf{Fact}_{\,#2}\else\mathsf{Fact}_{\,#1,#2}\fi}
\makeatother

\newcommand\Fact[1]{\mathsf{Fact}_{#1}}

\newcommand\TerminalCat{\mathbf{1}}

\newcommand\TerminalSpecies{\mathbb{N}}

\newcommand\elts{\mathop{\mathrm{con}}}

\newcommand\Nonterminals{N}
\newcommand\Productions{P}
\newcommand\Sentence{S}

\newcommand\Lang[2][]{\mathcal{L}_{#2}^{#1}}

\newcommand\src{\mathsf{src}}
\newcommand\tgt{\mathsf{tgt}}

\newcommand{\id}[1][]{id_{#1}}

\newcommand\ContourLang[2][]{\mathcal{C\!L}_{#2}^{#1}}
\newcommand\DyckLang[1]{\mathcal{D\!L}_{#1}}

\newcommand{\Species}[1]{\bigcat{Spec}_{#1}}
\newcommand{\Operads}[1]{\bigcat{Oper}_{#1}}

\newcommand\op{{\mathrm{op}}}

\newcommand\node[1][]{\mathbin{\mathrm{node}}_{#1}}

\newcommand{\speciesunit}[1]{\cspec{I}_{#1}}

\newcommand{\cod}{\mathrm{cod}}
\newcommand{\dom}{\mathrm{dom}}

\newcommand{\pull}[1]{\mathop{#1^*}}
\newcommand{\push}[1]{\mathop{#1_*}}
\newcommand{\spanmap}{{\,\longrightarrow\hspace{-1.15em}|\hspace{.95em}}}

\newcommand\bin{{\mathsf{bin}}}
\newcommand\nodes{{\mathsf{nodes}}}

\newcommand\colors{{\mathsf{colors}}}

\newcommand\tagU{u}
\newcommand\tagD{d}

\newcommand\drefs[1][]{\overset{#1}\sqsubset}
\newcommand\nrefs[1][]{\sqsubset^{#1}}

\newcommand\dseq[2][]{\overset{#1}{\underset{#2}{\Longrightarrow}}}
\newcommand\nseq[2][]{\Rightarrow_{#2}^{#1}}

\usepackage{scalerel}
\makeatletter
\newcommand\refs[1][]{\ThisStyle{\if D\m@switch\drefs[{#1}]\else\nrefs[{#1}]\fi}}
\newcommand\seq[2][]{\ThisStyle{\if D\m@switch\dseq[{#1}]{#2}\else\nseq[{#1}]{#2}\fi}}
\makeatother

\newcommand{\TraceCat}{\FreeCat{{\Bouquet[\Sigma]}}/{\sim_{\,\Dependency{}}}}

\newcommand{\els}[2][]{\smallint_{#1}#2}

\theoremstyle{plain} 
\theoremstyle{definition}\newtheorem{cexa}[thm]{Counterexample}
\theoremstyle{plain}

\begin{document}

\title[The categorical contours of the C-S representation theorem]{The categorical contours of the Chomsky-Sch{\"u}tzenberger representation theorem}
\titlecomment{This is a thoroughly revised and expanded version of a paper with a similar title \cite{MZ2022mfps} presented at MFPS 2022. 
    The authors are partially supported by ANR projects CoREACT (ANR-22-CE48-0015), LambdaComb (ANR-21-CE48-0017), and ReciProg (ANR-21-CE48-0019).}

\author[P.-A.~Melli\`es]{Paul-Andr\'e Melli\`es\lmcsorcid{0000-0001-6180-2275}}
\author[N.~Zeilberger]{Noam Zeilberger\lmcsorcid{0000-0002-5945-4184}}

\address{IRIF, Université Paris Cité, CNRS, Inria, Paris, France}
\email{paul-andre.mellies@cnrs.fr}

\address{LIX, École Polytechnique, Inria, Palaiseau, France}
\email{noam.zeilberger@polytechnique.edu}

\begin{abstract}
We develop fibrational perspectives on context-free grammars and on nondeterministic finite-state automata over categories and operads.
A generalized CFG is a functor from a free colored operad (aka multicategory) generated by a pointed finite species into an arbitrary base operad:
this encompasses classical CFGs by taking the base to be a certain operad constructed from a free monoid, as an instance of a more general construction of an \emph{operad of spliced arrows} $\Words{\cC}$ for any category $\cC$. 
A generalized NFA is a functor from an arbitrary bipointed category or pointed operad satisfying the unique lifting of factorizations and finite fiber properties:
this encompasses classical word automata and tree automata without $\epsilon$-transitions, but also automata over non-free categories and operads.
We show that generalized context-free and regular languages satisfy suitable generalizations of many of the usual closure properties, and in particular we give a simple conceptual proof that context-free languages are closed under intersection with regular languages.
Finally, we observe that the splicing functor $\Wordsonly : \Cat \to \Operads{}$ admits a left adjoint $\Contouronly:\Operads{}\to\Cat$,
which we call the \emph{contour category} construction since the arrows of $\Contour{\oO}$ have a geometric interpretation as oriented contours of operations of $\oO$.
A direct consequence of the contour / splicing adjunction is that every pointed finite species induces a universal CFG generating a language of \emph{tree contour words.}
This leads us to a generalization of the Chomsky-Schützenberger Representation Theorem, establishing that a subset of a homset $L \subseteq \cC(A,B)$ is a CFL of arrows if and only if it is a functorial image of the intersection of a $\Ccategory$-chromatic tree contour language with a regular language.

\end{abstract}

\maketitle

\section*{Introduction}
\label{section/introduction}

In ``Functors are Type Refinement Systems'' \cite{mz15popl}, we argued for the idea that rather than being modelled merely as categories, type systems
should be modelled as functors $p:\ccat{D}\to\ccat{T}$ from a category~$\ccat{D}$ whose morphisms may be considered as abstract typing derivations
to a category~$\ccat{T}$ whose morphisms are the terms corresponding to the underlying \emph{subjects} of those derivations.
One consequence of this fibrational point of view is that the notion of typing judgment receives a simple mathematical status, as a triple~$(R,f,S)$
consisting of two objects~$R,S$ in $\ccat{D}$ and a morphism $f$ in $\TBase$ such that $p(R)=\dom(f)$ and $p(S)=\cod(f)$.
The question of finding a derivation for such a judgment then reduces 
to the ``lifting problem'' of finding a morphism $\alpha:R\to S$ such that $p(\alpha)=f$.
Note that this problem could in general have zero, one, or more solutions.

We developed this perspective in a series of papers \cite{mz2013trmcb,mz15popl,mz16bifib,mz18isbell}, and believe that it may be usefully applied to a large variety of deductive systems, beyond type systems in the traditional sense.
In this work, we focus on derivability in context-free grammars and on recognition by nondeterministic finite-state automata -- two classic topics in formal language theory with wide applications in computer science.
Although grammars and automata are of a very different nature, both may be naturally described as certain kinds of functors.

Free operads\footnote{In this paper, when we say ``operad'' we always mean \emph{colored operad}, in the sense that our operads always carry a (potentially trivial) set of colors and operations are sorted $f : A_1,\dots,A_n \to B$ in addition to having an arity $n$.  An operad in this sense is also called a \emph{multicategory}, and we use the two words interchangeably.} play an important role in our treatment of context-free grammars.
Indeed, we will ultimately propose that a ``generalized CFG'' may be defined simply as
\emph{a functor from a free operad generated by a pointed finite species} into an arbitrary base operad.
In order to recover the classical definition of CFG, we rely on a construction $\Wordsonly : \Cat \to \Operads{}$ turning any category $\cC$ into an operad $\Words{\cC}$, which we call its \emph{operad of spliced arrows}.
One recovers the classical definition by taking $\cC = \FreeCat{\Bouquet[\Sigma]}$ to be the one-object category freely generated from the one-node ``bouquet'' graph with a loop $a : * \to *$ for every letter of the alphabet $a \in \Sigma$.
However, as we will see, it is very natural to consider context-free languages of arrows in an arbitrary category, or even generalized context-free languages of constants in an arbitrary operad.

Nondeterministic finite-state automata, it turns out, may be naturally modelled as \emph{functors satisfying the unique lifting of factorizations and finite fiber properties}.
One nice aspect of this formulation is that it captures both word automata and tree automata, which correspond to functors satisfying these two properties into free categories and into free operads respectively.
However, again we will see that it is very natural and at times essential to work with a more general definition of NFA over a non-free category or operad.

A benefit of expressing both context-free grammars and nondeterministic finite-state automata as functors is that they are placed within a common framework, which facilitates comparison and combination.
The germ of this paper was planted early on in the wake of our type refinement systems work, in 2017, when we began considering context-free grammars as functors of multicategories.
But ultimately it was only in 2022 that the key mathematical definitions emerged very quickly (over the span of a couple weeks), when we mused that it could be interesting to analyze the Chomsky-Sch{\"u}tzenberger Representation~Theorem,
a classic result in formal language theory 
that in its now-standard formulation
states that a language is context-free if and only if it is a homomorphic image of the intersection of a Dyck language of well-bracketed words with a regular language.
One reason the theorem is interesting is that it invokes a non-trivial closure property of context-free languages, namely closure under intersection with regular languages.
Another reason is that it suggests, intuitively, that Dyck languages are in some sense ``universal'' context-free languages.

We will see that both of these aspects of the theorem have good categorical explanations.
First, the intersection of a context-free language with a regular language may be computed by first taking the \emph{pullback} of the functor representing the corresponding context-free grammar along the corresponding nondeterministic finite-state automaton, which for non-trivial but elementary reasons defines a context-free grammar over the \emph{runs} of the automaton.
Second, the aforementioned ``spliced arrows'' construction $\Wordsonly$ extends to an \emph{adjunction} between the category of categories and the category of operads, 
\begin{equation*}
\begin{tikzcd}
\Operads{}
\arrow[rr,"{\Contouronly}",yshift=1.5ex]
& \bot &
\Cat\arrow[ll,"{\Wordsonly}",yshift=-1.5ex]
\end{tikzcd}
\end{equation*}
where the left adjoint builds what we call the ``contour category'' $\Contour{\oO}$ of an operad $\oO$.
This contour / splicing adjunction, which seems to be of independent interest, has as an immediate consequence that every pointed finite species induces a universal CFG generating a language of \emph{tree contour words.}
Such contour words are closely related to Dyck words, and allow us to prove a generalization of the representation theorem for context-free languages of arrows in an arbitrary category.

The rest of the paper is organized as follows.
Section~\ref{section/cfls} is devoted to context-free grammars, Section~\ref{section/nfas} to nondeterministic finite-state automata, and Section~\ref{section/csrep} to the representation theorem.
We discuss related work at the end of each of these sections, and we conclude in Section~\ref{section/conclusion}.

\section{Context-free languages of arrows in a category}
\label{section/cfls}

In this section, after some preliminaries on free operads and species, we describe the operad of spliced arrows construction, explain how to use it to define context-free grammars over arbitrary categories generating \emph{context-free languages of arrows,} and show that many standard properties of CFGs and CFLs may be naturally reformulated in this general setting.
We also spend some time developing the fibrational view of parsing as a lifting problem, explaining how any functor from a free operad $\FreeOper{\sS} \to \oO$ may be represented as a ``displayed free operad'' corresponding to a lax functor $\oO \to \Span\Set$ that admits a certain inductive definition.
Finally, we describe how the definition of categorical CFG admits a further natural generalization by allowing the base operad to be arbitrary, and show that this encompasses several extensions of the notion of context-free grammar from the literature.

\subsection{Preliminaries on free operads and species}
\label{section/preliminaries-operads-species}

We assume some basic familiarity with operads, otherwise known as multicategories, as introduced for example in Chapter~2 of Leinster's book \cite{LeinsterHOHC}.
We write $f \circ (g_1,\dots,g_n)$ for parallel composition of operations, and $f \circ_i g$ for partial composition \emph{after} the first $i$ inputs (so $f \circ_0 g$ denotes composition of $g$ into the 0th input of $f$).
Here we recall the notion of a (colored non-symmetric) species, and how one gives rise to a free operad.

\begin{defi}\label{definition/species}
A \emph{colored non-symmetric species,} which we abbreviate to \defin{species} for short, is a tuple $\sS = (C, V, i, o)$ consisting of a span of sets
$\begin{tikzcd}
  C^* & \arrow[l,"i"']V\arrow[r,"o"] & C
\end{tikzcd}$
with the following interpretation: $C$ is a set of \emph{colors,} $V$ is a set of \emph{nodes,} and $i : V \to C^*$ and $o : V \to C$ return respectively the list of input colors and the unique output color of each node.
\end{defi}
\noindent
Species in this sense could also be called ``signatures'', and are sometimes called ``multigraphs'' \cite{Lambek1987,Walters1989} since they bear the same relationship to multicategories as graphs do to categories (though that terminology has an unfortunate clash with a different concept in graph theory); we use ``species'' to emphasize the link with Joyal's theory of (symmetric) species \cite{Joyal1981} and with generalized species \cite{FGHW2008}.
Adopting the same notation as we use for operations of an operad, we write $x : R_1,\dots,R_n \to R$ to indicate that $x \in V$ is a node with list of input colors $i(x) = (R_1,\dots,R_n)$ and output color $o(x) = R$.
However, it should be emphasized that a species by itself only contains bare coloring information about the nodes, and does not say how to compose them as operations.
We will be primarily interested in \emph{finite} species.
\begin{defi}    
We say that a species is \defin{finite} (also called \emph{polynomial} \cite{Joyal1986}) just in case both sets $C$ and $V$ are finite.       
\end{defi}

A map of species $\phi : \sS \to \Tspecies$ from $\sS = (C,V,i,o)$ to $\Tspecies = (D,W,i',o')$ is given by a pair $\phi = (\phi_C,\phi_V)$ of functions $\phi_C : C \to D$ and $\phi_V : V \to W$ making the diagram commute:
\[
\begin{tikzcd}
  C^*\arrow[d,"\phi_C^*"] & \arrow[l,"i"']V\arrow[r,"o"]\arrow[d,"\phi_V"] & C \arrow[d,"\phi_C"]\\
  D^* & \arrow[l,"i'"']W\arrow[r,"o'"] & D
\end{tikzcd}
\]
Equivalently, writing $\phi$ for both $\phi_C$ and $\phi_V$, every node $x : R_1,\dots,R_n \to R$ of $\sS$ must be sent to a node $\phi(x) : \phi(R_1),\dots,\phi(R_n) \to \phi(R)$ of $\Tspecies$.
Just as every category has an underlying graph, every operad $\oO$ has an underlying species with the same colors and whose nodes are the operations of $\oO$, and this extends to a forgetful functor $\ForgetOperFunctor:\Operads{}\to\Species{}$ from the category of operads and functors of operads to the category of species and maps of species.
Moreover, this forgetful functor has a left adjoint
\begin{equation*}\label{equation/adjunction-species-operad}
\begin{tikzcd}
\Species{}
\arrow[rr,"\FreeOperFunctor",yshift=1.5ex]
& \bot &
\Operads{}\arrow[ll,"\ForgetOperFunctor",yshift=-1.5ex]
\end{tikzcd}
\end{equation*}
which sends any species $\sS$ to a free operad with the same set of colors and whose operations are generated from the nodes of $\sS$, analogous to the free category over a graph.
By the universal property of the adjoint pair, there is a natural isomorphism of homsets
\[
  \Operads{}(\FreeOper{\sS}, \oO) \cong \Species{}(\sS,{\ForgetOperFunctor}{\oO})
\]
placing functors of operads $p : \FreeOper{\sS} \to \oO$ and maps of species $\phi : \sS \to {\ForgetOperFunctor}{\oO}$ in one-to-one correspondence.
In the sequel, we will leave the action of the forgetful functor implicit, writing $\oO$ for both an operad and its underlying species ${\ForgetOperFunctor}\oO$.

\subsection{The operad of spliced arrows of a category}
\label{section/spliced-arrows}

\begin{figure}
  \begin{center}\includegraphics[width=1\textwidth]{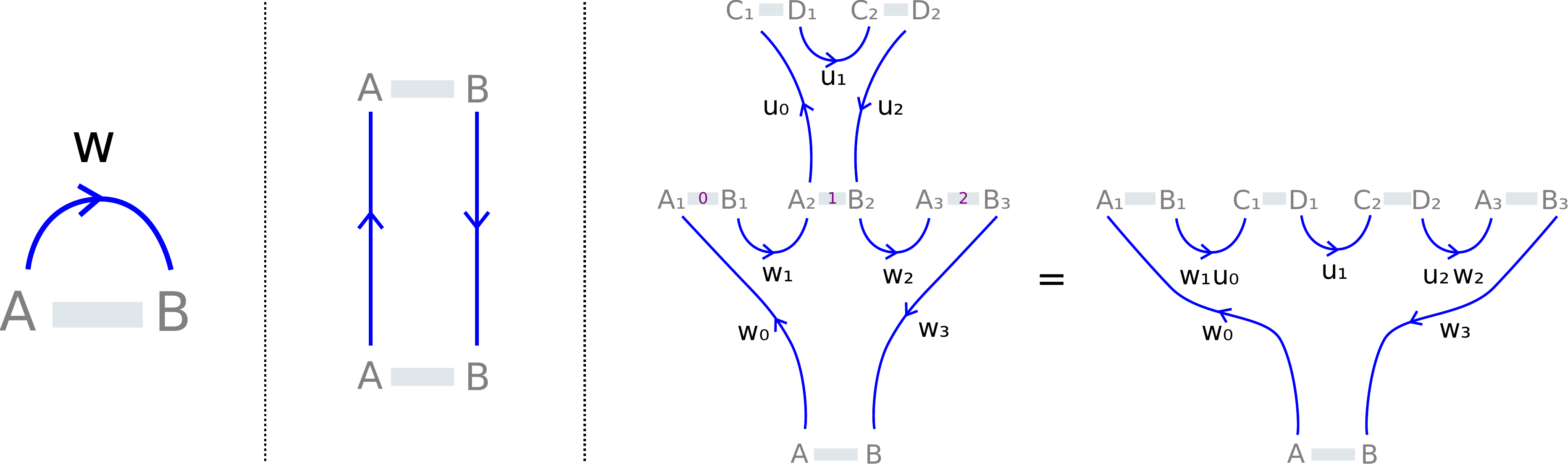}\end{center}
  \caption{Left: a constant of $\Words{\cC}$. Middle: an identity operation. Right: example of partial composition, plugging $g = u_0{-}u_1{-}u_2 : (C_1,D_1),(C_2,D_2) \to (A_2,B_2)$
    into gap~1 of $f = w_0{-}w_1{-}w_2{-}w_3 : (A_1,B_1),(A_2,B_2),(A_3,B_3) \to (A,B)$ to obtain $f \circ_1 g = w_0{-}w_1u_0{-}u_1{-}u_2w_2{-}w_3 : (A_1,B_1),(C_1,D_1),(C_2,D_2),(A_3,B_3) \to (A,B)$.}
  \label{fig:spliced-arrows}
\end{figure}

\begin{defi}\label{definition/spliced-arrows}
Let $\cC$ be a category.
The \defin{operad $\Words{\cC}$ of spliced arrows in $\cC$} is defined as follows:
\begin{itemize}
\item its colors are pairs $(A,B)$ of objects of $\cC$;
\item its $n$-ary operations $(A_1,B_1),\dots,(A_n,B_n) \to (A,B)$ consist of sequences $w_0{-}w_1{-}\dots{-}w_n$ of $n+1$ arrows in $\cC$ separated by $n$ gaps denoted by ``$-$''s, where each arrow must have type $w_i : B_i \to A_{i+1}$ for $0 \le i \le n$, under the convention that $B_0 = A$ and $A_{n+1} = B$;
\item composition of spliced arrows is performed by ``splicing into the gaps'': formally, the partial composition $f \circ_i g$ of a spliced arrow $g = u_0{-}\dots{-}u_m$ into another spliced arrow $f = w_0{-}\dots{-}w_n$ is defined by substituting $g$ for the $i$th occurrence of $-$ in $f$ (starting from the left using 0-indexing) and interpreting juxtaposition by sequential composition in $\cC$ (see Fig.~\ref{fig:spliced-arrows} for an illustration); 
\item the identity operation on $(A,B)$ is given by $\id[A]{-}\id[B]$.
\end{itemize}
\end{defi}
\noindent
It is routine to check that $\Words{\cC}$ satisfies the associativity and neutrality axioms of an operad, these reducing to associativity and neutrality of composition of arrows in $\cC$.
Indeed, the spliced arrows operad construction defines a functor
$ \Wordsonly : \Cat \to \Operads{} $
since any functor of categories $F : \cC \to \cD$ induces a functor of operads
$\Words{F} : \Words{\cC} \to \Words{\cD}$,
acting on colors by $(A,B) \mapsto (FA,FB)$ and on operations by
$w_0{-}\dots{-}w_n \mapsto Fw_0{-}\dots{-}Fw_n$.

\begin{exa}\label{example/spliced-words}
  As referenced in the Introduction, for any alphabet $\Sigma$ there is a one-object category $\FreeCat{\Bouquet[\Sigma]}$, freely generated from the bouquet graph with a loop $a : * \to *$ for every letter of $a \in \Sigma$, so that words $w \in \Sigma^*$ may be regarded as arrows $w : \ast \to \ast$ in $\FreeCat{\Bouquet[\Sigma]}$.
  In this case, the construction yields an operad $\Words{\Sigma} := \Words{\FreeCat{\Bouquet[\Sigma]}}$ that we like to call the \defin{operad of spliced words} in $\Sigma$.
  It has a single color, and its $n$-ary operations are sequences $w_0{-}\dots{-}w_n$ of $n+1$ words over $\Sigma$.
  For example, in the operad of spliced words over the ASCII alphabet, we have that
  $\texttt{\color{blue}"Hell"}{-}\texttt{\color{blue}", "}{-}\texttt{\color{blue}"rld!"} \circ (\texttt{\color{blue}"o"},\texttt{\color{blue}"Wo"}) = 
  \texttt{\color{blue}"Hello, World!"}$.
\end{exa}
\begin{rem}\label{remark/non-free}
  Although $\FreeCat{\Bouquet[\Sigma]}$ is a free category, this property of being freely generated does not extend to its operad of spliced words.
  Indeed, an operad of spliced arrows $\Words{\cC}$ is almost \emph{never} a free operad. 
  That is because any pair of objects $A$ and $B$ induces a binary operation $\id[A]{-}\id[A]{-}\id[B] : (A,A),(A,B) \to (A,B)$,
  and any arrow $w : A \to B$ of $\cC$ induces a corresponding constant $w : (A,B)$.
  Since $\id[A]{-}\id[A]{-}\id[B] \circ (\id[A],w) = w$, $\Words{\cC}$ cannot be a free operad except in the trivial case where $\cC$ has no objects and no arrows.
\end{rem}

\begin{exa}\label{example/spliced-ordinal-sum}
  Arrows in free categories with more than one object may be seen as ``typed'' words, in the sense that concatenation (= composition) is only a partial operation that combines a word $u : A \to B$ with a word of compatible type $v : B \to C$ to produce a word $uv : A \to C$.
  For example, consider the graph
  \[
  \Bracket{\Bouquet[\Sigma]} \quad\defeq\quad
  \begin{tikzcd}
    \bowOb \ar[r,"\bow"] & * \ar[loop above,"a\in\Sigma"] \ar[r,"\eow"] & \eowOb
  \end{tikzcd}
  \]
  obtained by extending the bouquet graph $\Bouquet[\Sigma]$ into a ``cylinder''.
  Arrows $w : \ast \to \ast$ of $\FreeCat{\Bracket{\Bouquet[\Sigma]}}$ correspond to ordinary words, which may be concatenated freely, but arrows $u : \bowOb \to \ast$ may only be prepended to the beginning of a word, while arrows $v : \ast \to \eowOb$ may only be appended to the end of a word.
  As we will see in examples later on, it becomes interesting to consider the spliced arrow operad over $\FreeCat{\Bracket{\Bouquet[\Sigma]}}$, which includes operations of the form
  $ f = w_0{-}\dots{-}w_n\eow : (\ast,\ast), \dots, (\ast,\ast) \to (\ast,\eowOb) $ and $g = \bow w_0{-}\dots{-}w_n : (\ast,\ast),\dots,(\ast,\ast) \to (\bowOb,\ast)$ that may be seen as spliced words with explicit ``end of input'' and ``beginning of input'' markers.
\end{exa}

\begin{rem}
  The operad $\Words{\TerminalCat}$ of spliced arrows over the terminal category is isomorphic to the terminal operad, with a single color $(\ast,\ast)$, and a single $n$-ary operation $\id{-}\dots{-}\id : (\ast,\ast),\dots,(\ast,\ast) \to (\ast,\ast)$ of every arity $n$.
  Likewise, the operad of spliced arrows over the product of two categories decomposes as a product of spliced arrow operads $\Words{(\cC\times \cD)} \cong \Words{\cC}\times\Words{\cD}$.
  This might suggest that the functor $\Wordsonly : \Cat \to \Operads{}$ is a right adjoint, and we will see in \S\ref{section/contour-category} that this is indeed the case.
\end{rem}


\subsection{Context-free grammars over a category}
\label{section/context-free-grammar}

\begin{figure}
  \includegraphics[width=\textwidth]{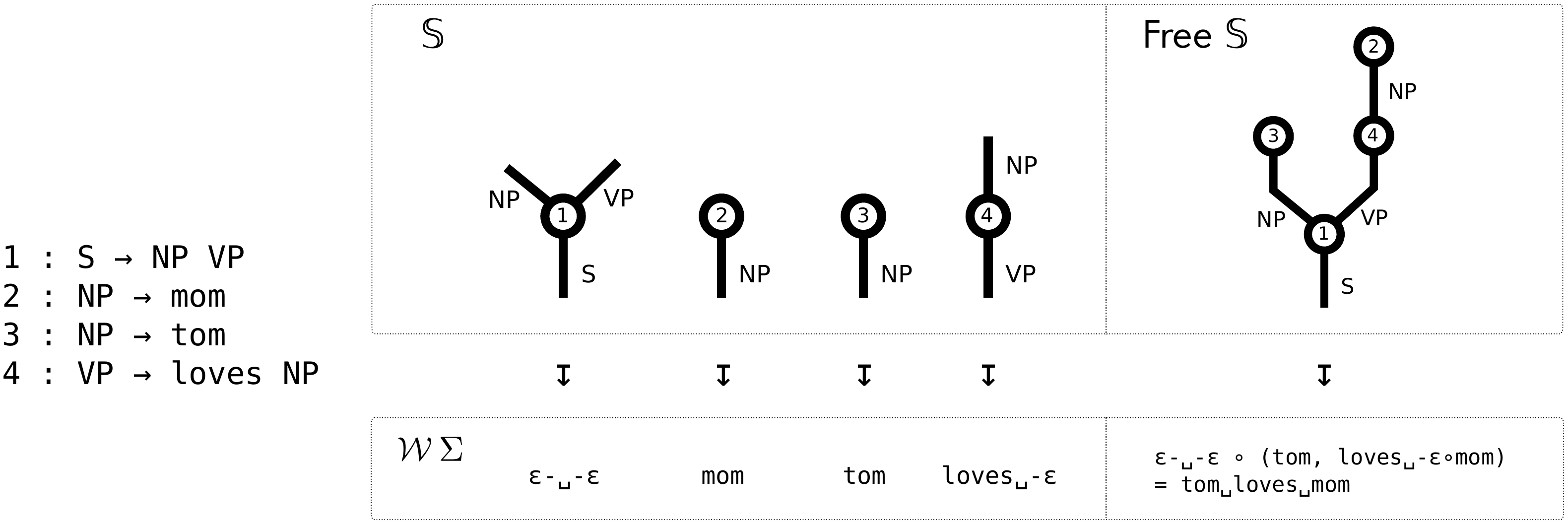}
  \caption{Example of a context-free grammar represented by a functor $\FreeOper{\sS} \to \Words{\Sigma}$, where we have indicated the action of the functor on the generators as well as the induced action on a closed derivation tree.}
  \label{fig:example-cfg}
\end{figure}

Classically, a context-free grammar is defined as a tuple $G = (\Sigma,\Nonterminals,\Sentence,\Productions)$ consisting of
a finite set $\Sigma$ of terminal symbols,
a finite set $\Nonterminals$ of nonterminal symbols,
a distinguished nonterminal $\Sentence \in \Nonterminals$ called the start symbol,
and a finite set $\Productions$ of production rules of the form $R\to\sigma$ where $R\in\Nonterminals$
and $\sigma\in(\Nonterminals\cup\Sigma)^{\ast}$ is a string of terminal or nonterminal symbols.
Observe that any sequence $\sigma$ on the right-hand side of a production can be factored as $\sigma=w_0 R_1 w_1 \dots R_n w_n$ where
$w_0,\dots,w_n$ are words of terminals and $R_1,\dots,R_n$ are nonterminal symbols, thus yielding 
\[R \to w_0 R_1 w_1 \dots R_n w_n\]
as the general form of a production rule.
We can use this simple observation that production rules factor into nonterminal and terminal components to capture derivations in context-free grammars by functors $p : \FreeOper{\sS}\to\Words{\Sigma}$, as illustrated in Figure~\ref{fig:example-cfg}, and more generally we can use it to define context-free grammars over arbitrary categories.

\begin{defi}\label{definition/cfg-over-cat}
  A \defin{context-free grammar over a category} (or \emph{categorical CFG}) is a tuple $G = (\cC,\sS,S,p)$ consisting of a category $\cC$, a finite species $\sS$ equipped with a distinguished color $S \in \sS$, and a
  functor of operads $p : \FreeOper{\sS} \to \Words{\cC}$.
  We then refer to the colors of~$\sS$ as \defin{nonterminals} and to the operations of~$\FreeOper{\sS}$ as \defin{derivations} or \defin{parse trees}.
  The \defin{context-free language of arrows} $\Lang{G}$ generated by the grammar $G$ is the subset of
  arrows in $\cC$ which, seen as constants of $\Words{\cC}$, are in the image of a closed derivation of $S$, that is, $\Lang{G} = \set{p(\alpha) \mid \alpha : S} \subseteq \cC(A,B)$, where $p(S) = (A,B)$.
\end{defi}
\noindent
As suggested in the Introduction, every context-free grammar in the classical sense $G = (\Sigma,\Nonterminals,\Sentence,\Productions)$ corresponds to a context-free grammar over $\FreeCat{\Bouquet[\Sigma]}$.
For the grammar of Figure~\ref{fig:example-cfg}, the corresponding species $\sS$ has three colors $\mathsf{NP},\mathsf{VP},\mathsf{S}$ and four nodes,
  \begin{align*}
    x_1 &: \mathsf{NP},\mathsf{VP} \to \mathsf{S} \\
    x_2,x_3 &: \mathsf{NP} \\
    x_4 &: \mathsf{NP} \to \mathsf{VP}
  \end{align*}
with the functor $p : \FreeOper{\sS} \to \Words{\Sigma}$ $(= \Words{\FreeCat{\Bouquet[\Sigma]}})$ uniquely defined by the action
\begin{align*}
  x_1 &\mapsto \varepsilon{-}\text{\textvisiblespace}{-}\varepsilon \\
  x_2 &\mapsto \mathsf{mom} \\
  x_3 &\mapsto \mathsf{tom} \\
  x_4 &\mapsto \mathsf{loves}\text{\textvisiblespace}{-}\varepsilon
\end{align*}
  on the generators as displayed in the middle of the figure.
Conversely, any finite species $\sS$ equipped with a color $\Sentence\in\sS$ and a functor of operads $p : \FreeOper{\sS} \to \Words{\Sigma}$ uniquely determines a context-free grammar over the alphabet $\Sigma$.
Indeed, the colors of $\sS$ give the nonterminals of the grammar and $\Sentence$ the distinguished start symbol, while the nodes of $\sS$ together with the functor $p$ give the production rules, with each node $x : R_1,\dots,R_n \to R$ such that $p(x) = w_0{-}\dots{-}w_n$ determining a context-free production rule $x : R \to w_0 R_1 w_1 \dots R_n w_n$.
\begin{prop}
  A language $\Lang{} \subseteq \Sigma^*$ is context-free in the classical sense if and only if it is the language of arrows of a context-free grammar over $\FreeCat{\Bouquet[\Sigma]}$.
\end{prop}


\noindent
An interesting feature of the general notion of categorical context-free grammar is that the nonterminals of the grammar $G = (\cC,\sS,\Sentence,p)$ are \emph{sorted} in the sense that every color of $\sS$ is mapped by $p$ to a unique color of $\Words{\cC}$, corresponding to a pair of objects of $\cC$.
Adapting the conventions from our work on type refinement systems, we sometimes write $R \refs[p] (A,B)$ or simply $R \refs (A,B)$ to indicate that $p(R) = (A,B)$ and say that the nonterminal $R$ refines the ``gap type'' $(A,B)$.
We emphasize that the language $\Lang{G}$ generated by a grammar with start symbol $S \refs (A,B)$ is a subset of the homset $\cC(A,B)$.

\begin{exa}\label{example/cfg-with-end-of-input}
  To illustrate some of the versatility afforded by the more general notion of context-free grammar, consider a CFG over the category $\FreeCat{\Bracket{\Bouquet[\Sigma]}}$ from Example~\ref{example/spliced-ordinal-sum}.
  Such a grammar may include production rules that can only be applied upon reaching the end of the input or only at the beginning, which is useful in practice, albeit usually modelled in an ad hoc fashion.
  For example, the grammar of arithmetic expressions defined by Knuth in the original paper on LR parsing \cite[example (27)]{Knuth1965} may be naturally described as a grammar over $\FreeCat{\Bracket{\Bouquet[\Sigma]}}$, which in addition to having three ``classical'' nonterminals $E, T, P \refs (\ast,\ast)$ contains a distinguished nonterminal $S \refs (\ast,\eowOb)$.
  Knuth's ``zeroth'' production rule $0 : S \to E\$$ is then just a unary node $0 : E \to S$ in $\sS$, mapped by $p$ to the operation $\varepsilon{-}\$ : (\ast,\ast) \to (\ast,\eowOb)$ in $\Words{\FreeCat{\Bracket{\Bouquet[\Sigma]}}}$.
\end{exa}
\noindent
More significant examples of context-free languages of arrows over categories with more than one object will be given in \S\ref{section/csrep}, including context-free grammars over the runs of finite-state automata.

\begin{rem}\label{remark/pointing}
It is at times useful to consider a categorical CFG equivalently as a triple of a \emph{pointed} finite species $(\sS,S)$, a \emph{bipointed} category $(\cC,A,B)$, and a functor of pointed operads $p_G : (\FreeOper{\sS},S) \to (\Words{\cC},(A,B))$.
Note that the operad of spliced arrows construction lifts to a functor
\[ \Wordsbulletonly : \Cat_{\bullet,\bullet} \to \Operads{\bullet} \]
sending a category $\cC$ equipped with a pair of objects $A$ and $B$ to the operad $\Words{\cC}$ equipped with the color $(A,B)$,
and that the free / forgetful adjunction \eqref{equation/adjunction-species-operad} likewise lifts to an adjunction
\begin{equation*}\label{equation/adjunction-species-operad-pointed}
\begin{tikzcd}
\Species{\bullet}
\arrow[rr,"{\FreeOperFunctor_{\bullet}}",yshift=1.5ex]
& \bot &
\Operads{\bullet}\arrow[ll,"{\ForgetOperFunctor_{\bullet}}",yshift=-1.5ex]
\end{tikzcd}
\end{equation*}
between pointed species and pointed operads.
This permits another way of understanding the language of arrows generated by a grammar considered as a functor of pointed operads: since the set of constants of color $(A,B)$ in $\Words{\cC}$ is in bijection with the set of arrows $A \to B$ in $\cC$, we have a natural isomorphism $\elts \circ\, (\Wordsbulletonly) \cong \hom$ for the evident functors $\hom : \Cat_{\bullet,\bullet} \to \Set$ and $\elts : \Operads{\bullet} \to \Set$, and the language $\Lang{G}$ is precisely the image of the function $\elts(p_G) : \elts(\FreeOper{\sS},S) \to \elts(\Words{\cC},(A,B)) \cong \cC(A,B)$.\footnote{It is also possible to take a more abstract and more general perspective on $\Lang{G}$, as building an \emph{initial model} of the grammar $p_G : \FreeOper{\sS} \to \Words{\cC}$ in some target universe $q : \coper{E} \to \coper{B}$ with sufficient fibrational structure. The classical case is recovered by taking $q$ to be the subset fibration $\Subset \to \Set$. See \cite{MZcflsmods}.}
\end{rem}

\subsection{Properties of a context-free grammar and its associated language}
\label{section/properties-of-cfg-and-language}

Standard properties of context-free grammars (cf.~\cite[Ch.~4]{SippuSoisoiParsingTheory1}), considered as categorical CFGs $G = (\FreeCat{\Bouquet[\Sigma]},\sS,\Sentence,p)$, may be reformulated as properties of either the species $\sS$, the operad $\FreeOper{\sS}$, or the functor $p : \FreeOper{\sS} \to \Words{\Sigma}$, with varying degrees of naturality:

\begin{itemize}
\item $G$ is \emph{linear} just in case $\sS$ only has nodes of arity $\le 1$.
  It is \emph{left-linear} (respectively, \emph{right-linear}) just in case it is linear and every unary node $x$ of $\sS$ is mapped by $p$ to an operation of the form $\varepsilon{-}w$ (resp.~$p(x) = w{-}\varepsilon$). 
\item $G$ is in \emph{Chomsky normal form} if $\sS$ only has nodes of arity 2 or 0, the color $\Sentence$ does not appear as the input of any node, every binary node is mapped by $p$ to $\varepsilon{-}\varepsilon{-}\varepsilon$ in $\Words{\FreeCat{\Bouquet[\Sigma]}}$, and every nullary node is mapped to a letter $a \in \Sigma$, unless $R = \Sentence$ in which case it is possible that $p(x) = \varepsilon$.
  (These conditions can be made a bit more more natural by considering $G$ as a context-free grammar over $\FreeCat{\Bracket{\Bouquet[\Sigma]}}$ with $S \refs (\ast,\eowOb)$, see~Example~\ref{example/cfg-with-end-of-input} above.)
\item
  $G$ is \emph{bilinear} (a generalization of Chomsky normal form \cite{LangeLeiss09,Leermakers1989}) iff $\sS$ only has nodes of arity $\le 2$.
\item
  $G$ is \emph{unambiguous} iff for any pair of constants $\alpha, \beta : S$ in $\FreeOper{\sS}$, if $p(\alpha) = p(\beta)$ then $\alpha = \beta$.
  Note that if $p$ is faithful then $G$ is unambiguous, although faithfulness is a stronger condition in general.
\item
  A nonterminal $R$ of $G$ is \emph{nullable} if there exists a constant $\alpha : R$ of $\FreeOper{\sS}$ such that $p(\alpha) = \varepsilon$.
\item
  A nonterminal $R$ of $G$ is \emph{useful} if there exists a pair of a constant $\alpha : R$ and a unary operation $\beta : R \to S$.
  Note that if $G$ has no useless nonterminals then $G$ is unambiguous iff $p$ is faithful.
\end{itemize}
Observe that almost all of these properties 
can be immediately translated to express properties of context-free grammars over \emph{any} category $\cC$. 
Basic closure properties \cite[\S7.3]{HMU2007} of classical context-free languages (CFLs) also generalize easily to CFLs of arrows.
\begin{prop}\label{proposition/basic-closure}\ 
  \begin{enumerate}
  \item If $\Lang{1},\dots,\Lang{k} \subseteq\cC(A,B)$ are CFLs of arrows, so is their union $\bigcup_{i=1}^k \Lang{i} \subseteq\cC(A,B)$.
  \item If $\Lang{1}\subseteq\cC(A_1,B_1), \dots,\Lang{n}\subseteq \cC(A_n,B_n)$ are CFLs of arrows, and if $w_0{-}\dots{-}w_n : (A_1,B_1),\dots,(A_n,B_n) \to (A,B)$ is an operation of $\Words{\cC}$, then the ``spliced concatenation'' $w_0 \Lang{1} w_1 \dots \Lang{n} w_n = \set{w_0u_1w_1\dots u_n w_n \mid u_1 \in \Lang{1},\dots, u_n \in \Lang{n}} \subseteq \cC(A,B)$ is also context-free.
  \item If $\Lang{} \subseteq \cC(A,B)$ is a CFL of arrows in a category $\cC$ and $F : \cC \to \cD$ is a functor of categories, then the functorial image $F(\Lang{}) \subseteq \cD(F(A),F(B))$ is also context-free.
\end{enumerate}
\end{prop}
\begin{proof}
  The proofs of (1) and (2) are just refinements of the standard proofs for context-free languages of words, keeping track of the underlying gap types:
\begin{enumerate}
\item
  Given grammars $G_i = (\cC, \sS_i, S_i, p_i)$ for $i=1\dots k$, where all of the $S_i$ refine the same gap type $(A,B)$, we define a new grammar $G = (\cC, \sS, S, p)$ that generates the union of the languages $\Lang{G} = \bigcup_{i=1}^k \Lang{G_i}$ by taking $\sS$ to be the disjoint union of the colors and nodes of $\sS_1,\dots,\sS_k$ combined with a distinguished color $S$ and unary nodes $x_i : S_i \to S$ for $i = 1\dots k$, and defining $\phi : \sS \to \ForgetOper{\Words{\cC}}$ to be the cotupling of $\phi_1,\dots,\phi_k$ extended by the mappings $S \mapsto (A,B)$ and $x_i \mapsto \id[A]{-}\id[B]$ for $i=1\dots k$.
\item
  Given grammars $G_1 = (\cC,\sS_1,S_1,p_1)$, \dots, $G_n = (\cC,\sS_n,S_n,p_n)$ where $S_i \refs (A_i,B_i)$ for each $1 \le i \le n$, together with an operation
  $w_0{-}\dots{-}w_n : (A_1,B_1),\dots,(A_n,B_n) \to (A,B)$ of $\Words{\cC}$, we construct a new grammar $G = (\cC,\sS,S,p)$ that generates the spliced concatenation $w_0 \Lang{G_1} w_1 \dots \Lang{G_n} w_n$ by taking $\sS$ to be the disjoint union of the colors and nodes of $\sS_1,\dots,\sS_n$ combined with a distinguished color $S$ and a single $n$-ary node $x : S_1,\dots,S_n \to S$, and defining $\phi : \sS \to \ForgetOper{\Words{\cC}}$ to be the cotupling of $\phi_1,\dots,\phi_n$ extended by the mappings $S \mapsto (A,B)$ and $x \mapsto w_0{-}\dots{-}w_n$.
\end{enumerate}
For (3), suppose given a grammar $G=(\cC,\sS,S,p)$ and a functor of categories $F : \cC \to \cD$.
  Then the grammar $F(G)$ generating the language $F(\Lang{G})$ is defined by postcomposing $p$ with $\Words{F} : \Words{\cC} \to \Words{\cD}$ while keeping the species $\sS$ and start symbol $S$ the same, $F(G)= (\cD,\sS,S,(\Words{F})\circ p)$. \qedhere

\end{proof}
\noindent
%
Finally, we note that there is a natural definition of morphisms between context-free grammars over the same category, which we refer to as \emph{translations}.
\begin{defi}\label{definition/translation}
  Let $G_1=(\cC,\sS_1,S_1,p_1)$ and~$G_2=(\cC,\sS_2,S_2,p_2)$ be two grammars over the same category.
  A \defin{translation} $T : G_1 \to G_2$ is given by a functor of operads $T : \FreeOper{\sS_1}\to\FreeOper{\sS_2}$ that commutes with the projection functors
  \[
    \begin{tikzcd}
      \FreeOper{\sS_1} \arrow[dr,"p_1"']\arrow[rr,"T"]  & &\FreeOper{\sS_2}\arrow[dl,"p_2"] \\
      & \Words{\cC} &
    \end{tikzcd}
  \]
  and preserves the start symbol $T(S_1) = S_2$.
  We say that the translation is \defin{full} and/or \defin{faithful} just in case the underlying functor $T : \FreeOper{\sS_1} \to \FreeOper{\sS_2}$ is full and/or faithful.
\end{defi}
\begin{prop}\label{translation-principle}
For any translation $T : G_1 \to G_2$, we have $\Lang{G_1} \subseteq \Lang{G_2}$.
If $T$ is full then $\Lang{G_1} = \Lang{G_2}$, and if $T$ is both full and faithful then the grammars moreover generate isomorphic sets of parse trees for every arrow of the language.
\end{prop}
\begin{proof}
  A translation induces a function from constants $\alpha : S_1$ of $\FreeOper{\sS_1}$ to constants $\beta = T(\alpha) : S_2$ of $\FreeOper{\sS_2}$ lying over the same arrow $p_1(\alpha) = p_2(\beta)$, and hence $\Lang{G_1} \subseteq \Lang{G_2}$.
  If the translation is full (respectively, fully faithful), then this function is surjective (respectively bijective), which implies that $\Lang{G_1} = \Lang{G_2}$, with a bijection between sets of derivations,
\[ \set{\alpha : S_1 \mid p_1(\alpha) = w} \cong \set{\beta : S_2 \mid p_2(\beta) = w} \]
for every arrow $w$, in the case that $T$ is fully faithful.
\end{proof}
\noindent
We refer to the conclusion of Proposition~\ref{translation-principle} under the strongest condition that the functor $T$ is fully faithful as the \defin{strong translation principle}.

\subsection{A fibrational view of parsing as a lifting problem}
\label{section/benabou}

We have seen how any context-free grammar $G=(\cC,\sS,S,p)$ gives rise to a language 
$\Lang{G} = \set{p(\alpha) \mid \alpha : S}$, corresponding to the arrows of $\cC$ which, seen as constants of $\Words{\cC}$, are in the image of some constant of color $S$ of the free operad $\FreeOper{\sS}$.
However, beyond characterizing the language defined by a grammar, in practice one is often confronted with a dual problem, namely that of parsing: given a word $w$, we want to compute the set of all its parse trees, or at least determine all of the nonterminals which derive it.
In our functorial formulation of context-free derivations, this amounts to computing the inverse image of $w$ along the functor $p$, i.e., the set of constants $p^{-1}(w) = \set{\alpha \mid p(\alpha) = w}$, or alternatively the set of colors in the image of $p^{-1}(w)$ along the output-color function.

To better understand this view of parsing as a lifting problem along a functor of operads, we find it helpful to first recall the correspondence between functors of categories $p : \ccat{D} \to \ccat{C}$ and lax functors $F : \ccat{C} \to \Span\Set$, where $\Span\Set$ is the bicategory whose objects are sets, whose 1-cells $S:X\,{\spanmap} Y$ are spans $X\leftarrow S\rightarrow Y$,
and whose 2-cells are morphisms of spans.
Suppose given such a functor $p : \ccat{D} \to \ccat{C}$.
To every object $A$ of $\ccat{C}$ there is an associated fiber $F(A) = p^{-1}(A)$ of objects in $\ccat{D}$ living over $A$, while to every arrow $w : A \to B$ of $\ccat{C}$ there is an associated fiber $F(w) = p^{-1}(w)$ of arrows in $\ccat{D}$ living over $w$, equipped with a pair of projection functions $F(A)\leftarrow F(w)\rightarrow F(B)$ mapping any lifting $\alpha : R \to S$ of $w : A \to B$ to its source $R \in F(A)$ and target $S \in F(B)$.
Moreover, given a pair of composable arrows
$u:A\to B$ and $v:B\to C$ in $\ccat{C}$, there is a morphism of spans
\begin{equation}\label{equation/lax-o}
\begin{tikzcd}[column sep = 1.2em, row sep = .8em]
{F(u)F(v)}
\arrow[double,-implies,rr]
&&
{F(u v)}
\quad
:
\quad
{F(A)} \arrow[spanmap,rr]
&&
{F(C)}
\end{tikzcd}
\end{equation}
from the composite of the spans 
$F({u}):F({A})\spanmap F({B})$
and 
$F({v}):F({B})\spanmap F({C})$ associated
to $u:A\to B$ and $v:B\to C$ 
to the span $F({u v}):F({A})\spanmap F({C})$ associated to the composite arrow $u v:A\to C$.
This morphism of spans is realized using composition in the category $\ccat{D}$,
namely by the function taking any pair of a lifting $\alpha : R \to S$ of $u$ and a lifting $\beta : S \to T$ of $v$
to the composite $\alpha \beta : R \to T$, which is a lifting of $u v$ by functoriality $p(\alpha\beta) = p(\alpha)p(\beta) = uv$.
Similarly, the identity arrows in the category $\ccat{D}$ define, for every object $A$ of the category $\ccat{C}$,
a morphism of spans
\begin{equation}\label{equation/lax-I}
\begin{tikzcd}[column sep = 1.2em, row sep = .8em]
{\id[F(A)]}
\arrow[double,-implies,rr]
&&
{F({\id[A]})}
\quad
:
\quad
{F({A})} \arrow[spanmap,rr]
&&
{F({A})}
\end{tikzcd}
\end{equation}
from the identity span $F(A)\leftarrow F(A)\rightarrow F(A)$
to the span associated to the identity arrow $\id[A]:A\to A$.
Associativity and neutrality of composition in $\ccat{D}$
ensure that the 2-cells \eqref{equation/lax-o} and \eqref{equation/lax-I} make the diagrams below commute:
\[
\begin{tikzcd}[column sep = .8em, row sep = 1.8em]
{F(u) F(v) F(w)}
\arrow[rr,double,-implies]
\arrow[dd,double,-implies]
&&
{F(u) F({v w})}
\arrow[dd,double,-implies]
\\
\\
{F({u v})F(w)}
\arrow[rr,double,-implies]
&&
{F({u v w})}
\end{tikzcd}
\quad
\begin{tikzcd}[column sep = .6em, row sep = .8em]
&& {F({u})}
\arrow[dddd,double,-]
\arrow[lldd,double,-implies]
&&
\\
\\
{F({\id[A]})F({u})}
\arrow[rrdd,double,-implies]
\\
\\
&& {F({u})}
\end{tikzcd}
\quad
\begin{tikzcd}[column sep = .6em, row sep = .8em]
{F({u})}
\arrow[dddd,double,-]
\arrow[rrdd,double,-implies]
&&
\\
\\
&&
{F({u}) F({\id[B]})}
\arrow[lldd,double,-implies]
\\
\\
{F({u})}
\end{tikzcd}
\]
for all triples of composable arrows $u:A\to B$, $v:B\to C$ and $w:C\to D$, and therefore that 
this collection of data defines what is called a lax functor $F:\ccat{C}\to\Span{\Set}$.
In general it is \emph{only} lax, in the sense that the 2-cells $F(u) F(v) \Rightarrow F({uv})$ and $\id[F(A)] \Rightarrow F({id_A})$ are not necessarily invertible.

Conversely, starting from the data provided by a lax functor $F:\ccat{C}\to\Span{\Set}$, we can define a category denoted $\els[\cC] F$ together with a functor $\pi:\els[\cC] F\to\ccat{C}$.
The category $\els[\cC] F$ has objects the pairs $(A,R)$ of an object~$A$ in $\ccat{C}$ and an element $R\in F(A)$,
and arrows $(w,\alpha):(A,R)\to (B,S)$ the pairs of an arrow $w:A\to B$ in $\ccat{C}$ and an element $\alpha\in F(w)$ mapped to $R\in F(A)$ and $S\in F(B)$ by the respective legs of the span:
\[
\begin{tikzcd}
& 1\arrow[ld,"R"']\arrow[rd,"S"]\arrow[d,"\alpha"] & \\
F(A)  &
F(w)\arrow[l]\arrow[r] & 
F(B)
\end{tikzcd}
\]
Here we slightly abuse notation by writing $F(w)$ for the apex of the span $F(A) \spanmap F(B)$ assigned to $w : A \to B$.
The composition and identity of the category $\els[\cC] F$ are then given by the morphisms of spans $F(u) F(v)\Rightarrow F({u v})$ and $\id[F(A)] \Rightarrow F({\id[A]})$ witnessing the lax functoriality of $F:\ccat{C}\to\Span{\Set}$.
The functor $\pi:\els[\cC] F\to\ccat{C}$ is given by the first projection.
This construction of a category $\els[\cC] F$ equipped with a functor $\pi:\els[\cC] F\to\ccat{C}$ starting from a lax functor $F : \ccat{C} \to \Span\Set$, apparently first made explicit in \cite[Prop.~4]{PavlovicAbramsky97}, is a mild variation of B\'enabou's construction \cite{Benabou1972,Benabou2000} of the same starting from a lax normal functor $F : \ccat{C}^\op \to \Dist$, which is itself a generalization of the well-known Grothendieck construction of a fibration starting from a pseudofunctor $F : \ccat{C}^\op \to \Cat$.
One can show that given a functor of categories 
$p:\ccat{D}\to\ccat{C}$, the construction
applied to the associated lax functor $F:\ccat{C}\to\Span{\Set}$
induces a category $\els[\cC] F$ isomorphic to $\ccat{D}$,
in such a way that~$p$ coincides with the isomorphism composed with $\pi$.
Recently, Ahrens and Lumsdaine \cite{AhrensLumsdaine2019} have introduced the useful terminology ``displayed category'' to refer to this way of presenting a category $\ccat{D}$ equipped with a functor $\ccat{D} \to \ccat{C}$ as a lax functor $\ccat{C} \to \Span\Set$, with their motivations coming from formalization of mathematics.

The constructions which turn a functor of categories $p:\ccat{D}\to\ccat{C}$ into a lax functor $F:\ccat{C}\to\Span{\Set}$ and back into a functor $\pi : \els[\cC] F \to \ccat{C}$ can be adapted smoothly to functors of operads, viewing $\Span\Set$ as a 2-categorical operad whose $n$-ary operations
$S:X_1,\dots,X_n\,{\spanmap} Y$ are multi-legged-spans
\[
\begin{tikzcd}[row sep = 0em]
X_1 \\
\vdots &
S
\arrow[lu]\arrow[ld]\arrow[r]
&
Y \\
X_n
\end{tikzcd}
\]
or equivalently spans $X_1\times\dots\times X_n \leftarrow S\rightarrow Y$,
and with the same notion of 2-cell.
We will follow Ahrens and Lumsdaine's suggestion and refer to the data of such a lax functor $F : \oO \to \Span\Set$ representing an operad $\oD \cong \els[\oO] F$ equipped with a functor $p : \oD \to \oO$ as a \defin{displayed operad}.

\subsection{An inductive formula and a sequent calculus for displayed free operads}
\label{section/magic-formula}

It is folklore that the free operad over a species $\sS = (C,V,i,o)$ may be described concretely as a certain family of trees: operations of $\FreeOper{\sS}$ are interpreted as (rooted planar) trees whose edges are colored by the elements of $C$ and whose nodes are labelled by the elements of $V$, subject to the constraints imposed by the functions $i : V \to C^*$ and $o : V \to C$.
The formal construction of the free operad may be viewed as a free monoid construction, adapted to a situation where the ambient monoidal product (in this case, the composition product of species) is only distributive on the left, see \cite[II.1.9]{MarklSchniderStasheff} and \cite[Appendix~B]{BJT1997}.

From the perspective of programming semantics, it is natural to consider the underlying species of $\FreeOper{\sS}$ as an inductive data type, corresponding to the initial algebra for the endofunctor $W_\sS$ on $C$-colored species defined by
\[
W_\sS = \Rspecies \mapsto \speciesunit{} + \sS \circ \Rspecies
\]
where $+$ denotes the coproduct of $C$-colored species which is constructed by taking the disjoint union of nodes, while $\circ$ and $\speciesunit{}$ denote respectively the composition product of $C$-colored species and the identity species, defined as follows.
Given two $C$-colored species $\sS$ and $\Rspecies$, 
the $n$-ary nodes $R_1,\dots,R_n\to S$
of $\sS\circ \Rspecies$ 
are formal composites $g\bullet (f_1,\dots,f_k)$ consisting of a node
$g:S_1,\dots,S_k\to S$ of $\sS$ and of a tuple of nodes
$f_1:\Gamma_1\to S_1$, $\dots$, $f_k:\Gamma_k\to S_k$ of $\Rspecies$,
such that the concatenation of the lists of colors $\Gamma_1,\dots,\Gamma_k$ is equal to the list $R_1,\dots,R_n$.
The unit $\speciesunit{}$ is the $C$-colored species with a single unary node $\ast_R:R\to R$ for every color $R\in C$,
and no other nodes.

As the initial $W_\sS$-algebra, the free operad over $\sS$ is equipped with a map of species
$\speciesunit{} + \sS \circ \FreeOper{\sS} \longrightarrow \FreeOper{\sS}$,
which by the Lambek Lemma is invertible,
with the following interpretation: any operation of $\FreeOper{\sS}$ is either an identity operation, or the parallel composition of a node of $\sS$ with a list of operations of $\FreeOper{\sS}$.
Note that this interpretation also corresponds to a canonical decomposition
\begin{equation}
  \includegraphics[scale=1]{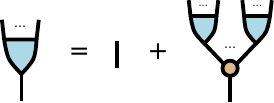}
\end{equation}
of rooted planar trees, by iterated examination of the root.
An operation of $\FreeOper{\sS}$ may be seen as a rooted planar tree whose edges and nodes are labelled by the species $\sS$, and we sometimes refer to such an operation as an \emph{$\sS$-rooted tree} (cf.~\cite[\S3.2]{BeLaLe1998Species}).

It is possible to derive an analogous inductive characterization of \emph{functors} $p : \FreeOper{\sS} \to \oO$ from a free operad into an arbitrary operad $\oO$ considered as displayed free operads, i.e., as lax functors $F : \oO \to \Span\Set$ generated by an underlying map of species $\phi : \sS \to \ForgetOper{\oO}$.
Two subtleties arise.
First, that the species $\sS$ and the operad $\oO$ may in general have a different set of colors, related by the change-of-color function $\phi_C$.
To account for this, rather than restricting the operations $+, \circ,\speciesunit{}$ to the category of $C$-colored species,
one should consider them as global functors
\[ +,\circ : \Species{} \times_{\Set} \Species{} \to \Species{}
  \qquad
  \speciesunit{} : \Set \to \Species{}
\]
on the ``polychromatic'' category of species, which respect the underlying sets of colors in a functorial way.
Second, and more significantly, the above functor $W_\sS$ transports a species $\Rspecies$ living over $\ForgetOper\oO$ to a species living over $\speciesunit{} + \ForgetOper\oO \circ \ForgetOper\oO$, so that in order to obtain again a species living over $\ForgetOper\oO$ (and thus define an endofunctor) one needs to ``push forward'' along the canonical $W_{\ForgetOper\oO}$-algebra $[e,m] : \speciesunit{} + \ForgetOper\oO \circ \ForgetOper\oO \longrightarrow \ForgetOper\oO$ that encodes the operad structure of $\oO$, seen as a monoid in $(\Species{},\circ,\speciesunit{})$.
A detailed discussion is beyond the scope of this paper, but we nevertheless state the following proposition, whose intuitive content should be clear:
\begin{prop}\label{prop/magic-formula} Let $\phi : \sS \to \ForgetOper{\oO}$ be a map of species from a species $\sS$ into an operad $\oO$, and let $p : \FreeOper{\sS} \to \oO$ be the corresponding functor from the free operad.
  Then the associated lax functor $F : \oO \to \Span\Set$ computing the fibers of $p$ is given by $F(A) = \phi^{-1}(A)$ on colors of $\oO$, and by the least family of spans $F(f)$ indexed by operations $f : A_1,\dots,A_n \to A$ of $\oO$ such that
\begin{equation}\label{equation/magic-formula}
  F(f) \cong \sum_{{f=\id[A]}\atop {\phi(R) = A}}\id[R]\ \ + \sum_{f=g\circ (h_1,\dots,h_k)} \phi^{-1}(g)\bullet(F({h_1}),\dots,F({h_k}))
\end{equation}
  where we write $\circ$ for composition in the operad $\oO$ and $\bullet$ for formal composition of nodes in $\sS$ with operations in $\FreeOper{\sS}$.
  Specializing the formula to constant operations, the left summand disappears and \eqref{equation/magic-formula} simplifies to:
\begin{equation}\label{equation/magic-formula0}
  F(c) \cong \sum_{c=g\circ (c_1,\dots,c_k)} \phi^{-1}(g)\bullet(F({c_1}),\dots,F({c_k}))
\end{equation}
\end{prop}
\noindent
For the benefit of readers with a background in logic, we note that the inductive characterization of displayed free operads may be equivalently expressed as a sequent calculus parameterized by an underlying map of species $\phi : \sS \to \oO$, see Figure~\ref{fig:sequent-calculus}.
Adapting a convention from our type refinement systems work, we write
\[ R_1,\dots,R_k \seq[\phi]{g} R \]
for the judgment whose evidence is provided by a node $x : R_1,\dots,R_k \to R$ in $\sS$ such that $\phi(x) = g$, and similarly 
\[ S_1,\dots,S_n \seq[p]{f} R \]
for the judgment whose evidence is provided by an operation $f : S_1,\dots,S_n \to R$ in $\FreeOper{\sS}$ such that $p(\alpha) = f$.
Again, we write $R \refs[p] A$ to indicate that $R$ is a color of $\FreeOper{\sS}$ with image $p(R) = A$ in $\oO$.
As standard, we treat the inference rules as an inductive definition, in the sense that the set of derivations of the calculus is the least set closed under the inference rules.
Then formulas~\eqref{equation/magic-formula} and \eqref{equation/magic-formula0} are equivalent to the following:
\begin{prop}\label{prop/sequent-calculus}
  For the sequent calculus of Fig.~\ref{fig:sequent-calculus}, parameterized with respect to a map of species $\phi : \sS \to \oO$ and the corresponding functor $p : \FreeOper{\sS} \to \oO$, derivations of $\Gamma \seq[p]{f} R$ are in one-to-one correspondence with operations $\alpha : \Gamma \to R$ in $\FreeOper{\sS}$ such that $p(\alpha) = f$.
  As a special case, constants $\alpha : R$ such that $p(\alpha) = c$ are in one-to-one correspondence with derivations of $\seq[p]{c} R$, which necessarily only use instances of the \textsc{node} rule.
\end{prop}
\begin{figure}
  \[
    \inferrule*[vcenter,Left={leaf}]{R \refs[p] A}{R \seq[p]{\id[A]} R}
    \qquad\qquad
    \inferrule*[vcenter,Left={node}]{f = g\circ (h_1,\dots,h_k) \\\\ R_1,\dots,R_k \seq[\phi]{g} R \\ \Gamma_1 \seq[p]{h_1} R_1 \quad \dots \quad \Gamma_k \seq[p]{h_k} R_k}{\Gamma_1,\dots,\Gamma_k \seq[p]{f} R}
  \]
  \caption{A sequent calculus for displayed free operads.}
  \label{fig:sequent-calculus}
\end{figure}

\subsection{Application to parsing} 
\label{section/application-to-parsing}

Instantiating \eqref{equation/magic-formula0} with the underlying functor 
of a categorical CFG generated by a map of species $\phi : \sS \to \Words{\ccat{C}}$, we immediately obtain the following characteristic formula for the family of sets of parse trees $F(w)$ of an arrow $w$ in $\ccat{C}$, seen as liftings of the constant $w$ in $\Words{\ccat{C}}$ to a constant in $\FreeOper{\sS}$:
\begin{equation}\label{equation/magic-formula0-words}
  F(w) \cong \sum_{w=w_0 u_1 w_1\dots u_n w_k} \phi^{-1}(w_0{-}\dots{-}w_k)\bullet(F({u_1}),\dots,F({u_k}))
\end{equation}
Equivalently, specializing the sequent calculus of Figure~\ref{fig:sequent-calculus} and translating the generic judgment forms into traditional CFG notation by
\begin{align*}
  R_1,\dots,R_k \seq[\phi]{w_0{-}\dots{-}w_n} R &\quad\leadsto\quad R \to w_0 R_1 w_1\dots R_n w_n \\
   \seq[p]{w} R &\quad\leadsto\quad R \to^+ w
\end{align*}
we obtain the following inference rule:
\begin{equation}\label{equation/magic-formula0-words-sequent-calculus}
  \inferrule{w = w_0u_1 w_1 \dots u_k w_k \\ R \to w_0R_1w_1 \dots R_kw_k \\ R_1 \to^+ u_1 \quad \dots \quad R_k \to^+ u_k}{R \to^+ w}
\end{equation}
inductively defining the set of closed context-free derivations.
For any given word $w$, let $N_w = \set{R \mid R \to^+ w}$ be the set of nonterminals that derive it, corresponding to the image of $F(w)$ along the function returning the root label of a parse tree.
We have that:
\begin{equation}\label{equation/magic-formula0-words-subset}
\mathtoolsset{multlined-width=0.9\displaywidth}
\begin{multlined}
  R \in N_w \iff w=w_0 u_1 w_1\dots u_k w_k \ \wedge \\
  R \to w_0R_1w_1\dots R_k w_k \wedge R_1 \in N_{u_1} \wedge \dots \wedge R_k \in N_{u_k}
\end{multlined}
\end{equation}
This equation is essentially the characteristic formula expressed by Leermakers \cite{Leermakers1989} for the defining relation of the ``C-parser'', which generalizes the well-known Cocke-Younger-Kasami (CYK) algorithm.
Presentations of the CYK algorithm are usually restricted to grammars in Chomsky normal form (cf.~\cite{LangeLeiss09}), but as observed by Leermakers, the relation $N_w$ defined by \eqref{equation/magic-formula0-words-subset} can be solved effectively for any context-free grammar $G$ and given word $w = a_1 \dots a_n$ by building up a parse matrix $N_{i,j}$ indexed by the subwords $w_{i,j} = a_{i+1}\dots a_j$ for all $1 \le i \le j \le n$, and yields a cubic complexity algorithm in the case that $G$ is bilinear in the sense of \S\ref{section/properties-of-cfg-and-language}.
Moreover, by adding nonterminals, it is always possible to transform a CFG into a bilinear CFG that generates the same language, even preserving the original derivations up to isomorphism.
\begin{prop}\label{proposition/bilinear-normal-form}
  For any categorical context-free grammar $G = (\ccat{C}, \sS, S, p)$, there is a bilinear CFG $G_\bin = (\ccat{C}, \sS_\bin, S, p_\bin)$ together with a fully faithful functor of operads
  $B : \FreeOper{\sS} \to \FreeOper{\sS_\bin}$ such that $p = B p_\bin$.
  In particular, by the strong translation principle (Prop.~\ref{translation-principle}) we have that $\Lang{G} = \Lang{G_\bin}$ with isomorphic sets of parse trees.
\end{prop}
\begin{proof}
Given $G = (\ccat{C},\sS,S,p)$, a bilinear grammar $G_\bin = (\ccat{C},\sS_\bin,S,p_\bin)$ is constructed as follows.
$\sS_\bin$ includes all of the colors and all of the nullary nodes of $\sS$, with $\phi_\bin(R) = \phi(R)$ and $\phi_\bin(c) = \phi(c)$.
Additionally, for every node $x : R_1,\dots,R_n \to R$ of $\sS$ of positive arity $n>0$, where $\phi(x) = w_0{-}\dots{-}w_n : (A_1,B_1),\dots,(A_n,B_n) \to (A,B)$ in $\Words{\ccat{C}}$, we include:
\begin{itemize}
\item $n$ new colors $I_{x,0},\dots,I_{x,n-1}$, with $\phi_\bin(I_{x,i-1}) = (A,A_i)$ for $1 \le i \le n$;
\item one nullary node $x_0: I_{x,0}$, with $\phi_\bin(x_0) = w_0$;
\item $n$ binary nodes $x_1,\dots,x_n$, where $x_i : I_{x,i-1},R_i \to I_{x,i}$ and $\phi_\bin(x_i) = \id[A]{-}\id[A_i]{-}w_i$ for all $1 \le i \le n$, under the convention that $I_{x,n} = R$.
\end{itemize}
We define the functor $B : \FreeOper{\sS} \to \FreeOper{\sS_\bin}$ on colors by $B(R) = R$, on nullary nodes by $B(c) = c$, and on nodes $x : R_1,\dots,R_n \to R$ of positive arity by
$B(x) = x_n\circ_0 \dots \circ_0 x_1 \circ_0 x_0$.
By induction on $n$, there is a one-to-one correspondence between nodes $x : R_1,\dots,R_n \to R$ of $\sS$ and operations $B(x) : R_1,\dots,R_n \to R$ of $\FreeOper{\sS_\bin}$, so the functor $B$ is fully faithful.
\end{proof}

\subsection{Generalized context-free grammars}
\label{section/gCFGs}

We conclude this section by observing that for many purposes, a functor of operads from a free operad generated by a finite species $\sS$ into an arbitrary operad $\coper{O}$ may be considered (after selection of a start symbol in $\sS$) as a ``generalized'' context-free grammar, generating a language of constants in $\coper{O}$.
\begin{defi}\label{definition/gcfg}
  A \defin{generalized context-free grammar (over an operad)} is a tuple $G = (\coper{O},\sS,S,p)$ consisting of an operad $\coper{O}$, a finite species $\sS$ equipped with a distinguished color $S \in \sS$, and a functor of operads $p : \FreeOper{\sS} \to \coper{O}$.
  The \defin{generalized context-free language of constants} generated by the grammar $G$ is given by the subset of constants $\Lang{G} = \set{p(\alpha) \mid \alpha : S} \subseteq \oO(A)$, where $p(S) = A$.
\end{defi}
\noindent
Hence a classical context-free word grammar can be seen as a generalized CFG (gCFG) over a spliced word operad $\coper{O} = \Words{\Sigma}$, and a categorical CFG as a gCFG over a spliced arrow operad $\coper{O} = \Words{\Ccategory}$.
By selecting $\coper{O}$ appropriately, Definition~\ref{definition/gcfg} encompasses other extensions of the notion of context-free grammar from the literature.
\begin{exa}\label{example/mcfg}
  For any operad $\coper{P}$, one can construct a new operad $\AffMon{\coper{P}}$ whose colors are lists $[A_1,\dots,A_k]$ of colors in $\coper{P}$ and whose $n$-ary operations $[\Gamma_1],\dots,[\Gamma_n] \to [A_1,\dots,A_k]$ are pairs $([f_1,\dots,f_k],\sigma)$ of a list of operations
  $f_1 : \Omega_1 \to A_1, \dots, f_k : \Omega_k \to A_k$
  in $\coper{P}$ and an injection $\sigma : \Omega_1,\dots,\Omega_k \hookrightarrow \Gamma_1,\dots,\Gamma_n$ from the concatenation of the inputs of the $k$ operations to the concatenation of the $n$ input lists.
  Observe that if $\coper{P}$ is an uncolored operad (i.e., has a unique color, like $\Words{\Sigma}$), then the colors of $\AffMon{\coper{P}}$ are isomorphic to natural numbers.
  Let us also remark that $\AffMon{\coper{P}}$ is in fact the \emph{free semi-cartesian (or ``affine'') strict monoidal operad} over $\coper{P}$.\footnote{See \cite[p.~56--58]{LeinsterHOHC} for discussion of the similar constructions of the free PRO(P) over a colored operad.}

  A gCFG over $\AffMon{\Words{\Sigma}}$ with start symbol $S \refs[] 1$ is precisely a \emph{multiple context-free grammar}~\cite{seki+mcfg1991}, while more generally a gCFG over $\AffMon{\Words{\Ccategory}}$ with $S \refs[] [(A,B)]$ (for some objects $A$ and $B$ of $\Ccategory$) could be called a ``multiple categorical CFG''.
  Such a grammar 
  is a \emph{$k$-multiple CFG} (again in the sense of \cite{seki+mcfg1991}) just in case every nonterminal $R \refs[] [\Gamma]$ refines a list of length $|\Gamma| \le k$.
  
  For instance, let $\Ccategory$ be the free category over the following graph:
  \[\begin{tikzcd}
      A \ar[out=120,in=60,loop,"a"]\ar[r,"\#"] &
      B \ar[out=120,in=60,loop,"b"]\ar[r,"\#'"] &
      C \ar[out=120,in=60,loop,"c"]
    \end{tikzcd}\]
  We can define a 3-mCFG generating the language $a^n\# b^n \#' c^n$ with a pair of nonterminals
  \begin{align*}
    S &\refs[] [(A,C)] \\
    R &\refs[] [(A,A),(B,B),(C,C)]
  \end{align*}
  and production rules specified by the triple of nodes in $\sS$
  \begin{align*}
    x_1 &: R \\
    x_2 &: R \to R \\
    x_3 &: R \to S
  \end{align*}
  mapped respectively to the following operations in $\AffMon{\Words{\Ccategory}}$
  \begin{align*}
    [\id[A],\id[B],\id[C]] &: [(A,A),(B,B),(C,C)] \\
    [a{-}\id[A],b{-}\id[B],c{-}\id[C]] &: [(A,A),(B,B),(C,C)] \to [(A,A),(B,B),(C,C)] \\
    [{-}\#{-}\#'{-}] &: [(A,A),(B,B),(C,C)] \to [(A,C)]
  \end{align*}
  where we have elided the injection component and only shown the lists of operations in $\Words{\Ccategory}$, since the injection $\sigma$ used is always the identity.
\end{exa}
\begin{exa}\label{example/pmcfg}
  The \emph{free cartesian strict monoidal operad} $\CartMon{\coper{P}}$ over an operad $\coper{P}$ is constructed just as $\AffMon{\coper{P}}$, but replacing injections by arbitrary functions $\sigma : \Omega_1,\dots,\Omega_k \to \Gamma_1,\dots,\Gamma_n$.
  A gCFG over $\CartMon{\Words{\Sigma}}$ is precisely a \emph{parallel multiple context-free grammar} in the sense of \cite{seki+mcfg1991}, and more generally a gCFG over $\CartMon{\Words{\Ccategory}}$ could be called a ``parallel multiple categorical CFG''.
\end{exa}
\noindent
Examples of generalized context-free languages of a more semantic flavor can also be derived by taking $\coper{O}$ to be $\Set$ or similar.
Recall that $\Set$ is the (large) operad whose colors are sets and whose operations $f:X_1,\dots,X_n \to X$ are $n$-ary functions $f:X_1\times\cdots\times X_n\to X$, with composition defined as expected.
\begin{exa}\label{example/series-parallel}
  (Adapted from \cite[\S1.1.3]{CourcelleEngelfriet2012}.)
  Let $\DiGr_{\bullet,\bullet}$ be the set of finite directed graphs equipped with two distinct vertices labeled $\src$ and $\tgt$.
  The \emph{parallel composition} $G \parallel H$ of a pair of graphs $(G,H) \in \DiGr_{\bullet,\bullet} \times \DiGr_{\bullet,\bullet}$ is formed by starting from the disjoint union of the two graphs and identifying their two sources and two targets to form a graph $G \parallel H \in \DiGr_{\bullet,\bullet}$ with $\src(G\parallel H) = \src(G) = \src(H)$ and $\tgt(G\parallel H) = \tgt(G) = \tgt(H)$.
  The \emph{series composition} $G \series H$ is formed by similarly starting from the disjoint union of $G$ and $H$ but then identifying $\tgt(G) = \src(H)$ and taking $\src(G\series H) = \src(G)$ and $\tgt(G\series H) = \tgt(H)$.
  We can define a gCFG over $\Set$ with one nonterminal $S \refs \DiGr_{\bullet,\bullet}$ and three production rules corresponding to a pair of binary nodes $par : S,S \to S$ and $seq : S,S \to S$ mapped to the operations $\parallel$ and $\series$ respectively, plus a constant $e : S$ mapped to the unique directed graph $\bullet\to \bullet$ with one edge and two vertices.
  The language of constants generated by this gCFG is exactly the set of \emph{series-parallel graphs}.
\end{exa}

Many of the properties of categorical CFGs discussed in \S\ref{section/properties-of-cfg-and-language} extend directly to gCFGs, as do all of the closure properties of CFLs and the translation principles.
For the record, we state here the appropriate generalizations for gCFLs of the closure properties listed in Proposition~\ref{proposition/basic-closure}, omitting the proofs since they are essentially identical.
\begin{prop}\label{proposition/basic-closure-gCFLs}\ 
  \begin{enumerate}
  \item If $\Lang{1},\dots,\Lang{k} \subseteq\oO(A)$ are gCFLs, then so is their union $\bigcup_{i=1}^k \Lang{i} \subseteq\oO(A)$.
  \item If $\Lang{1}\subseteq\oO(A_1), \dots,\Lang{n}\subseteq \oO(A_n)$ are gCFLs, and if $f : A_1,\dots,A_n \to A$ is an operation of $\oO$, then the application $f(\Lang{1},\dots,\Lang{n}) = \set{f(u_1,\dots,u_n) \mid u_1 \in \Lang{1},\dots, u_n \in \Lang{n}} \subseteq \oO(A)$ is also a gCFL.
  \item If $\Lang{} \subseteq \oO(A)$ is a gCFL of constants in an operad $\oO$ and $F : \oO \to \oP$ is a functor of operads, then the functorial image $F(L) \subseteq \oP(F(A))$ is also a gCFL.
\end{enumerate}
\end{prop}

\subsection{Related work}
\label{section/cfls/related}

Although the notion of context-free grammar over a category seems to be new, the idea of viewing classical context-free grammars fibrationally and representing them as functors appears in various guises in the literature.
Walters briefly described in~\cite{Walters1989} how context-free grammars may be represented by certain morphisms of multigraphs that we would call maps of species (see comment immediately after Definition~\ref{definition/species}) and might notate $\phi : \sS \to \tilde\Sigma$, where $\tilde\Sigma$ is defined similarly to the underlying species of $\Words{\Sigma}$ but with two distinct kinds of nodes: nullary nodes $a : *$ for every letter $a \in \Sigma$, and $n$-ary nodes $\mu_n : *^n \to *$ standing for concatenation.
Note that taking $\tilde\Sigma$ so-defined as the target of $\phi$ imposes a restriction on the grammar (similar to Chomsky normal form), that all production rules are either of the form $R \to a$ or $R \to R_1\dots R_n$.
In order to recover the language generated by a grammar with start symbol $S$, Walters next considers the functor $\mathcal{F}_\times\phi : \mathcal{F}_\times \sS \to \mathcal{F}_\times \tilde{\Sigma}$ between the \emph{free categories with products} generated by the species, then composes with the canonical functor $\psi : \mathcal{F}_\times\tilde{\Sigma} \to \textrm{Mon}_\Sigma$ into the Lawvere theory of monoids containing $\Sigma$, and finally takes the image $\psi (\mathcal{F}_\times\phi(\mathcal{F}_\times\sS(1,S)))$.
All in all, we find the encoding of classical CFGs by functors of operads $\FreeOper{\sS} \to \Words{\Sigma}$ to be simpler, more direct and more general, but Walters' approach is certainly similar in spirit.

The methodology of \emph{abstract categorial grammars} introduced by De~Groote~\cite{deGroote2001acgs} is probably even more closely related to our approach.
Based on linear lambda-calculus, the approach consists of expressing different kinds of grammars by specifying a signature $\Sigma_1$ as an ``abstract vocabulary'', a signature $\Sigma_2$ as an ``object vocabulary'', and a ``lexicon'' given by a morphism of signatures $\Sigma_1 \to \Sigma_2$.
In particular, context-free grammars are represented by taking the abstract vocabulary $\Sigma_1$ to have a type for each nonterminal and a constant $x : R_1 \multimap \dots \multimap R_n \multimap R$ for every production rule $x : R \to w_0R_1\dots R_n w_n$, the object vocabulary $\Sigma_2$ to have a single type $*$ and a constant $a : * \multimap *$ for every letter, and the lexicon to map every nonterminal $R$ to $*\multimap *$, and every constant $x$ as above to the lambda-term $\lambda u_1\dots u_n.w_0 \circ u_1 \circ \dots \circ u_n \circ w_n$.
Although the definition of ACGs was originally given type-theoretically rather than categorically, the framework could very clearly be reformulated naturally in an operadic setting, using free closed operads, and seems to have deep connections with the approach described here.
One way of making this connection more precise could be to prove that the interpretation of pairs $(A,B)$ as types $A \multimap B$ extends to a full and faithful embedding $\Words{\cC} \to \FreeOper[\multimap]{\cC}$ of the operad of spliced arrows $\Words{\cC}$ into the free closed operad over $\cC$, where the operations of $\FreeOper[\multimap]{\cC}$ may be represented by $\beta\eta$-normal (ordered) linear lambda terms.
It is worth emphasizing that ACGs were introduced to capture a wide variety of grammatical formalisms, including categorial grammars in the style of Lambek calculus but also extensions of context-free grammars like we discussed in \S\ref{section/gCFGs} such as multiple~CFGs \cite{deGrootePogodalla2004}.


\section{Finite-state automata over categories and operads}
\label{section/nfas}

We have seen how the classical notion of context-free grammar generalizes naturally to define context-free grammars over any category $\cC$, or even over any operad $\oO$, presented by a functor of operads $\FreeOper{\sS} \to \Words{\cC}$ or $\FreeOper{\sS} \to \oO$ together with a distinguished start symbol in $\Sspecies$.
In this section, we carry on in the same fibrational vein and begin by explaining how classical nondeterministic finite-state word automata may be generalized to define NFAs over any category, presented by functors of categories $\cQ \to \cC$ that are both \emph{finitary} and have the \emph{unique lifting of factorizations} (ULF) property, together with distinguished initial and accepting states in $\cQ$.
The base category $\cC$ is often freely generated, as in the case of word automata, but we also give natural examples of automata over non-free categories.
One important aspect of this approach is that it adapts smoothly when one shifts from word automata to tree automata, simply by replacing categories with operads.
Indeed, a key fact is that the splicing construction transports any finitary ULF functor of categories $p:\cQ\to\cC$ to a finitary ULF functor of operads $\Words{p}:\Words{\cQ}\to\Words{\cC}$, thereby transforming an NFA over a category into an NFA over its operad of spliced arrows.
As we will see in the following section (\S\ref{section/csrep}), this observation greatly facilitates the proof
that context-free languages are closed under intersection with regular languages.

\subsection{Nondeterministic word automata as finitary ULF functors over categories}\label{section/warmup}
\begin{figure}
  \begin{center}\includegraphics[width=\textwidth]{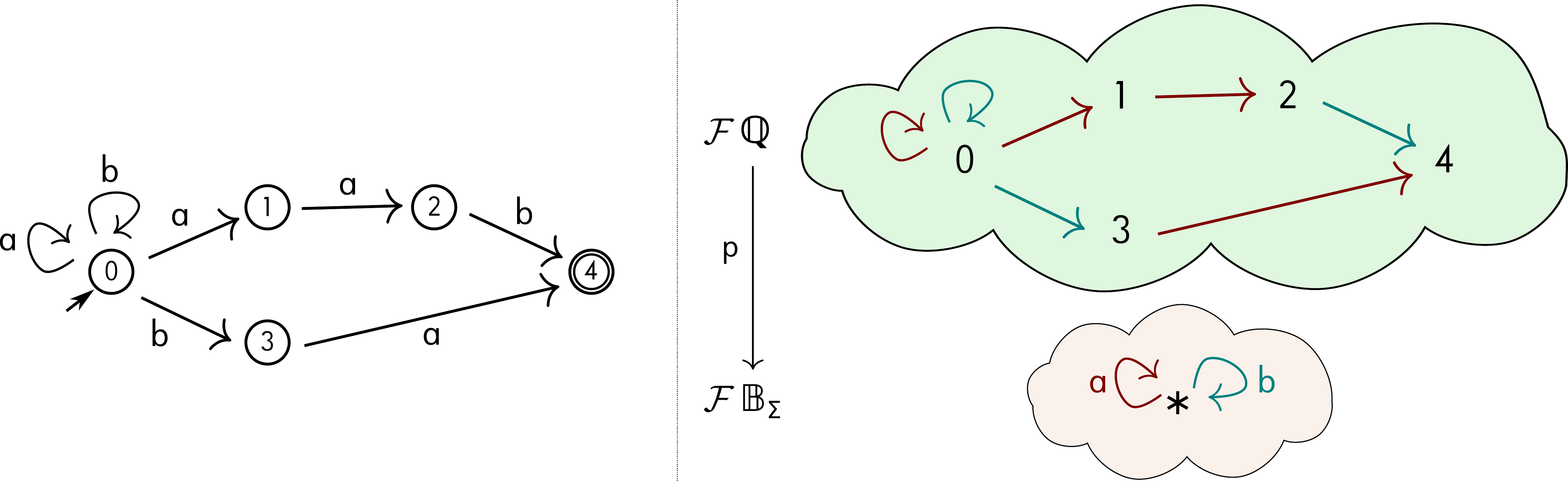} \end{center}
  \caption{Left: an NFA represented with a traditional state-diagram. Right: the bare NFA as a finitary ULF functor $p : \FreeCat{\cspec{Q}} \to \FreeCat{\Bouquet[\Sigma]}$. (We do not label the generators of $\FreeCat{\cspec{Q}}$ or indicate composite arrows, and we use colors to indicate the images of the generators in $\FreeCat{\Bouquet[\Sigma]}$.)}
  \label{fig:nfa}
\end{figure}

Classically, a nondeterministic finite-state automaton may be represented by a finite alphabet~$\Sigma$, a finite set~$Q$ of states, a finite set~$\TranSet$ of transitions together with a labelling function $\delta:\TranSet\to Q\times\Sigma\times Q$, and a choice of initial and accepting states.
We will focus first on the underlying ``bare'' automaton~$M = (\Sigma,Q,\delta)$ before the choice of initial and accepting states.
Every such bare automaton~$M$ induces a functor of categories $p:\cQ \to \FreeCat{\Bouquet[\Sigma]}$
where $\cQ$ is the category with set of objects $Q$ and with arrows freely generated
by arrows of the form $t:q\to q'$ for every transition $t\in\TranSet$ such that $\delta(t)=(q,a,q')$, 
and where the functor $p:\cQ\to \FreeCat{\Bouquet[\Sigma]}$ transports every such generator $t:q\to q'$ 
to the arrow $a:\ast\to\ast$ representing the letter~$a\in\Sigma$ in the category~$\FreeCat{\Bouquet[\Sigma]}$.
See Figure~\ref{fig:nfa} for an example.
Under this formulation,
any composite arrow $\alpha:q_0\to q_f$ of the category~$\cQ$ describes 
a run of the automaton~$M$ over the word $w=p(\alpha):\ast\to\ast$ which starts in state~$q_0\in Q$ and ends in state~$q_f\in Q$, 
as depicted below:
\begin{center}
\begin{tikzcd}[column sep=2em,row sep=.8em]
{q_0}\arrow[dd,|->]\arrow[rrrr,dashed,"{\alpha}"] &&&& {q_f}\arrow[dd,|->]
&&
\cQ\arrow[dd,"{p}"]
\\
\\
\ast\arrow[rrrr,"w"] &&&& \ast
&&
\FreeCat{\Bouquet[\Sigma]}
\end{tikzcd}
\end{center}
One distinctive property of the functor $p:\cQ \to \FreeCat{\Bouquet[\Sigma]}$ is that it has the \emph{unique lifting of factorizations} (ULF) property~\cite{Lawvere1986,BungeNiefeld2000,BungeFiore2000,LawvereMenni2010}.
Recall that a functor of categories has the ULF property (or ``is ULF'') when any factorization of the image of an arrow lifts uniquely to a factorization of that arrow.
\begin{defi}\label{definition/ULF-cat}
  A functor $p : \cD \to \cC$ is said to be \defin{ULF} if for every arrow $\alpha$ of $\cD$, if $p(\alpha)= uv$ for some pair of arrows~$u$ and~$v$ of $\cC$,
there exists a unique pair of arrows~$\beta$ and~$\gamma$ in $\cD$ such that $\alpha=\beta\gamma$ and $p(\beta)=u$ and $p(\gamma)=v$.
\end{defi}
\noindent
Such functors are also known as \emph{discrete Conduché fibrations} \cite{Johnstone1999,Guetta2020}, and they provide a common generalization of both discrete fibrations and discrete opfibrations.
The ULF condition arises naturally when considering the correspondence recalled in \S\ref{section/benabou} between functors $\cD \to \cC$ and lax functors $\cC\to\Span{\Set}$.
\begin{prop}\label{proposition/ULF-is-pseudo-Span}
  A functor of categories $p : \cD \to \cC$ is ULF iff the associated lax functor $F:\cC\to\Span{\Set}$ is a pseudofunctor.
\end{prop}
\begin{proof}
  The definition of ULF says precisely that the structure maps $F(u)F(v) \Longrightarrow F({uv})$ of the lax fiber functor are invertible, thus establishing the right-to-left direction.
  Moreover, it is an immediate consequence of the ULF property that $p(\alpha) = \id[A]$ implies $\alpha$ is an identity arrow, since otherwise there would be two distinct liftings $\alpha = \id\,\alpha$ and $\alpha = \alpha\,\id$ of the factorization $\id[A] = \id[A]\,\id[A]$.
  Hence the structure maps $\id[F(A)] \Longrightarrow F({\id[A]})$ are also invertible and $F$ is a pseudofunctor, thus establishing the left-to-right implication.
\end{proof}
\noindent
The ULF property reflects an important structural property of nondeterministic finite-state automata:
that every arrow $\alpha:q_0\to q_f$ lying above some arrow $p(\alpha)=w$ corresponding to a run of the automaton
can be factored uniquely as a sequence of transitions along the letters of the word~$w$.
In the case of the ULF functor $p : \cQ \to \FreeCat{\Bouquet[\Sigma]}$ that we construct from a classical NFA,
both the categories $\cQ$ and $\FreeCat{\Bouquet[\Sigma]}$ are freely generated, and 
the functor $p$ is generated from a homomorphism of graphs.
(In effect, on the right-hand side of Figure~\ref{fig:nfa} we visualize the functor by this generating graph homomorphism.)
This is not a coincidence, as \emph{any} ULF functor into a free category is necessarily of this form.
\begin{prop}[\cite{Street1996}, \cite{Guetta2020}]\label{ULF-into-free}
  Let $p : \cD \to \cC$ be a functor into a category $\cC = \FreeCat{\cgraph{G}}$ freely generated by some graph~$\cgraph{G}$.
  Then $p$ is ULF iff $\cD = \FreeCat{\cgraph{H}}$ is free over some graph $\cgraph{H}$ and $p = \FreeCat{\phi}$ is generated by a graph homomorphism $\phi : \cgraph{H} \to \cgraph{G}$.
\end{prop}
\begin{proof}
  The right-to-left direction is immediate. 
  For the left-to-right implication, we take $\cgraph{H} = p^{-1}(\cgraph{G})$ to be the graph containing just those arrows of $\cD$ lying over the generators of $\cC$.
  Since the image of any arrow $\alpha$ in $\cD$ uniquely decomposes as a composition of edges $p(\alpha) = e_1\cdots e_n$ in $\cgraph{G}$, the ULF property implies that $\alpha$ uniquely decomposes as a composition of edges in $\cgraph{H}$, and $p$ is generated by its restriction to the generators $\phi : \cgraph{H} \to \cgraph{G}$.
\end{proof}

Of course another important property of NFAs is their finite nature, which may also be expressed as a property of the corresponding functor $p:\cQ\to\FreeCat{\Bouquet[\Sigma]}$.
\begin{defi}
  We say that a functor $p : \cD \to \cC$ is \defin{finitary} if either of the following equivalent conditions hold:
  \begin{itemize}
  \item the fiber $p^{-1}(A)$ as well as the fiber $p^{-1}(w)$ is finite for every object $A$ and arrow $w$ in the category $\cC$;
  \item the associated lax functor $F:\cC\to\Span{\Set}$ factors via $\Span{\FinSet}$.
  \end{itemize}
\end{defi}
\noindent
The unique lifting of factorizations and finitary fiber properties characterize classical NFAs in the following sense.
\begin{prop}\label{finitary-into-finite-graph}
  Let $\phi : \cgraph{H} \to \cgraph{G}$ be a homomorphism into a finite graph $\cgraph{G}$. Then $\FreeCat{\phi}:\FreeCat{\cgraph{H}}\to\FreeCat{\cgraph{G}}$ is finitary iff $\cgraph{H}$ is finite.
\end{prop}
\begin{proof}By the infinite pigeonhole principle.
\end{proof}
\begin{prop}\label{proposition/ULF-nfa}
  A functor $p:\cQ\to\FreeCat{\Bouquet[\Sigma]}$ corresponds to the underlying bare automaton of a NFA 
  iff $p$ is ULF and finitary.
\end{prop}
\begin{proof}
  By Props.~\ref{ULF-into-free} and \ref{finitary-into-finite-graph}, $p$ is ULF and finitary iff $\cQ \cong \FreeCat{\cspec{Q}}$ is generated by a finite graph with $p$ induced by a homomorphism, thereby defining the transition graph of a bare NFA.
\end{proof}
\noindent
This leads us to the following definition of NFA over an arbitrary category.
\begin{defi}\label{definition/NFA-over-a-category}
  A \defin{nondeterministic finite-state automaton over a category} is given by a tuple
  $M = (\cC, \cQ, p, q_0, q_f)$
  consisting of two categories $\cC$ and $\cQ$, a finitary ULF functor $p:\cQ\to\cC$,
  and a pair $q_0,q_f$ of objects of $\cQ$.
  An object of~$\cQ$ is then called a \defin{state} and an arrow of~$\cQ$ is called a \defin{run} of the automaton.
  The \defin{regular language of arrows} $\Lang{M}$ recognized by the automaton
  is the set of arrows $w$ in $\cC$ that can be lifted along~$p$ to an arrow $\alpha:q_0\to q_f$ in $\cQ$, that is $\Lang{M} = \set{p(\alpha) \mid \alpha : q_0 \to q_f} \subseteq \cC(A,B)$, where $p(q_0) = A, p(q_f) = B$.
\end{defi}
\noindent
Similarly to the way a categorical CFG may be considered as a functor of pointed operads (see Remark~\ref{remark/pointing}), it is also sometimes helpful to view a categorical NFA as a finitary ULF functor of bipointed categories $p_M : (\cQ,q_0,q_f) \to (\cC,A,B)$.
We emphasize that the definition does not require the base category~$\Ccategory$ to be free, which we will see in examples later on.

Before continuing to develop this definition, let us briefly comment on our treatment of initial and accepting states, as well as the treatment of determinism and $\epsilon$-transitions.
Rabin and Scott's original definition of NFA \cite{RabinScott1959} (also adopted by Eilenberg \cite{Eilenberg1974}), takes a set of initial states together with a set of accepting states, while some more recent authors \cite{HMU2007,Sipser2013} take a single initial state together with a set of accepting states.
By Rabin and Scott's determinization theorem, a language $\Lang{} \subseteq \Sigma^*$ may be recognized by a NFA with a set of initial and accepting states just in case it may be recognized by a DFA with a single initial state and a set of accepting states, so the choice of having a single initial state or a set of initial states does not matter in terms of recognition power.
Classically, the same is not quite true for NFAs over a fixed alphabet with a single initial state and a single accepting state (in the absence of $\epsilon$-transitions), since the language recognized by such an automaton is necessarily closed under concatenation if it contains the empty word.
Nevertheless, it is easy to see that any NFA with a set of initial and accepting states may be converted to an equivalent NFA with a single initial state and a single accepting state, over an alphabet extended with distinguished symbols marking the beginning and the end of the word.
This construction is particularly natural to express using the category $\FreeCat{\Bracket{\Bouquet[\Sigma]}}$ introduced in Example~\ref{example/spliced-ordinal-sum}.
\begin{prop}\label{proposition/classical-NFAs}
  A language $\Lang{} \subseteq \Sigma^*$ is regular in the classical sense if and only if $\bow \Lang{}\eow$  is the language of arrows
  recognized by a nondeterministic finite-state automaton over $\FreeCat{\Bracket{\Bouquet[\Sigma]}}$.
\end{prop}
\begin{proof}
  Given a classical NFA $M$ with underlying bare automaton $\cQ \to \FreeCat{\Bouquet[\Sigma]}$, we construct an NFA $\Bracket{M} = (\FreeCat{\Bracket{\Bouquet[\Sigma]}}, \cQ', p' : \cQ' \to \FreeCat{\Bracket{\Bouquet[\Sigma]}}, q_0, q_f)$ where $\cQ'$ and the functor $p'$ are  constructed from $\cQ$ and $p$ by freely adjoining a pair of objects $q_0$ and $q_f$ lying over $\bot$ and $\top$ respectively, together with an arrow $q_0 \to q$ over $\bow : \bot \to *$ for every initial state $q$ and an arrow $q' \to q_f$ over $\eow : * \to \top$ for every accepting state $q'$ of $M$.
  By construction, every accepting run of $M$ over the word $w$ corresponds to an accepting run of $\Bracket{M}$ over the arrow $\bow w\eow : \bot \to \top$.
\end{proof}
\noindent
Although Definition~\ref{definition/NFA-over-a-category} could be adapted to take a set of initial and accepting states, we find the convention of taking a single pair $(q_0,q_f)$ of initial and accepting states to align well with the definition of context-free languages of arrows generated by a single nonterminal (Definition~\ref{definition/cfg-over-cat}). 
Let us also observe here that the category $\FreeCat{\Bracket{\Bouquet[\Sigma]}}$ in the statement of Proposition~\ref{proposition/classical-NFAs} is isomorphic to the \emph{input category} $\mathcal{I}_{\mathsf{word}}$ used in Colcombet and Petri\c{s}an's definition of categorical automata. (See bottom of Figure~1 of \cite{ColcombetPetrisan2020}; we will give a more detailed comparison with Colcombet and Petri\c{s}an's framework in \S\ref{section/automata/related} below.)

\begin{defi}\label{definition/determinism}
  An NFA $M = (\cC, \cQ, {p : \cQ \to \cC}, q_0, q_f)$ is \defin{deterministic} if the functor $p$ is a \emph{discrete opfibration}, i.e., for any arrow $w : A \to B$ in $\cC$ and any state $q \in \cQ$ such that $p(q) = A$, there is a unique arrow $\alpha : q \to q'$ such that $p(\alpha) = w$. 
  $M$ is \defin{codeterministic} if $p$ is a \emph{discrete fibration}, i.e., for any arrow $w : A \to B$ in $\cC$ and any state $q' \in \cQ$ such that $p(q') = B$, there exists a unique arrow $\alpha : q \to q'$ such that $p(\alpha) = w$. 
  $M$ is \defin{bideterministic} if it is both deterministic and codeterministic, i.e., $p$ is a discrete bifibration.
  Finally, $M$ is a \defin{partial} (co/bi)deterministic automaton if the functor $p : \cQ \to \cC$ is a partial discrete op/bi/fibration, in the sense that for any $w : A \to B$ in $\cC$ and $q \in \cQ$ such that $p(q) = A$ (or $q'$ such that $p(q') = B$), there is at most one $\alpha : q \to q'$ such that $p(\alpha) = w$.
\end{defi}
\begin{rem}\label{remark/epsilon}
  With this way of representing finite-state automata by functors $p : \cQ \to \cC$, $\epsilon$-transitions may be naturally interpreted as non-identity arrows $\alpha : q \to q'$ in $\cQ$ mapped by $p$ to an identity arrow $\id[A] : A \to A$ in $\cC$.
  However, it is an immediate consequence of the ULF property that an NFA in the sense of Definition~\ref{definition/NFA-over-a-category} cannot contain any non-trivial fiber categories $p^{-1}(A)$ (Proposition~\ref{proposition/ULF-is-pseudo-Span}).
  One would therefore need a more general notion of ULF functor in order to accommodate categorical NFAs with $\epsilon$-transitions.
  This seems to raise some subtle issues.
  In particular, the standard treatment of $\epsilon$-removal \cite[\S2.5]{HMU2007} seems to indicate that one needs to design a notion between ULF functors and general Conduché fibrations.
  Recall that the standard notion of Conduché fibration weakens the ULF property by only requiring uniqueness of factorizations up to zigzag in the fiber, which means here zigzag of $\epsilon$-transitions.
  The right definition of categorical NFA with $\epsilon$-transitions should probably ensure more, that there is a universal (initial or terminal) solution in the space of possible factorizations.
\end{rem}
\noindent
As we have seen, our definition of NFA over a category does not require the base category $\cC$ to be free.
We now list some examples of such automata and the regular languages of arrows that they generate.
\begin{exa}[Product automata]\label{example/product-automata}
  The product of two finitary ULF functors $p : \cQ \to \cC$ and $p' : \cQ' \to \cC'$ is again a finitary ULF functor $p \times p' : \cQ \times \cQ' \to \cC \times \cC'$.
  This enables us to construct the product of two NFAs $M$ and $M'$ as an NFA $M\times M'$ whose underlying bare automaton is $p\times p'$ and whose initial and accepting states are given by the pairing of the respective initial states $(q_0,q_0')$ and $(q_f,q_f')$.
  Note that even if $M$ and $M'$ are NFAs over free categories $\cC = \FreeCat{\cgraph{G}}$ and $\cC' = \FreeCat{\cgraph{G}'}$, the product automaton $M\times M'$ will not be since the product of free categories is not a free category.
  Intuitively, $M\times M'$ may be interpreted as an automaton that reads pairs of arrows in parallel.
  It recognizes the cartesian product $\Lang{}\times \Lang[']{}$ of the languages recognized by $M$ and $M'$.
\end{exa}
\begin{exa}[Singleton automata]\label{example/singleton}
  For any word $w = a_1 \dots a_n$ of length $n$, there is a partial bideterministic $(n+1)$-state automaton $M_w$ that recognizes the singleton language $\set{w}$, with initial state $0$, accepting state $n$, and transitions of the form $(i,a_{i+1},i+1)$ for each $0 \le i < n$.
  This construction may be generalized to define a (not necessarily deterministic or codeterministic) NFA $M_w$ recognizing any arrow $w$ in any category satisfying a certain finitariness assumption.
  Indeed, for any arrow $w : A \to B$ of a category $\cC$, the \emph{category of factorizations} $\Fact{w}$ of $w$ (also called the \emph{interval category} of $w$ in \cite{LawvereMenni2010}) is defined by taking objects to be triples $(X,u,v)$ of an object $X \in \cC$ and a pair of arrows $u : A \to X$, $v : X \to B$ such that $w = uv$, and arrows  $(X,u,v) \to (X',u',v')$ to be arrows $x : X \to X'$ such that $u' = ux$ and $v = xv'$.
  This category has an initial object $(\id[A],w)$ and a terminal object $(w,\id[B])$, and it comes equipped with an evident forgetful functor $\mu_w : \Fact{w} \to \cC$ returning the middle of a factorization, which is always ULF \cite{Johnstone1999,LawvereMenni2010}.
  Hence 
  \begin{equation}\label{equation/nfa-Mw}
    M_w = (\cC,\Fact{w},\mu_w,(\id[A],w),(w,\id[B]))
  \end{equation}
  may be considered as a nondeterministic automaton recognizing the singleton language $\set{w}$, although it is not necessarily a finite-state automaton since the functor $\mu_w$ need not be finitary in general.
  Let us say that a category $\cC$ \defin{has finitary factorizations} if this is always the case, in other words, if the ULF functor $\mu_w : \Fact{w} \to \cC$ is finitary for every arrow $w$ of $\cC$.
  Then we see that to any arrow $w$ of a category $\cC$ with finitary factorizations, \eqref{equation/nfa-Mw} defines an NFA $M_w$ recognizing exactly that arrow.

  Free categories always have finitary factorizations, and indeed the category of factorizations of any arrow of a free category is isomorphic to a non-empty finite ordinal, as in the example above of the $(n+1)$-state automaton associated to a word of length $n$.
  
  A category with finitely many objects has finitary factorizations just in case every arrow has finitely many factorizations $w = uv$ of length 2.\footnote{Note that having finitary factorizations is weaker than the \emph{Möbius} condition on a category \cite{LawvereMenni2010}, which asks that the set of non-trivial decompositions $\set{(u_1,\dots,u_n) \mid n \ge 0, w = u_1\cdots u_n (\forall 1\le i \le n, u_i \ne \id)}$ of every arrow $w$ is finite (cf.~\cite[Prop.~2.6]{LawvereMenni2010}).
  As a property of monoids, having finitary factorizations is standardly called the \emph{finite decomposition} property, which is weaker than what Eilenberg called the \emph{local finiteness} property \cite[\S VII.4, p.170]{Eilenberg1974} corresponding to the Möbius condition.}
  For example, consider the one-object category $\cE$ containing a single non-identity arrow $f : * \to *$ that is idempotent $f f = f$.
  The NFA $M_f$ recognizing the language $\set{f} \subseteq \cE(\ast,\ast)$ 
  has three states $(\id,f)$, $(f,\id)$, and $(f,f)$, as depicted by the following transition diagram in which all paths of length $\ge 1$ commute:
  \[
    \begin{tikzpicture}[thick,scale=0.75, every node/.style={transform shape}]
      \node[state,initial] (q0) at (0,0) {$\id,f$};
      \node[state,right of=q0,accepting] (q1) {$f,\id$};
      \node[state,below of=q1] (q2) {$f,f$};
      \draw
      (q0) edge node[above]{$f$} (q1)
      (q0) edge node[below left]{$f$} (q2)
      (q1) edge [bend left] node[right]{$f$} (q2)
      (q2) edge node[left]{$f$} (q1)
      (q2) edge [loop right] node{$f$} (q2);
    \end{tikzpicture}
    \]
\end{exa}

\begin{exa}[Total automata]\label{example/total-automaton}
  The identity functor $\id[\cC] : \cC \to \cC$ on a category is trivially ULF and finitary, since every object and arrow has a unique lifting.
  Hence for any pair of objects $A$ and $B$ of $\cC$, there is an NFA $M_{\cC(A,B)} = (\cC,\cC,\id[\cC],A,B)$ recognizing the total language $\cC(A,B)$, i.e., accepting every arrow $w : A \to B$.
  In the case of the one-object category $\cC = \FreeCat{\Bouquet[\Sigma]}$, this is just the familiar automaton with a single (initial and accepting) state and a transition on that state for every letter of $\Sigma$.
  On the other hand, it is worth emphasizing that $M_{\cC(A,B)}$ defines an NFA in our sense even when $\cC$ has infinitely many objects.
  In particular, this example shows that what is important is finiteness of the \emph{fibers} of the functor $p : \cQ \to \cC$ defining an NFA, and not of the total category of states $\cQ$.
\end{exa}

\begin{exa}[Asynchronous automata]
One good reason for considering NFAs whose base category~$\cC$
is not free is that they provide for interesting notions of automata 
when the category~$\cC$ is presented by generators and relations.
Let us say that an \emph{equation} in a graph~$\cgraph{G}$ is a pair $(f,g)$ of paths $f,g:A\to B$ 
with same source and target.
Every set~$\mathcal{R}$ of such equations
induces a congruence~$\sim_{\mathcal{R}}$
on the paths of~$\cgraph{G}$,
defined as the smallest equivalence relation 
between paths of $\cgraph{G}$ with same source and target, 
compatible to composition, and
such that $f\sim_{\mathcal{R}}g$
for every equation $(f,g) \in \mathcal{R}$.
%
%
We say that a category $\cC$ is \emph{presented} by a graph $\cgraph{G}$
of generators and a set $\mathcal{R}$ of equations
when it comes equipped with an isomorphism 
$c:\FreeCat{\cgraph{G}}/{\sim_{\,\mathcal{R}}} \to \cC$ 
from the category $\FreeCat{\cgraph{G}}/{\sim_{\,\mathcal{R}}}$ whose objects are nodes of $\cgraph{G}$ and whose arrows are $\sim_{\,\mathcal{R}}$-equivalence classes of paths.
%
%

A remarkable property of ULF functors $p : \cD \to \cC$ into a category $\cC$
presented by a graph~$\cgraph{G}$ and a set~$\mathcal{R}$ of equations,
with isomorphism $c:\FreeCat{\cgraph{G}}/{\sim_{\,\mathcal{R}}} \to\cC$,
is that the category~$\cD$ inherits a presentation by a graph $\pull{p}\cgraph{G}$ 
and a set $\pull{p}\mathcal{R}$ of equations derived from $\cgraph{G}$ and $\mathcal{R}$.
The graph $\pull{p}\cgraph{G}$ has nodes the pairs $(x,A)$
of a node $x$ in $\cgraph{G}$ and an object~$A$ of~$\cD$
such that $c(x)=p(A)$, and edges the pairs $(e, u)$ of an edge $e:x\to y$ in $\cgraph{G}$
and an arrow $u:A\to B$ in $\cD$ such that $c(e)=p(u)$.
The set $\pull{p}\mathcal{R}$ contains all equations $(f,g)$
between paths of $\pull{p}\cgraph{G}$ of the form
$$
\begin{tikzcd}[cells={nodes={font=\normalsize}}]
f=(e_1,u_1)\cdots (e_m,u_m) \, , \, g=(e'_1,v_1)\cdots (e'_n,v_n)
\quad : \quad (x,A)\arrow[r] & (y,B)
\end{tikzcd}
$$
such that $(e_1\cdots e_m, e'_1\cdots e'_n) \in \mathcal{R}$ 
and the composite arrows
$u_1\cdots u_m = v_1\cdots v_n$
are equal in the category~$\cD$.
The isomorphism $d : \pull{p}\cgraph{G}/{\sim_{\,\pull{p}\mathcal{R}}} \to \cD$ is given by the second projection $(x,A) \mapsto A$ on objects, and by the map $[e_1\cdots e_m,u_1 \cdots u_m] \mapsto u_1\cdots u_m$ sending an equivalence class of lists of edge/arrow pairs to the composition of the arrows in $\cD$.

%
We illustrate the property by giving a direct description of NFAs
on a category~$\cC$ defined by a Mazurkiewicz alphabet.
A concurrent alphabet $(\Sigma,\Dependency{\Sigma})$
as defined by Mazurkiewicz \cite{Mazurkiewicz89}
is a set~$\Sigma$ of letters, equipped with a reflexive
and symmetric binary relation~$\Dependency{\Sigma}\subseteq \Sigma\times\Sigma$.
%
Two letters $a,b\in\Sigma$ are called \emph{dependent} when $(a,b)\in\Dependency{\Sigma}$
and \emph{independent} otherwise.
%
%
The \emph{permutation equivalence} on finite words is defined as the smallest 
equivalence relation ${\permeq{\Sigma}}$
such that $w_1\,a\,b\,w_2\permeq{\Sigma}w_1\,b\,a\,w_2$ for every pair of words 
$w_1,w_2\in\Sigma^{\ast}$ and independent letters $a,b\in\Sigma$.
A \emph{trace} is an equivalence class of the permutation equivalence.
Traces $w\in\Sigma^{\ast}\!/\!\permeq{\Sigma}$ may be regarded
as the arrows $w : \ast\to\ast$ of a one-object category $\TraceCat$
obtained by quotienting the category $\FreeCat{\Bouquet[\Sigma]}$
with an equation between arrows $a\,b=b\,a:\ast\to\ast$
for every pair of independent letters $a,b\in\Sigma$.
%

%
%

Unpacking the definitions, a NFA over $\TraceCat$
may be equivalently described as follows, starting from the data of
a finite set $Q$ of states, a finite set~$\TranSet$ of transitions 
together with a labelling function $\delta:\TranSet\to Q\times\Sigma\times Q$, 
and a choice of initial and accepting states.
%
The transition graph is moreover equipped with a set $\diamond$ of permutation tiles 
defined as pairs $(\alpha,\beta)$ of paths $\alpha=st$ and $\beta=t's'$
of length $2$ of the form 
\begin{center}
$\delta(s)=(q,a,q')$, $\delta(t)=(q',b,q'')$,
$\delta(t')=(q,b,q''')$ and $\delta(s')=(q''',a,q'')$
\end{center}
where $q,q',q'',q'''\in Q$ are states and the two letters~$a,b\in\Sigma$ are independent.
We write $\alpha\diamond\beta$ when $(\alpha,\beta)\in \diamond$.
%
%
The set of permutation tiles is required to satisfy the following properties:

\begin{enumerate}
\item {permutation tiles are deterministic:}
$\alpha\diamond \beta$ and $\alpha\diamond\beta'$ implies $\beta=\beta'$,

\item {permutation tiles are symmetric:} 
$\alpha\diamond \beta$ implies $\beta\diamond \alpha$,

%
%
%
%

\item {successive transitions with independent letters can permute:} 
for every transition path $\alpha=s\, t$ of length~$2$
with $\delta(s)=(q,a,q')$ and $\delta(t)=(q',b,q'')$
where the two letters~$a,b\in\Sigma$ are independent,
there exists a permutation tile $\alpha\diamond\beta$.
\end{enumerate}
%
%
%
%
%

\medbreak
The correspondence works as follows.
Given a finite transition system $(Q,\TranSet,\delta)$ 
equipped with a set $\diamond$ of permutation tiles,
one defines the category~$\cC$ presented 
by the corresponding graph~$\cgraph{G}$ 
with an equation $(\alpha,\beta)$ for every permutation tile $\alpha\diamond\beta$ ;
the three properties (1), (2), and (3) ensure then that 
the labelling functor $q:\cC\to \TraceCat$
is finitary ULF.
%
Conversely, every finitary ULF functor
$p:\cC\to\TraceCat$
induces a finite graph $\cgraph{G}=\pull{p}{\Bouquet[\Sigma]}$ 
and a set of equations $\diamond = \pull{p}\mathcal{R}_{\Dependency{}}$
where $\mathcal{R}_{\Dependency{}}$ denotes the set of equations
$(ab,ba)$ for $a$ and $b$ two independent letters ;
one checks that the nodes and edges of~$\cgraph{G}$
define a transition system $(Q,\TranSet,\delta)$
which satisfies properties (1), (2), and (3).
Note that property (3) reflects the ULF property 
for the factorisation 
of the arrow $ab$ in $\TraceCat$ as $b$ followed by $a$,
for any pair of independent letters $a$ and $b$ in $\Sigma$.

This establishes that the notion of NFA over 
the category $\TraceCat$
of Mazurkiewicz traces captures a natural notion of asynchronous automaton,
corresponding to a nondeterministic version of Bednarczyk's 
asynchronous systems \cite{Bednarczyk88,HussonMorin00,Morin05,BaudruMorin06}
where property (3) is called the Independent Diamond (ID) property.
A detailed comparison of this class of asynchronous automaton
with the more liberal notion of automaton with concurrency relations~\cite{DrosteShortt02,GoubaultMimram12} appears in~\cite{Morin05}.
The ID property implies that
every NFA over $\FreeCat{{\Bouquet[\Sigma]}}/\Dependency{}$
satisfies the cube property (CP) discussed in \cite{HussonMorin00}
and thus defines what is called an asynchronous graph in~\cite{MelliesMimram07,Mellies21}. 
We leave for future work the connection with Zielonka's theorem
which is based on a more concrete and process-based notion
of concurrent automaton~\cite{Zielonka1987,BaudruMorin06}.
%
%
%
\end{exa}

\subsection{From word automata to tree automata and beyond}
\label{section/tree-automata}

Reassuringly, the foregoing fibrational analysis of nondeterministic finite-state automata based on finitary ULF functors adapts smoothly when one shifts from word automata to tree automata \cite{tata2008}.
As a first step in that direction, we describe how the ULF and finite fiber properties may be extended to functors of operads.
\begin{defi}\label{definition/ULF}
  A functor of operads $p:\oD\to\oO$ has the \defin{unique lifting of factorizations property} (or \defin{is ULF}) if any of the following equivalent conditions hold:
  \begin{enumerate}
  \item
    for any operation $\alpha$ of $\oD$, if $p(\alpha) = g \circ (h_1,\dots,h_n)$ for some operation $g$ and list of operations $h_1,\dots,h_n$ of $\oO$, there exists a unique operation $\beta$ and list of operations $\gamma_1,\dots,\gamma_n$ of $\oD$ such that $\alpha = \beta \circ (\gamma_1,\dots,\gamma_n)$ and $p(\beta) = g$, $p(\gamma_1) = h_1, \dots, p(\gamma_n) = h_n$;
  \item
    for any operation $\alpha$ of $\oD$, if $p(\alpha) = g \circ_i h$ for some operations $g$ and $h$ of $\oO$ and index $i$, there exists a unique pair of operations $\beta$ and $\gamma$ of $\oD$ such that $\alpha = \beta \circ_i \gamma$ and $p(\beta) = g$ and $p(\gamma) = h$;
  \item
    the structure maps of the associated lax functor of operads $F:\oO\to\Span{\Set}$ discussed in \S\ref{section/benabou} are invertible.
  \end{enumerate}
\end{defi}
\begin{defi}
  We say that a functor of operads $p : \oD \to \oO$ is \defin{finitary} if either of the following equivalent conditions hold:
  \begin{itemize}
  \item the fiber $p^{-1}(A)$ as well as the fiber $p^{-1}(f)$ is finite for every color $A$ and operation $f$ of the operad $\oO$;
  \item the associated lax functor of operads $F:\oO\to\Span{\Set}$ factors via $\Span{\FinSet}$.
  \end{itemize}
\end{defi}
\noindent
Just as we saw in the previous section, 
one can easily check that the underlying bare automaton $M=(\Sigma,Q,\delta)$ of a (bottom-up or top-down) nondeterministic finite-state tree automaton (cf.~\S1.1 and \S1.6 of \cite{tata2008}) gives rise to a finitary ULF functor of operads $p:\oQ\to\FreeOper{\Sigma}$,
where 
\begin{itemize}
\item $\FreeOper{\Sigma}$ is the free operad generated by the ranked alphabet $\Sigma$ (seen as an uncolored species);
\item the operad~$\oQ$ has states of the automaton as colors and operations freely generated by $n$-ary nodes of the form $t : q_1,\dots,q_n \to q$ for every transition $t\in\TranSet$ of the form $\delta(t)=(q_1,\dots,q_n,a,q)$, with $a \in \Sigma_n$;
\item the functor $p$ transports every such $n$-ary transition $t : q_1,\dots,q_n \to q$ to the underlying $n$-ary letter $a:\ast,\dots,\ast\to\ast$.
\end{itemize}
Moreover, this correspondence is one-to-one.
\begin{prop}\label{ULF-into-free-operad}
  Let $p : \oD \to \oO$ be a functor into an operad $\oO = \FreeOper{\Sigma}$ freely generated by some species~$\Sigma$.
  Then $p$ is ULF iff $\oD = \FreeOper{\Sspecies}$ is free over some species $\Sspecies$ and $p = \FreeOper{\phi}$ is generated by a map of species $\phi : \Sspecies \to \Sigma$.
\end{prop}
\begin{prop}\label{proposition/ULF-nfa-trees}
  A functor of operads $p:\oQ\to\FreeOper{\Sigma}$ corresponds to a bare nondeterministic finite-state tree automaton
  iff $p$ is ULF and finitary.
\end{prop}
\noindent
This motivates us to proceed as for word automata and propose a more general notion 
of finite-state automaton over an arbitrary operad.

\begin{defi}\label{definition/NFA-over-an-operad}
  A \defin{nondeterministic finite-state automaton over an operad} is given by a tuple
  $M = (\oO, \oQ, p : \oQ \to \oO, q_r)$
  consisting of two operads $\oO$ and $\oQ$, a finitary ULF functor of operads $p:\oQ\to\oO$,
  and a color $q_r$ of $\oQ$.
  A color of~$\oQ$ is called a \defin{state}, and an operation of~$\oQ$ is called a \defin{run tree} of the automaton~$p:\oQ\to\oO$.
  The \defin{regular language of constants} $\Lang{M}$ recognized by the automaton
  is the set of constants $c$ in $\oO$ that can be lifted along~$p$ to a constant $\alpha:q_r$ in $\oQ$, that is $\Lang{M} = \set{p(\alpha) \mid \alpha : q_r} \subseteq \oO(A)$, where $p(q_r) = A$.
\end{defi}
\noindent
Let us refer to NFAs over categories in the sense of Definition~\ref{definition/NFA-over-a-category} as \emph{categorical automata} for short, and NFAs over operads in the foregoing sense as \emph{operadic automata}.

We now state a simple property of ULF and finitary functors establishing a useful connection
between categorical automata and operadic automata.
\begin{prop}\label{proposition/Words-ULF}
Suppose that $p:\cQ\to\cC$ is a functor of categories.
If $p$ is ULF (respectively, finitary) then so is the functor of operads~$\Words{p}:\Words{\cQ}\to\Words{\cC}$.
\end{prop}
\begin{proof}
  Follows immediately from the fact that a lifting of an operation $f_0{-}\dots{-}f_n$ in $\Words{\cC}$ along $\Words{p}$ is a sequence of liftings of the arrows $f_0,\dots,f_n$ in $\cC$ along $p$.
\end{proof}
\begin{cor}\label{corollary/WM-automaton}
  For any categorical automaton $M = (\cC,\cQ,p,q_0,q_f)$ there is an associated operadic automaton $\Words{M} = (\Words{\cC},\Words{\cQ},\Words{p},(q_0,q_f))$ generating the same language $\Lang{M}=\Lang{\Words{M}}$.
\end{cor}
\begin{proof}
  By Proposition~\ref{proposition/Words-ULF}, and because constants of $\Words{\cC}$ are exactly arrows of $\cC$.
\end{proof}
\noindent
We will see that the $\Words{M}$ construction transforming an automaton over a category into an automaton over its operad of spliced arrows is fundamental for computing the intersection of a context-free language with a regular language of arrows. 
Moreover, we highlight that since operads of spliced arrows are not freely generated (Remark~\ref{remark/non-free}), the $\Words{M}$ construction provides an important example of how operadic automata truly generalize classical tree automata.
\begin{rem}\label{remark/DFA-over-spliced-arrows}
As an aside, we mention that it is also interesting to consider \emph{deterministic} finite-state automata over spliced arrow operads.
To see why, let us first observe that any discrete opfibration over $\Words{\cC}$, which may be viewed equivalently as a functor of operads $\Words{\cC} \to \Set$, induces a category $\cD$ equipped with an identity-on-objects functor of categories $\cC \to \cD$, and conversely.
Indeed, given a functor $F : \Words{\cC} \to \Set$ one can construct a new category $\cD_F$ with same objects as $\cC$ and whose arrows $\alpha : A \to B$ are given by the elements of $F(A,B)$, with composition and identities defined by the interpretations
\begin{align*}
  F(\id[A]{-}\id[B]{-}\id[C]) &: F(A,B)\times F(B,C) \longrightarrow F(A,C) \\
  F(\id[A]) &: F(A,A)
\end{align*}
of the corresponding spliced arrow operations.
Moreover, there is an identity-on-objects functor $\cC \to \cD_F$ given by sending every arrow $w : A \to B$ of $\cC$, seen as a constant $w : (A,B)$ of $\Words{\cC}$, to its interpretation $F(w) : F(A,B)$.
Conversely, given a category $\cD$ and an identity-on-objects functor $h : \cC \to \cD$, one can define $F_h : \Words{\cC} \to \Set$ by taking $F_h(A,B) = \cD(A,B)$ and defining the interpretation of the spliced arrow operations by mapping the arrows of $\cC$ into $\cD$ along $h$ and then composing:
\[
F_h(w_0{-}w_{1}\dots{-}w_n) = (u_1,\dots,u_n) \mapsto h(w_0)u_1 h(w_1)\dots u_n h(w_n)
\]
From the description of this correspondence, it is clear it restricts to one between \emph{finitary} discrete opfibrations over $\Words{\cC}$ (equivalent to functors $\Words{\cC} \to \FinSet$) and \emph{locally finite} categories $\cD$ equipped with an identity-on-objects functor $\cC \to \cD$.
Thus what is a DFA over $\Words{\cC}$?
Nothing more than a locally finite category $\cD$, an identity-on-objects functor $h : \cC \to \cD$, and an arrow $q_r : A \to B$ of $\cD$, with the DFA accepting an arrow $w : A \to B$ of $\cC$ just in case $h(w) = q_r$.
In the case that $\cC$ is a one-object category, this is essentially the notion of \emph{recognition by a finite monoid,} which is one of the classical equivalent characterizations of regular languages \cite[\S III.10, p.62]{Eilenberg1974}.\footnote{A language $\Lang{} \subseteq \Sigma^*$ is said to be recognizable by a finite monoid $X$ if there exists a subset $S \subseteq X$ and a homomorphism $h : \Sigma^* \to X$ such that $\Lang{} = h^{-1}(S)$.  Our formulation requires $S$ to be a singleton, cf.~Proposition~\ref{proposition/classical-NFAs} above and the preceding discussion.}
\end{rem}

\subsection{Additional properties of ULF and finitary functors}
\label{section/properties-ULF-finitary}

We state here a few more basic properties of ULF and finitary functors, in both categorical and operadic versions, that will be important in the rest of the paper.
\begin{prop}\label{proposition/ULF-fin-pullback-closed}
  If $p : Y \to X$ is a ULF (resp.~finitary) functor of categories or operads, and $G : Z \to X$ is an arbitrary functor, then the functor $\pull{G}p : Z\times_{X} Y \to Z$ defined by pulling back $p$ along $G$ is ULF (resp.~finitary).
  \[
  \begin{tikzcd}
    Z\times_XY \arrow[d,"\pull{G}p\ \text{ULF/finitary}"']\arrow[r]
    \arrow[rd, phantom, "\lrcorner" very near start]
    & Y \arrow[d,"p\ \text{ULF/finitary}"] \\
    Z\arrow[r,"G"'] & X
  \end{tikzcd}
  \]
\end{prop}
\noindent
Closure of both classes of functors under pullback along arbitrary functors may be verified directly, or may be seen as an immediate consequence of the fact that both classes of functors $p : Y \to X$ admit ``Grothendieck-type'' correspondences with pseudo or lax functors $F : X \to \mathcal{U}$ into a certain universe (namely, ULF functors as pseudo functors into $\mathcal{U} = \Span{\Set}$, and finitary functors as lax functors into $\mathcal{U} = \Span{\FinSet}$), so that the pullback along an ordinary functor $G : Z \to X$ is simply the composite $F \circ G : Z \to \mathcal{U}$.

As we already saw in Example~\ref{example/total-automaton}, the identity functor $\id[\cC] : \cC \to \cC$ on any category is trivially finitary and ULF, as is the identity functor $\id[\oO] : \oO \to \oO$ on any operad.
It is also easy to check that both classes of functors are closed under composition.
Moreover, ULF functors satisfy a property analogous to the pasting law for pullbacks in a category, while finitary functors behave similarly to monomorphisms.
\begin{prop}\label{proposition/ULF-finitary-triangle} Consider a commutative triangle of functors of categories or operads:
\[
  \begin{tikzcd}
    X \ar[rd,"H"'] \ar[rr,"F"] && Y \ar[ld,"G"] \\
    & Z
  \end{tikzcd}
\]
\begin{enumerate}
\item if $F$ and $G$ are ULF, then so is $H$
\item if $F$ and $G$ are finitary, then so is $H$
\item if $G$ and $H$ are ULF then so is $F$
\item if $H$ is finitary then so is $F$
\end{enumerate}
\end{prop}
\begin{proof}
  For concreteness, we take $X$, $Y$, and $Z$ to be categories. (The case where they are operads is completely analogous.)
  \begin{enumerate}
  \item Let $\alpha$ be an arrow in $X$. Any factorization $H(\alpha) = uv$ in $Z$ lifts uniquely to a factorization of $\alpha = \beta\gamma$ with $H(\beta) = u$, $H(\gamma) = v$, by first lifting along $G$ to obtain a unique factorization of $F(\alpha)$, and then lifting along $F$.
  \item All of the fibers of $H$ are finite (for both the objects and arrows of $Z$) since they may be decomposed as a finite union of finite fibers
  $H^{-1}(z) = \bigcup_{y \in G^{-1}(z)} F^{-1}(y)$.
\item Let $\alpha$ be an arrow in $X$ and $F(\alpha) = st$ a factorization in $Y$.
  By projecting down along $G$ and then lifting the resulting factorization $H(\alpha) = G(s)G(t)$ along $H$, we obtain unique $\beta$ and $\gamma$ such that $\alpha = \beta\gamma$, $H(\beta) = G(s)$, $H(\gamma) = G(t)$.
  But since $F(\alpha) = F(\beta)F(\gamma)$ is a lifting of the factorization $G(st) = G(s)G(t)$ along $G$, this implies that $F(\beta) = s$ and $F(\gamma) = t$ by the ULF property of $G$.
  \item All of the fibers $F^{-1}(y)$ are finite, since they are contained in the fibers $H^{-1}(G(y))$. \qedhere
  \end{enumerate}
  
\end{proof}
\noindent
We remark in passing that finitary ULF functors may be arranged into various categories and bicategories of interest, which we describe now for the case of categorical automata (the operadic case again being completely analogous), although we will not make explicit use of these (bi)categories in the paper.
For any category $\cC$, finitary ULF functors over $\cC$ and commutative triangles
\[
  \begin{tikzcd}
    \cQ \ar[rd,"\text{finULF}"'] \ar[rr,"F"] && \cQ' \ar[ld,"\text{finULF}"] \\
    & \cC
  \end{tikzcd}
\]
define a category of bare NFAs over $\cC$ and morphisms between them, where the functor $F : \cQ \to \cQ'$ may be seen as a \defin{functional simulation} since it sends any run $\alpha : q_0 \to q_f$ in $\cQ$ over an arrow $w : A \to B$ of $\cC$ to a run $F(\alpha) : F(q_0) \to F(q_f)$ in $\cQ'$ over the same arrow.
By Proposition~\ref{proposition/ULF-finitary-triangle}, $F$ is itself finitary ULF, and this category is just the slice $\finULF/\cC$ of the category $\finULF$ of small categories and finitary ULF functors between them.
It is possible to formulate an analogue of Proposition~\ref{translation-principle}: if $M$ and $M'$ are two NFAs over the same category, and $F : \cQ \to \cQ'$ is a functional simulation between their underlying bare automata that preserves initial and accepting states, then $\Lang{M} \subseteq \Lang{M'}$, where equality holds $\Lang{M} = \Lang{M'}$ when $F$ is full, extending to an isomorphism of runs when $F$ is fully faithful.

In a different direction of abstraction, a span
\[
  \begin{tikzcd}
    & \cQ \ar[ld,"\text{finULF}"']\ar[rd,"O"]& \\
    \cC && \cD
  \end{tikzcd}
\]
of a finitary ULF functor $\cQ \to \cC$ and an arbitrary functor $O : \cQ \to \cD$ may be considered as a \defin{nondeterministic finite-state transducer} $\cC \spanmap \cD$, since it associates to every run $\alpha : q_0 \to q_f$ in $\cQ$ over an arrow $w : A \to B$ of $\cC$ a corresponding output arrow $O(\alpha) : O(q_0) \to O(q_f)$ in $\cD$.
By Propositions~\ref{proposition/ULF-fin-pullback-closed} and \ref{proposition/ULF-finitary-triangle}, ordinary composition of spans
\[
  \begin{tikzcd}
     &  & \cQ\times_{\cD}\cQ'\ar[ld,"\text{finULF}"']\ar[rd]\arrow[dd, phantom, "\rotatebox{-45}{$\lrcorner$}" very near start] && \\
    & \cQ \ar[ld,"\text{finULF}"']\ar[rd]&  & \cQ' \ar[ld,"\text{finULF}"']\ar[rd]&\\
    \cC && \cD && \cE
  \end{tikzcd}
\]
combines a transducer $\cC \spanmap \cD$ with a transducer $\cD \spanmap \cE$ to produce a transducer $\cC \spanmap \cE$, and one obtains a bicategory by defining a 2-cell between transducers of the same type as a functional simulation between the underlying bare automata that commutes with the output functors.

\subsection{Closure properties of regular languages}
\label{section/regular-closure}

Closure of finitary ULF functors under pullback along arbitrary functors (Proposition~\ref{proposition/ULF-fin-pullback-closed}) implies that regular languages of arrows and constants are closed under inverse functorial image, restricted to a suitable choice of objects in the source category or operad.
\begin{prop}\label{proposition/regular-invimage-closure-cat}
  Let $\Lang{} \subseteq \cC(A,B)$ be a regular language of arrows in $\cC$, let $F : \cD \to \cC$ be a functor of categories, and let $R,S$ be objects of $\cD$ such that $F(R) = A, F(S) = B$.
  Then the homset-restricted inverse image $F^{-1}(\Lang{}) \cap \cD(R,S)$ is regular.
\end{prop}
\begin{prop}\label{proposition/regular-invimage-closure-oper}
  Let $\Lang{} \subseteq \oO(A)$ be a regular language of constants in $\oO$, let $F : \oP \to \oO$ be a functor of operads, and let $R$ be a color of $\oP$ such that $F(R) = A$.
  Then the conset-restricted inverse image $F^{-1}(\Lang{}) \cap \oP(R)$ is regular.
\end{prop}
\noindent
Note that restricting the inverse image to a homset/conset is necessary due to our conventions for the initial and accepting states of categorical and operadic automata (although variations of those conventions could certainly be envisaged, see discussion above after Def.~\ref{definition/NFA-over-a-category}).
Indeed, given a categorical automaton $M = (\cC,\cQ,p,q_0,q_f)$, a functor $F : \cD \to \cC$, and a pair of objects $R,S$ of $\cD$ such that $F(R) = p(q_0), F(S) = p(q_f)$, the restricted pullback automaton is given by $\pull{F}M|_{R,S} = (\cD,\cD\times_{\cC}\cQ,\pull{F}p,(R,q_0),(S,q_f))$.
The restricted pullback of an operadic automaton is constructed similarly.

Conversely, since finitary ULF functors are closed under composition (Proposition~\ref{proposition/ULF-finitary-triangle}), we immediately obtain closure of regular languages under functorial image along a finitary ULF functor.
\begin{prop}\label{proposition/regular-finULF-image-closure-cat}
 If $\Lang{} \subseteq \cC(A,B)$ is a regular language of arrows in $\cC$ and $F : \cC \to \cD$ is a finitary ULF functor of categories, then the image $F(\Lang{}) \subseteq \cD(F(A),F(B))$ is also regular.  
\end{prop}
\begin{prop}\label{proposition/regular-finULF-image-closure-oper}
 If $\Lang{} \subseteq \oO(A)$ is a regular language of constants in $\oO$ and $F : \oO \to \oP$ is a finitary ULF functor of operads, then the image $F(\Lang{}) \subseteq \oP(F(A))$ is also regular.  
\end{prop}
\noindent
For example, the pushforward of a categorical automaton $M = (\cC,\cQ,p,q_0,q_f)$ along a finitary ULF functor $F : \cC \to \cD$ is the automaton $\push{F}M = (\cD,\cQ, F\circ p, q_0,q_f)$.

On the other hand, regular languages are \emph{not} closed under image along an arbitrary functor, as the following counterexamples demonstrate.
\begin{cexa}\label{counterex/non-finitary}
  Let $\ccat{W}$ be the \emph{category of Dyck walks} freely generated by the graph whose nodes are natural numbers and with a pair of edges 
  \[ \begin{tikzcd}[cells={nodes={font=\normalsize}}]
n \ar[r,"u_n",yshift=1.1ex] & n+1\ar[l,"d_n",yshift=-1.1ex]\end{tikzcd}\]
  for every $n \in \N$.
  Let $\Sigma$ be the alphabet containing a pair of brackets `$[$' and `$]$' and consider the functor $F : \ccat{W} \to \FreeCat{\Bouquet[\Sigma]}$ defined by
  \[ n \mapsto * \qquad u_n \mapsto \mathop{[}\qquad d_n \mapsto \mathop{]} \]
  for every $n\in\N$.
  The total language $\ccat{W}(0,0)$ is regular (Example~\ref{example/total-automaton}), but its image $F(\ccat{W}(0,0)) \subseteq \FreeCat{\Bouquet[\Sigma]}(*,*)$ is the Dyck language of well-bracketed words, which is not a regular language in the classical sense \cite[p.~195]{HMU2007} and hence cannot be a regular language of arrows in our sense by Prop.~\ref{proposition/ULF-nfa}.
  Observe in this case that the functor $F$ is ULF since it is generated from a graph homomorphism, but it is manifestly not finitary.
\end{cexa}
\begin{cexa}\label{counterex/non-ULF}
  Let $2 = \begin{tikzcd}[cells={nodes={font=\normalsize}}]0 \ar[r,"f"] & 1\end{tikzcd}$ be the walking arrow category, and consider any functor $F : 2 \to \cC$, which picks out an arrow $F(f) = w : A \to B$ in $\cC$.
  Observe that the functor $F$ is finitary, but will typically not be ULF unless $w$ is indecomposable.
  The singleton language $\set{f}$ is regular (Example~\ref{example/singleton}), and its image $F(\set{f}) = \set{w} \subseteq \cC(A,B)$ will be as well if the category $\cC$ has finitary factorizations, but otherwise need not be regular.
  In particular, let $\cC = \Monoid[\mathbb{R}_{\ge 0}]$ be the monoid of non-negative real numbers $(\mathbb{R}_{\ge 0},+,0)$ considered as a one-object category, and let $w = 1$ be the unit value considered as an arrow $1 : * \to *$.
  We claim that the language $\set{1} \subseteq \mathbb{R}_{\ge0}$ is not regular.
  Suppose by way of contradiction that there is an automaton $M$ recognizing $\set{1}$ equipped with a finitary ULF functor $p : \cQ \to \Monoid[\mathbb{R}_{\ge0}]$ and initial and accepting states $q_0$ and $q_f$.
  Since 1 factors in the monoid as
  \[
    1 = \frac{1}{n} + \frac{n-1}{n}
  \]
  for all $n\ge 1$, the accepting run $\alpha : q_0 \to q_f$ of 1 factors as
  \[
    \alpha = \begin{tikzcd}q_0 \ar[r,"\beta_n"] & q_n \ar[r,"\gamma_n"] & q_f \end{tikzcd}
  \]
  for some unique arrows $\beta_n$ and $\gamma_n$ lying over $\frac 1n$ and $\frac{n-1}n$, by assumption that $p$ is ULF. 
  Now there are two possibilities:
  \begin{enumerate}
  \item The objects $q_n$ are distinct for all $n \ge 1$. In this case the fiber $p^{-1}(*)$ is infinite, contradicting the assumption that $p$ is finitary.
  \item There exist $i < j$ such that $q_i = q_j$. In this case
  \[
    \alpha' = \begin{tikzcd}q_0 \ar[r,"\alpha_i"] & q_i \ar[r,equals] & q_j \ar[r,"\gamma_j"] & q_f \end{tikzcd}
  \]
  is another accepting run of the automaton over the element $p(\alpha') = \frac{1}{i} + \frac{j-1}{j} > 1$, contradicting the assumption that $M$ recognizes just the singleton language $\set{1}$.
  \end{enumerate}
\end{cexa}
\noindent
Finally, we can easily derive that regular languages are closed under intersection.
\begin{prop}
  If $\Lang{1},\Lang{2} \subseteq\cC(A,B)$ (resp.~$\Lang{1},\Lang{2} \subseteq \oO(A)$) are regular languages of arrows (resp.~constants), then so is their intersection $\Lang{1} \cap \Lang{2}$.
\end{prop}
\begin{proof}
  Given two automata $M_1$ and $M_2$ with underlying bare automata $p_1 : \cQ_1 \to \cC$ and $p_2 : \cQ_2 \to \cC$ over the same category $\cC$, and with initial and accepting states mapped to the same pair of objects $A$ and $B$, the intersection automaton $M_1 \cap M_2$ can be constructed in either of two equivalent ways:
  \begin{enumerate}
  \item By pulling back $p_1$ along $p_2$, and then taking the image along $p_2$.
  \item By forming the product automaton $M_1\times M_2$ (Example~\ref{example/product-automata}) over $\cC\times \cC$ and then pulling back along the diagonal functor $\cC \to \cC\times \cC$.
  \end{enumerate}
  The construction of the intersection of operadic automata is completely analogous.
\end{proof}
\noindent
At this point, we should mention that we expect regular languages of arrows and constants to satisfy some additional closure properties just as in the classical case, including suitably formulated versions of closure under union, concatenation, and Kleene star.
However, the usual proof of these facts is greatly facilitated by the use of $\epsilon$-transitions, which are ruled out by our working definition of categorical and operadic NFAs based on ULF functors (see Remark~\ref{remark/epsilon}).
We therefore defer consideration of these additional closure properties to the future introduction of an appropriate notion of categorical/operadic automaton with $\epsilon$-transitions.
In any case, we will not need them to prove the representation theorem.

\subsection{Related work}
\label{section/automata/related}

It is standard to represent finite-state automata by their transition graphs, and natural to consider such transition graphs as ordinary (finite directed) graphs equipped with a homomorphism into the bouquet graph on the underlying alphabet.
Such a perspective on nondeterministic finite-state automata without $\epsilon$-transitions was developed by Steinberg \cite{Steinberg2001}, motivated by monoid and semigroup theory, and with special prominence given to graph homomorphisms that are \emph{immersions,} which correspond to partial bideterministic automata in the sense of Definition~\ref{definition/determinism}. 
Walters' short note on context-free languages \cite{Walters1989} that we mentioned in \S\ref{section/cfls/related} also opens with a treatment of automata as graph homomorphisms, but taking morphisms of \emph{reflexive graphs} in order to allow $\epsilon$-transitions.
As discussed in Remark~\ref{remark/epsilon}, an important question is how to adapt our approach based on finitary ULF functors to allow for $\epsilon$-transitions, although we are confident that this can be done by appropriately generalizing the notion of ULF.
Walters' approach to automata and context-free grammars was further developed by Rosenthal~\cite{Rosenthal1995}, considering also tree automata, and more recently by Earnshaw and Soboci\'nski \cite{EarnshawSobocinski2022} to define regular languages in free monoidal categories.

The fibrational perspective on nondeterministic finite-state automata over categories and operads that we introduced here is closely related to the ``functorial approach'' to automata theory developed by Colcombet and Petri\c{s}an \cite{ColcombetPetrisan2020}.
Roughly speaking, the approaches are dual and related by Grothendieck-type correspondences: whereas we consider automata as functors into a category (or operad) with certain fibrational properties, Colcombet and Petri\c{s}an consider them as functors \emph{out of} a category into different large categories of interest, such as $\Set$ for deterministic automata, $\Rel$ for nondeterministic automata, or $\mathrm{Vec}$ for weighted automata.
Given that in many cases such Grothendieck-type correspondences rise to the level of equivalences, generally the approaches provide two complementary perspectives on the same objects.
Still, it is worth drawing attention to a few small but important differences between our approach and the approach described in \cite{ColcombetPetrisan2020}.
First, although many of the definitions and constructions given by Colcombet and Petri\c{s}an are stated for an arbitrary ``input category'' $\mathcal{I}$, in practice they consider automata over a fixed category $\mathcal{I}_{\mathsf{word}}$, isomorphic to the category we call $\FreeCat{\Bracket{\Bouquet[\Sigma]}}$, and they do not impose any finiteness condition on automata.
Second, rather than having a single initial and accepting state, Colcombet and Petri\c{s}an take a more abstract definition of the ``behavior'' of an automaton given by precomposing with a full and faithful functor $\iota : \mathcal{O} \to \mathcal{I}$ from some ``output category'' $\mathcal{O}$, typically the full subcategory of $\mathcal{I}_{\mathsf{word}} \cong \FreeCat{\Bracket{\Bouquet[\Sigma]}}$ spanned by the arrows $\bot \to \top$.
The idea is that precomposing with $\iota$ abstracts away the internal states of the automaton, remembering just the language that it recognizes.
Colcombet and Petri\c{s}an use this to give an elegant analysis of automata minimization based on factorization systems.
 \begin{wrapfigure}{l}{3cm}
  \begin{tikzpicture}[thick,scale=0.618, every node/.style={transform shape}]
    \node[state,initial] (q0) at (0,0) {$q_0$};
    \node[state,above right of=q0] (q1) {$q_1$};
    \node[state,below right of=q0] (q2) {$q_2$};
    \node[state,below right of=q1,accepting] (qf) {$q_f$};
    \draw
    (q0) edge node[above left]{$a$} (q1)
    (q0) edge node[below left]{$a$} (q2)
    (q1) edge node[above right]{$b$} (qf)
    (q2) edge node[below right]{$b$} (qf);
  \end{tikzpicture}
\end{wrapfigure}
Finally, and although it is by no means forced by the functorial approach, let us observe that modelling nondeterministic automata as functors into $\Rel$ (rather than as pseudofunctors into $\Span\Set$ or $\Span\FinSet$) means equating all of the runs over a given word between a given pair of states.
For instance, the NFA shown on the left has two distinct runs $q_0 \to q_f$ over the word $ab$, but these are identified when considering it as a $\Rel$-automaton.

The idea of generalizing regular languages from recognizable subsets of free monoids $\Lang{} \subseteq \Sigma^*$ to recognizable subsets of arbitrary monoids $\Lang{} \subseteq X$ is classic.
For instance, it is discussed in Eilenberg's book (see \cite[\S III.12, p.68]{Eilenberg1974}) using several of the equivalent characterizations of regular languages, including recognition by a finite monoid. 
He mentions that the definition of automaton could also be suitably generalized, but does not do so explicitly ``since all the known properties of recognizable sets already follow from the definition adopted above''.
There is also precedent for the idea of generalizing from languages of words in free monoids to languages of arrows in free categories, having been proposed at least as early as the mid-1980s in the work of Thérien and collaborators \cite{TherienSznajder1988,WeissTherien1986}, who argued that this was well-motivated by automata-theoretic considerations:
\begin{quote}
  For example, we are often interested in decompositions of automata.
  In such situations a component may receive its input from the output of some other component.
  This ``preprocessing'' imposes restrictions on the possible input sequences that need to be considered.
  A simple way to take into account these restrictions is to view a machine as processing input sequences that are paths in a finite directed multigraph. \cite[p.395]{TherienSznajder1988}
\end{quote}
Again, though, rather than taking NFAs or DFAs as a starting point, these authors started from the alternative characterization of regular word languages as unions of equivalence classes of a finite index congruence on a free monoid, and explained how to generalize it to free categories in a way that allowed proving Kleene's theorem among other results.
We take one more opportunity to emphasize that our definition of regular language is not restricted to free categories and operads, and that non-free categories and operads are required in some key examples.
(The work by Thérien et al.~was, incidentally, revisited by Jones who took a similar approach to define regular languages in a certain class of \emph{profinite categories} \cite{Jones1996profinite}.)
We do \emph{not} expect Kleene's theorem on the equivalence between regular and rational languages  to hold for regular languages of arrows in arbitrary categories -- indeed classically it is a theorem about recognizable subsets of free monoids with counterexamples in the non-free case \cite[\S IV.5, p.175--178]{Eilenberg1974}, and Counterexample~\ref{counterex/non-ULF} demonstrates that some singleton languages of arrows cannot be recognized by any NFA.
On the other hand, we mention that Fahrenberg et al.~have recently proved a Kleene theorem for higher-dimensional automata (without $\epsilon$-transitions) recognizing languages of interval pomsets \cite{FahrenbergJSZ2022}, which it would be interesting to analyze within our framework.

Although it is beyond the scope of this paper, we expect that regular languages of arrows or constants in the sense of Definitions~\ref{definition/NFA-over-a-category} and \ref{definition/NFA-over-an-operad} should also have other equivalent characterizations related to the aforementioned classical ones -- in particular there seem to be natural notions of \emph{recognition by a locally finite category} and \emph{recognition by a locally finite operad} (cf.~Remark~\ref{remark/DFA-over-spliced-arrows}).
It may be fruitful to connect this with Salvati's notion of regular language of $\lambda$-terms \cite{Salvati09} (in part inspired by \cite{deGroote2001acgs}), which may be formulated as recognition by a locally finite and well-pointed cartesian closed category (cf.~\cite{MoreauNguyen24}).

Finally, we only briefly touched upon the subject of constructing categories of automata and transducers in \S\ref{section/properties-ULF-finitary}, but this also has precedent in the literature.
Our definition of a nondeterministic finite-state transducer as a span $\cC \twoheadleftarrow \cQ \rightarrow \cD$ where the left leg satisfies a fibrational condition (namely, is finitary and ULF) is similar to the notion of ``Mealy morphism'' defined by Paré using slightly different conditions \cite{Pare2012mealy}. 
Another line of work by Katis, Sabadini, and Walters \cite{KatisSabadiniWalters1997} models automata as spans of graphs, with sequential composition given by composition of spans.
This graph-theoretic approach is certainly related to ours, since spans of graphs are in one-to-one correspondence with spans of ULF functors between the associated free categories (Prop.~\ref{ULF-into-free}), and moreover the two notions of composition agree.
Hyland defined a \emph{traced monoidal category} $\mathrm{Aut}$ whose objects are natural numbers and whose morphisms $m \to n$ are finite-state automata (over a fixed alphabet $\Sigma$) equipped with $m$ distinct initial states and $n$ distinct accepting states \cite{Hyland2008}.
He then exhibited a traced monoidal functor $\mathrm{Aut} \to \mathrm{Rat}$ to a category whose arrows are matrices of rational languages, as one direction of Kleene's theorem.
It could be interesting to adapt Hyland's analysis to our setting.

\section{The Chomsky-Schützenberger Representation Theorem}
\label{section/csrep}

In this section, we give a relatively simple and conceptual proof of the representation theorem, generalized to context-free languages of arrows in an arbitrary category.
We begin in \S\ref{section/pulling-back} by explaining how to represent the intersection of a context-free language with a regular language by first taking the \emph{pullback} of a CFG along an NFA.
Next, in \S\ref{section/contour-category}, we exhibit a left adjoint to the spliced arrow operad construction, which we call the \emph{contour category} construction.
As a direct consequence of this adjunction, we show in \S\ref{section/unigram} that every pointed finite species induces a CFG that is universal in a precise sense and that generates a language of \emph{tree contour words,} which are closely related to Dyck words.
Finally, in \S\ref{section/represenation-theorem} we state and prove an appropriate generalization of the representation theorem, which also relies in a crucial way on the strong translation principle.

\subsection{Pulling back context-free grammars along finite-state automata}\label{section/pulling-back}

To compute the pullback, we apply two lemmas.
\begin{lem}\label{lemma/pullback-free-along-ULF}
Suppose given a species~$\Sspecies$, a functor of operads ${p:\FreeOper{\Sspecies}\to\Ooperad}$
and a ULF functor of operads $p_{\Qoperad}:\Qoperad\to\Ooperad$.
In that case, 
the pullback of $p$ along $p_{\Qoperad}$ in the category of operads is obtained from
a corresponding pullback of $\phi : \Sspecies \to \ForgetOper{\Ooperad}$ along $\ForgetOper{p_{\Qoperad}} : \ForgetOper{\Qoperad}\to\ForgetOper{\Ooperad}$ in the category of species:
\begin{equation*}
\begin{tikzcd}[column sep=4em,row sep=2em]
  {\Sspecies'}\arrow[d,"{\phi' = \pull{p_{\Qoperad}}\phi}"{swap}]\arrow[r,"\psi'"]
  \arrow[rd, phantom, "\lrcorner" very near start]
  & {\Sspecies}\arrow[d,"\phi"]
\\
{\ForgetOper{\Qoperad}}\arrow[r,"{\ForgetOper{p_{\Qoperad}}}"{swap}] & {\ForgetOper{\Ooperad}}
\end{tikzcd}\ \ \in\ \ \Species{}
\qquad\vrule\ \vrule\qquad
\begin{tikzcd}[column sep=4em,row sep=2em]
  \FreeOper{\Sspecies'}\arrow[d,"{p' = \pull{p_{\Qoperad}}p}"{swap}]\arrow[r,"\FreeOper{\psi'}"]
  \arrow[rd, phantom, "\lrcorner" very near start]
  & \FreeOper{\Sspecies}\arrow[d,"p"]
\\
\Qoperad\arrow[r,"{p_{\Qoperad}}"{swap}] & \Ooperad
\end{tikzcd}\ \ \in\ \ \Operads{}
\end{equation*}
\end{lem}
\begin{proof}
  By Proposition~\ref{proposition/ULF-fin-pullback-closed}, the pullback of the ULF functor $p_{\Qoperad}$ along the functor of operads $p$ is a ULF functor $\Qoperad' \to \FreeOper{\Sspecies}$, which must be of the form $\FreeOper{\psi'} : \FreeOper{\Sspecies'} \to \FreeOper{\Sspecies}$ for some map of species $\psi' : \Sspecies' \to \Sspecies$ by Proposition~\ref{ULF-into-free-operad}.
We therefore have a pullback in $\Operads{}$ as on the right above, which is sent to a pullback in $\Species{}$ by the right adjoint forgetful functor.
The unit of the adjunction induces another pullback diagram
\[
\begin{tikzcd}[column sep=4em,row sep=2em]
  \Sspecies'\arrow[r,"{\psi'}"]\arrow[d,"{\eta_{\Sspecies}}"']
  \arrow[rd, phantom, "\lrcorner" very near start]
  & \Sspecies\arrow[d,"\eta_{\Sspecies}"] \\
  \FreeOper{\Sspecies'}\arrow[r,"\FreeOper{\psi'}"] & \FreeOper{\Sspecies}
\end{tikzcd}
\]
and by composing these
\[
\begin{tikzcd}[column sep=4em,row sep=2em]
  \Sspecies'\arrow[r,"{\psi'}"]\arrow[d,"{\eta_{\Sspecies}}"']
  \arrow[rd, phantom, "\lrcorner" very near start]
  & \Sspecies\arrow[d,"\eta_{\Sspecies}"] \\
  \FreeOper{\Sspecies'}\arrow[d,"{p'}"{swap}]\arrow[r,"\FreeOper{\psi'}"]
  \arrow[rd, phantom, "\lrcorner" very near start]
  & \FreeOper{\Sspecies}\arrow[d,"p"]
\\
\Qoperad\arrow[r,"{p_{\Qoperad}}"{swap}] & \Ooperad
\end{tikzcd}  \quad=\quad
\begin{tikzcd}[column sep=4em,row sep=2em]
  {\Sspecies'}\arrow[d,"{\phi'}"{swap}]\arrow[r,"\psi'"]
  \arrow[rd, phantom, "\lrcorner" very near start]
  & {\Sspecies}\arrow[d,"\phi"]
\\
{\ForgetOper{\Qoperad}}\arrow[r,"{\ForgetOper{p_{\Qoperad}}}"{swap}] & {\ForgetOper{\Ooperad}}
\end{tikzcd}
\]
we obtain the desired pullback in the category of species.
\end{proof}
\noindent
\begin{lem}\label{lemma/pullback-finite-species-along-finitary}
  Consider a pullback in the category of species:
\[
\begin{tikzcd}[column sep=4em,row sep=2em]
  \Sspecies'\arrow[r,"{\psi'}"]\arrow[d,"\phi'"']
  \arrow[rd, phantom, "\lrcorner" very near start]
  & \Sspecies\arrow[d,"\phi"] \\
  \Tspecies'\arrow[r,"\psi"'] & \Tspecies
\end{tikzcd}
\]
If $\Sspecies$ is finite and $\psi$ is finitary (i.e., has finite fibers) then $\Sspecies'$ is finite.
\end{lem}
\begin{proof}
The pullback species $\Sspecies' = \Sspecies \times_{\Tspecies} \Tspecies'$ admits a simple description: its colors (resp.~nodes) are pairs of colors (resp.~nodes) from $\Sspecies$ and $\Tspecies'$ lying over the same color (resp.~node) in $\Tspecies$.
Hence finiteness of $\Sspecies'$ follows immediately from the assumption that $\Sspecies$ is finite and that $\psi$ has finite fibers.
\end{proof}
\begin{thm}\label{theorem/pullback-of-cfg-along-nfa}
For any CFG $G=(\cC,\Sspecies,S,p_G)$ and NFA $M = (\cC,\cQ,p_M,q_0,q_f)$ over the same category, with $p_G(S) = (p_M(q_0),p_M(q_f))$, there is a CFG $G' = \pull{M}G$ generating the language $\Lang{G'} = p_M^{-1}(\Lang{G}) \cap \cQ(q_0,q_f)$.
\end{thm}
\begin{proof}
By Prop.~\ref{proposition/Words-ULF}, the finitary ULF functor of categories~$p_M : \cQ \to \cC$ induces a finitary ULF functor of operads ${\Words{p_M}} : \Words{\Qcategory}\to \Words{\Ccategory}$.
By Lemma~\ref{lemma/pullback-free-along-ULF}, the pullback of $p_G$ along $p_M$ may therefore be computed as a pullback of the underlying maps of species:
\[
\begin{tikzcd}[column sep=4em,row sep=2em]
  {\Sspecies'}\arrow[d,"{\phi_G'}"{swap}]\arrow[r,"\psi'"]
  \arrow[rd, phantom, "\lrcorner" very near start]
  & {\Sspecies}\arrow[d,"\phi_G"]
\\
{\ForgetOper{\Words{\cQ}}}\arrow[r,"{\ForgetOper{\Words{p_M}}}"{swap}] & {\ForgetOper{\Words{\cC}}}
\end{tikzcd}\ \ \in\ \ \Species{}
\qquad\vrule\ \vrule\qquad
\begin{tikzcd}[column sep=4em,row sep=2em]
  \FreeOper{\Sspecies'}\arrow[d,"{p_G'}"{swap}]\arrow[r,"\FreeOper{\psi'}"]
  \arrow[rd, phantom, "\lrcorner" very near start]
  & \FreeOper{\Sspecies}\arrow[d,"p_G"]
\\
\Words{\cQ}\arrow[r,"{\Words{p_M}}"{swap}] & \Words{\cC}
\end{tikzcd}\ \ \in\ \ \Operads{}
\]
Moreover, by Lemma~\ref{lemma/pullback-finite-species-along-finitary}, the pullback species $\Sspecies' = \Sspecies\times_{\Words{\cC}}\Words{\cQ}$ is finite,
so that $\pull{M}G = (\cQ,\Sspecies',(S,(q_0,q_f)),p'_G)$ defines a CFG.
Finally, we have $\Lang{G'} = p_M^{-1}(\Lang{G}) \cap \cQ(q_0,q_f)$ because by the universal property of the pullback in $\Operads{}$, $G'$ derives a run $\alpha : q_0 \to q_f$ of $M$ just in case there exists a constant $\beta : S$ in $\FreeOper{\Sspecies}$ such that $p_G(\beta) = p_M(\alpha)$.
\end{proof}
\noindent
We emphasize that the pullback grammar $\pull{M}G$ generates a language of arrows in $\cQ$, corresponding to the runs of the automaton~$M$ over the arrows in the language generated by the original grammar $G$.
The pullback grammar admits the following concrete description in traditional CFG syntax:
\begin{itemize}
\item its nonterminals are pairs $(R,(q,q'))$ -- which for geometric intuition we prefer to visualize as triples $(q,R,q')$ -- for every nonterminal $R$ of the grammar $G$ and pair of states $q,q'$ of the automaton $M$ such that $p_G(R) = (p_M(q),p_M(q'))$; 
\item it has a production rule $(q,R,q')\to\alpha_0(q_1,R_1,q'_1)\alpha_1\dots (q_n,R_n,q'_n)\alpha_n$
for every production rule $R\to w_0R_1w_1\dots R_n w_n$ of $G$
and sequence of $n+1$ runs $\alpha_0:q\to q_1$, $\alpha_1:q'_1\to q_2$, $\dots$, $\alpha_n:q'_{n}\to q'$ of $M$ over the respective arrows $w_0,\dots,w_n$;
\item the start symbol is $(q_0,S,q_f)$, where $S$ is the start symbol of $G$ and $q_0$ and $q_f$ are the initial and accepting states of $M$.
\end{itemize}
\begin{exa}\label{example/parsing-intervals}
  Recall from Example~\ref{equation/nfa-Mw} that to any arrow $w : A \to B$ of a category with finitary factorizations $\Ccategory$, there is associated an NFA $M_w$ recognizing exactly the singleton language $\set{w}$.
  If $G$ is a CFG over $\Ccategory$, then by pulling back $G$ along $M_w$ we obtain a new grammar $\pull{M_w}G$ that may be seen as a specialization of $G$ to the factors of $w$.
  In particular, in the classical case where $w = a_1\dots a_n$ is a word of length $n$ over $\Sigma$, the nonterminals of $\pull{M_w}G$ may be taken as triples $(i,R,j)$ of a nonterminal $R$ of $G$ together with a pair of indices $0\le i,j \le n$, generating all the parses of the subword $w_{i,j} = a_{i+1}\dots a_j$ as an $R$.
  If we then view this grammar as defining a displayed operad over $\Words{\Fact{w}}$, as discussed in \S\ref{section/benabou}--\S\ref{section/magic-formula}, we obtain precisely the classical parse matrices $N_{i,j}$ referenced in \S\ref{section/application-to-parsing}.
\end{exa}
As an immediate corollary of Theorem~\ref{theorem/pullback-of-cfg-along-nfa}, we obtain the closure of CFLs of arrows under intersection with regular languages of arrows.
\begin{cor}\label{corollary/intersection-of-languages}
  If $\Lang{} \subseteq \cC(A,B)$ is context-free and $\Lang{M} \subseteq \cC(A,B)$ is regular then $\Lang{} \cap \Lang{M} \subseteq \cC(A,B)$ is context-free.
\end{cor}
\begin{proof}
  By taking the functorial image (Prop.~\ref{proposition/basic-closure}(3)) of the CFL $p_M^{-1}(\Lang{G}) \cap \cQ(q_0,q_f)$ along the functor $p_M$.
\end{proof}
\noindent
Explicitly, the grammar generating $L \cap L_M$ is defined by starting with the pullback grammar $G' = \pull{M}G$ from Theorem~\ref{theorem/pullback-of-cfg-along-nfa} above and postcomposing the functor $p'_G : \FreeOper{\Sspecies'} \to \Words{\cQ}$ with $\Words{p_M} : \Words{\cQ} \to \Words{\cC}$.
This grammar admits almost the same concrete description as $G'$ but where the production rules instead take the form
\[(q,R,q')\to w_0(q_1,R_1,q'_1) w_1\dots (q_n,R_n,q'_n)w_n\]
for every production rule $R\to w_0R_1w_1\dots R_n w_n$ of $G$
and sequence of $n+1$ runs $q\to q_1$, $q'_1\to q_2$, $\dots$, $q'_{n}\to q'$ of $M$ over the respective arrows $w_0,\dots,w_n$.
The reader familiar with the classic ``Bar-Hillel construction'' for combining a CFG with an NFA \cite[Theorem~8.1]{BarHillel+1961} may recognize that our construction is very similar, but completely uniform and with the benefit of being derived systematically.

Although it will not be relevant to the representation theorem for categorical CFGs, let us observe that the construction of the pullback of a CFG along a categorical NFA works just as well to construct the pullback of a generalized CFG (in the sense of \S\ref{section/gCFGs}) along an operadic NFA.
\begin{thm}
For any gCFG $G=(\oO,\Sspecies,S,p_G)$ and operadic NFA $M = (\oO,\oQ,p_M,q_r)$ over the same operad, with $p_G(S) = p_M(q_r)$, there is a gCFG $G' = \pull{M}G$ generating the language $\Lang{G'} = p_M^{-1}(\Lang{G}) \cap \oQ(q_r)$.
\end{thm}
\noindent
Indeed, the proof is almost identical to the proof of Theorem~\ref{theorem/pullback-of-cfg-along-nfa}, simply omitting the first step where we applied the functor $\Wordsonly$ before applying Lemma~\ref{lemma/pullback-free-along-ULF}.
We immediately obtain that gCFLs are closed under intersection with regular languages of constants.
\begin{cor}\label{corollary/intersection-gCFL-regular}
  If $\Lang{} \subseteq \oO(A)$ is a gCFL and $\Lang{M} \subseteq \oO(A)$ is regular then $L\cap \Lang{M} \subseteq \oO(A)$ is a gCFL.
\end{cor}

\subsection{The contour category of an operad and the contour / splicing adjunction}
\label{section/contour-category}

In \S\ref{section/spliced-arrows}, we began by explaining how to construct a functor
$\Wordsonly : \Cat \to \Operads{}$
transforming any category $\Ccategory$ into an operad $\Words{\Ccategory}$ of spliced arrows of arbitrary arity, which played a central role in our definition of context-free language of arrows in a category.
We construct now a left adjoint functor
\begin{equation}\label{equation/adjunction-contour-splice}
\begin{tikzcd}
\Operads{}
\arrow[rr,"{\Contouronly}",yshift=1.5ex]
& \bot &
\Cat\arrow[ll,"{\Wordsonly}",yshift=-1.5ex]
\end{tikzcd}
\end{equation}
which extracts from any given operad $\Ooperad$ a category $\Contour{\Ooperad}$ whose arrows may be interpreted as ``oriented contours'' along the boundary of the operations of the operad.

\begin{figure}
  \begin{center}\includegraphics[width=\textwidth]{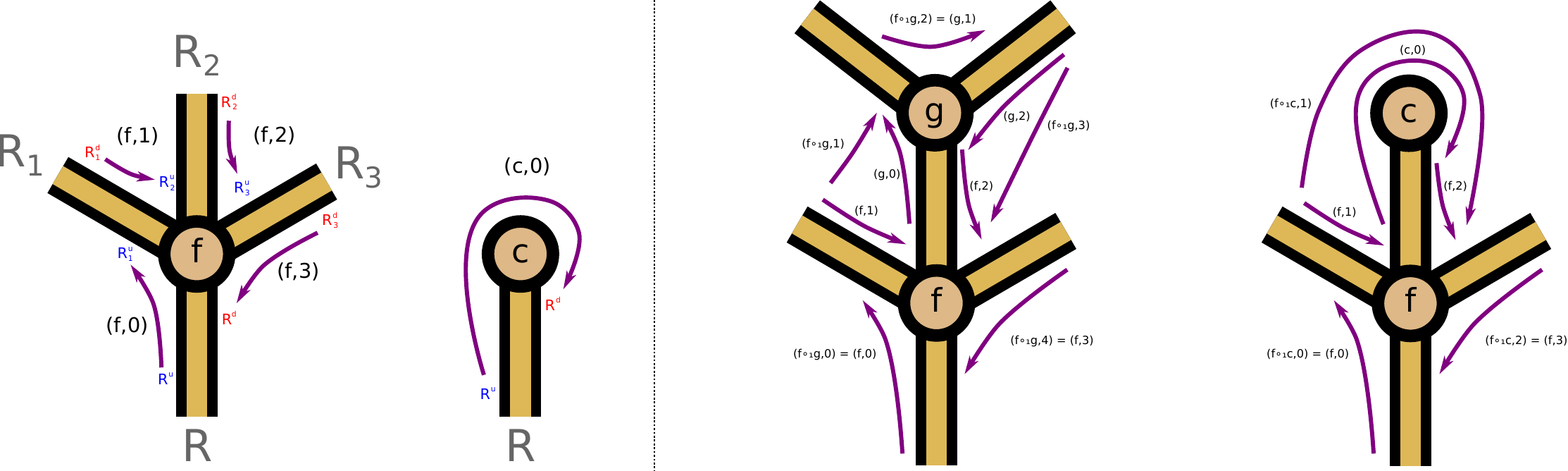}\end{center}
  \caption{Left: interpretation of the generating arrows of the contour category $\Contour{\Ooperad}$. Right: interpretation of equations \eqref{equation/contour1} and \eqref{equation/contour2}.}
  \label{fig:contour-category}
\end{figure}

\begin{defi}\label{definition/contour-category}
The \defin{contour category $\Contour{\Ooperad}$} of an operad $\Ooperad$ is defined as a quotient of the following free category:
\begin{itemize}
\item objects are given by \emph{oriented colors} $R^\epsilon$ consisting of a color $R$ of $\Ooperad$ and an orientation $\epsilon \in \set{\tagU,\tagD}$ (``up'' or ``down'');
\item arrows are generated by pairs $(f,i)$ of an operation $f : R_1,\dots,R_n \to R$ of $\Ooperad$ and an index $0 \le i \le n$, defining an arrow $R_i^\tagD \to R_{i+1}^\tagU$ under the conventions that $R_0^\tagD = R^\tagU$ and $R_{n+1}^\tagU = R^\tagD$;
\end{itemize}
subject to the conditions that $\id[R^\tagU] = (\id[R], 0)$ and $\id[R^\tagD] = (\id[R], 1)$ as well as the following equations:
\begin{align}
  (f \circ_i g, j)  &= \begin{cases}
                         (f,j) & \text{if }j < i  \\ 
                         (f,i)(g,0)\hphantom{(f,i+1)} & \text{if }j = i \\
                         (g,j-i) & \text{if }i < j < i+m \\
                         (g,m)(f,i+1) & \text{if }j = i+m \\
                         (f,j-m+1) & \text{if }j > i+m
                       \end{cases} \label{equation/contour1}\\ 
  (f \circ_i c, j) &= \begin{cases}
                         (f,j) & \text{if }j < i \\
                         (f,i)(c,0)(f,i+1) & \text{if }j = i \\
                         (f,j+1) & \text{if }j > i
                       \end{cases} \label{equation/contour2}
\end{align}
whenever the left-hand side is well-formed,
for every operation $f$, operation $g$ of positive arity $m > 0$, constant $c$,
and indices $i$ and $j$ in the appropriate range.

We refer to each generating arrow $(f,i)$ of the contour category $\Contour{\Ooperad}$ as a \defin{sector} of the operation~$f$.
See Fig.~\ref{fig:contour-category} for a graphical interpretation of sectors and of the equations on contours seen as compositions of sectors.
\end{defi}

The contour construction provides a left adjoint to the spliced arrow construction
because a functor of operads
$\Ooperad \to \Words{\Ccategory}$
is entirely described by the data of a pair of objects $(A,B) = (R^\tagU, R^\tagD)$ in $\Ccategory$ for every color $R$ in $\Ooperad$ together with a sequence $f_0, f_1, \dots, f_n$ of $n+1$ arrows in $\Ccategory$, where $f_i : R_i^\tagD \to R_{i+1}^\tagU$ for $0 \le i \le n$ for each operation $f : R_1,\dots,R_n \to R$ of $\Ooperad$, under the same conventions as above.
The equations \eqref{equation/contour1} and \eqref{equation/contour2} on the generators of $\Contour{\Ooperad}$ reflect the equations imposed by the functor of operads $\Ooperad \to \Words{\Ccategory}$ on the spliced arrows of $\Ccategory$ appearing as the image of operations in $\Ooperad$.
In that way we transform any functor of operads $\Ooperad \to \Words{\Ccategory}$ into a functor
$\Contour{\Ooperad} \to \Ccategory$
which may be seen as an interpretion of the contours of the operations of $\Ooperad$ in $\Ccategory$.

The unit and counit of the contour / splicing adjunction also have nice descriptions.
The unit of the adjunction defines, for any operad $\Ooperad$, a functor of operads $\Ooperad \to \Words{\Contour{\Ooperad}}$ that acts on colors by $R \mapsto (R^\tagU, R^\tagD)$,
and on operations by sending an operation $f : R_1,\dots,R_n \to R$ of $\Ooperad$ to the spliced word of sectors $(f,0){-}\dots{-}(f,n) : (R_1^\tagU,R_1^\tagD),\dots,(R_n^\tagU,R_n^\tagD) \to (R^\tagU,R^\tagD)$.
The counit of the adjunction defines, for any category $\Ccategory$, a functor of categories $\Contour{\Words{\Ccategory}} \to \Ccategory$ that acts on objects by $(A,B)^\tagU \mapsto A$ and $(A,B)^\tagD \mapsto B$, and on arrows by sending the $i$th sector of a spliced word to its $i$th word, $(w_0{-}\dots{-}w_n,i) \mapsto w_i$.

\begin{rem}\label{remark/contour-free-is-free}
In the case of a free operad over a species $\sS$,
the contour category $\Contour{\FreeOper{\sS}}$
admits an even simpler description as a free category $\Contour{\FreeOper{\sS}} \cong \FreeCat{\Contour{\sS}}$ generated by the arrows $(x,i) : R_i^\tagD \to R_{i+1}^\tagU$ for every node $x : R_1,\dots,R_n \to R$ of the species $\Sspecies$.
Indeed, it is easy to verify that there is an analogous adjunction between species and graphs, and that this pair of contour / splicing adjunctions commutes with the pair of respective free / forgetful adjunctions:
\begin{equation}\label{equation/adjunction-contour-splice-extended}
\begin{tikzcd}[row sep=3em]
\Operads{}
\arrow[rr,"{\Contouronly}",yshift=1.5ex]\arrow[d,xshift=1.5ex]
& \bot &
\Cat\arrow[ll,"{\Wordsonly}",yshift=-1.5ex]\arrow[d,xshift=1.5ex] \\
\Species{}
\arrow[rr,"{\Contouronly}",yshift=1.5ex]\arrow[u,xshift=-1.5ex]\arrow[u,phantom,"\dashv"]
& \bot &
\Graph\arrow[ll,"{\Wordsonly}",yshift=-1.5ex]\arrow[u,xshift=-1.5ex]\arrow[u,phantom,"\dashv"]
\end{tikzcd}
\end{equation}
We refer to the generating arrows $(x,i)$ of the contour graph $\Contour{\Sspecies}$ as \defin{corners} in the sense of the theory of planar maps \cite{Schaeffer2015map}, since geometrically they indeed correspond to the corners of $\Sspecies$-rooted trees seen as rooted planar maps.
By \eqref{equation/adjunction-contour-splice-extended}, every sector of an operation of $\FreeOper{\Sspecies}$ factors uniquely in the contour category~$\Contour{\FreeOper{\Sspecies}}$ as a sequence of corners.
\end{rem}
\noindent
In contrast to the situation for $\Wordsonly$ (Prop.~\ref{proposition/Words-ULF}), it is \emph{not} the case that $\Contouronly$ always preserves the ULF property.
\begin{rem}
Consider the category $\mathbf{1}+\mathbf{1}$ with two objects $A$ and $B$ and only identity arrows,
and the unique functor $p$ to the terminal category $\mathbf{1}$.
We claim that the associated ULF functor of operads $\Words{p}$ induces a functor of categories $\Contour{\Words{p}}$ which is not ULF.
%
Consider the two binary operations $f = \id[A]{-}\id[A]{-}\id[A]$ and $g = \id[A]{-}\id[A]{-}\id[B]$
and the constant $c = \id[A]$ in $\Words{(\mathbf{1}+\mathbf{1})}$,
as well as the binary operations $h=\id[\ast]{-}\id[\ast]{-}\id[\ast]$ and the constant $d=\id[\ast]$
in $\Words{\mathbf{1}}$.
The category~$\Contour{\Words{(\mathbf{1}+\mathbf{1})}}$ has the sequence of sectors $\alpha = (f,0)(c,0)(g,1)$ as an arrow, which
is different from the identity. 
On the other hand, it is mapped by $\Contour{\Words{p}}$ to the sequence $w = (h,0)(d,0)(h,1)$, which is equal thanks to Equation~\eqref{equation/contour2}
to the sector $(h\circ_0 d,0)$ of the unary operation~$h\circ_0 d = \id[\ast]{-}\id[\ast]$ of $\Words{\mathbf{1}}$,
and hence $w = \id[(\ast,\ast)^\tagU]$.
Since the factorization $\id = \id\, \id$ in $\Words{\mathbf{1}}$ lifts to two distinct factorizations $\alpha = \id\, \alpha = \alpha\, \id$ in $\Words{(\mathbf{1}+\mathbf{1})}$, $p$ is not ULF.
\end{rem}
\noindent
Still, we can verify that maps of species induce ULF functors between their contour categories, as an immediate consequence of Prop.~\ref{ULF-into-free} and the factorization \eqref{equation/adjunction-contour-splice-extended}.
\begin{prop}\label{proposition/Contour-species-ULF}
If $\psi:\Sspecies\to \Tspecies$ is a map of species, then $\Contour{\FreeOper{\psi}} : \Contour{\FreeOper{\Sspecies}} \to \Contour{\FreeOper{\Tspecies}}$ is a ULF functor of categories.
\end{prop}

\subsection{The universal context-free grammar of a pointed species, and its associated tree contour language}
\label{section/unigram}

\begin{figure}
  \begin{center}\includegraphics[width=0.35\textwidth]{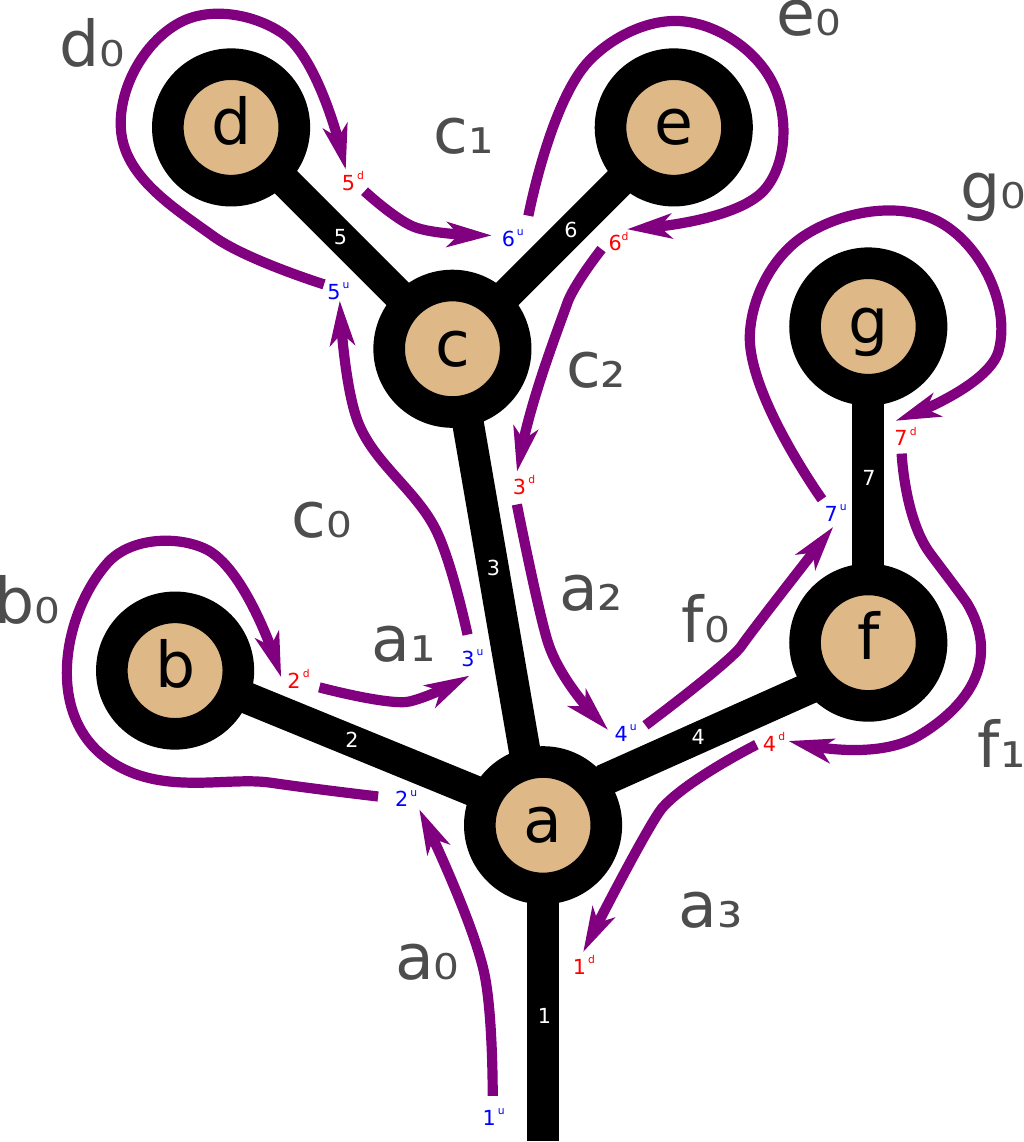} \qquad \vrule \qquad \includegraphics[width=0.35\textwidth]{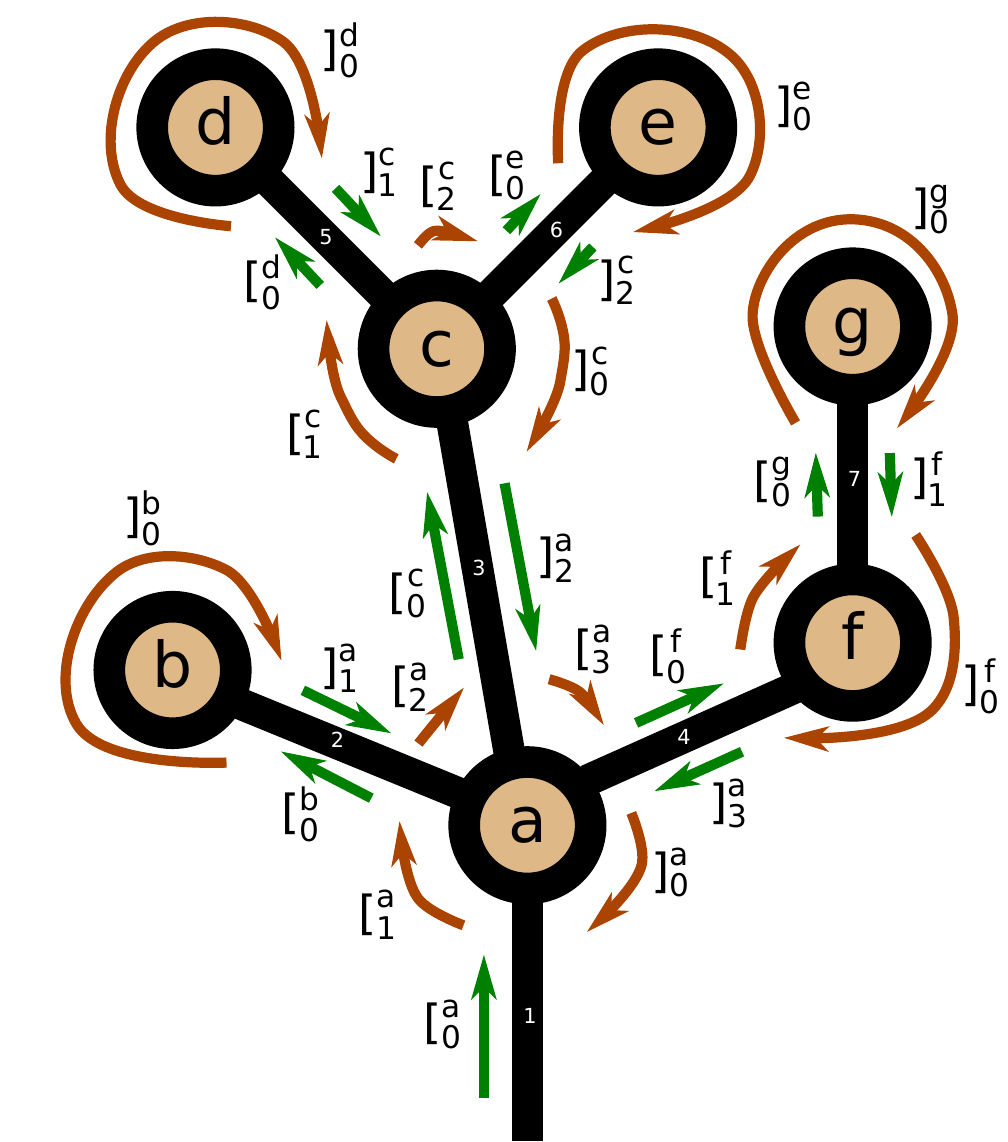} \end{center}
  \caption{Left: an $\Sspecies$-rooted tree of root color 1 and its corresponding contour word $\mathsf{a_0 b_0 a_1 c_0 d_0 c_1 e_0 c_2 a_2 f_0 g_0 f_1 a_3} : 1^\tagU \to 1^\tagD$.
    Right: the corresponding Dyck word obtained by splitting each corner of the contour in two pieces. 
  }
  \label{fig:contour-word}
\end{figure}

Every finite species~$\Sspecies$ equipped with a color~$S$ 
comes with a \emph{universal} context-free grammar~$(\Contour{\FreeOper{\Sspecies}},\Sspecies,S,p_\Sspecies)$, characterized by the fact that $p_\sS : {\FreeOper{\Sspecies}}\to{\Words{\Contour{\FreeOper{\Sspecies}}}}$ is the unit of the contour / splicing conjunction.
We denote the universal grammar associated to a pointed species $\dot\sS := (\sS,S)$ by $\UnivGrammar{\dot\sS}$.
By ``universal'' context-free grammar, we mean that any categorical context-free grammar
$G=(\Ccategory,\Sspecies,S,p)$ with the same underlying species and start symbol factors uniquely through~$\UnivGrammar{\dot\sS}$ in the sense
that there exists a unique functor $q_G:\Contour{\FreeOper{\Sspecies}}\to\Ccategory$ satisfying the equation
\begin{equation}\label{equation/universal-CFG-factorization}
\begin{tikzcd}[column sep=2em]
{\FreeOper{\Sspecies}}
\arrow[rr,"{p}"]
&&
\Words{\Ccategory}
\quad = \quad
{\FreeOper{\Sspecies}}
\arrow[rr,"{p_{\Sspecies}}"] 
&& 
{\Words{\Contour{\FreeOper{\Sspecies}}}}
\arrow[rr,"{\Words{q_G}}"] && \Words{\Ccategory}
\end{tikzcd}
\end{equation}
We refer to the language of arrows of a universal grammar $\UnivGrammar{\dot\sS}$
as the \defin{tree contour language} associated to a pointed species, denoted $\ContourLang{\dot\sS}$, and to its arrows as \defin{tree contour words}, or more precisely as $\dot\sS$-tree contour words.
Tree contour words in $\ContourLang{\dot\sS}$ describe the contours of $\Sspecies$-rooted trees with root color $S$, see left side of Fig.~\ref{fig:contour-word} for an illustration.
The factorization above shows that any context-free grammar $G$ is the functorial image of the universal grammar $\UnivGrammar{\dot\sS}$ along the functor of categories~$q_G$, whose purpose is to transport each corner of a node in $\Sspecies$ to the corresponding arrow in~$\Ccategory$ as determined by $G$.
At the level of languages, we have 
\begin{prop}\label{proposition/image-of-a-contour-language}
The language of a context-free grammar $G$ coincides with the image of its tree contour language, $\Lang{G} = q_G\,\ContourLang{\dot\sS}$.
\end{prop}
\begin{exa}
  Consider the context-free grammar of Figure~\ref{fig:example-cfg}.
  The production rules of the universal grammar $\UnivGrammar{\dot\sS}$ may be depicted in traditional CFG syntax as
  \begin{align*}
    \mathsf{S} &\to 1_0\,\mathsf{NP}\,1_1\,\mathsf{VP}\,1_2 \\
    \mathsf{NP} &\to 2_0 \\
    \mathsf{NP} &\to 3_0 \\
    \mathsf{VP} &\to 4_0\,\mathsf{NP}\,4_1
  \end{align*}
  where we have written $i_j$ as shorthand for the sector $(x_i,j)$.
  The language $\Lang{G}$ is then obtained as the image of $\ContourLang{\dot\Sspecies}$ under the functor $q_G : \Contour{\FreeOper\Sspecies} \to \FreeCat{\Bouquet[\Sigma]}$ defined by:
  \[
    1_0 \mapsto \id \quad
    1_1 \mapsto \text{\textvisiblespace} \quad
    1_2 \mapsto \id
  \]
  \[
    2_0 \mapsto \mathsf{mom}\qquad
    3_0 \mapsto \mathsf{tom} \qquad
    4_0 \mapsto \mathsf{loves}\text{\textvisiblespace}\quad
    4_1 \mapsto \id
  \]
\end{exa}

\begin{rem}
  The notion of tree contour language makes sense even for infinite species, although in that case the resulting universal grammar $\UnivGrammar{\dot\sS}$ is no longer a CFG, having infinitely many nonterminals or infinitely many productions.
  Still, it may be an interesting object of study.
  In particular, the tree contour language $\UnivGrammar{(\TerminalSpecies,\ast)}$ generated by the \emph{terminal species} $\TerminalSpecies$ with one color $\ast$ and a single operation of every arity appears to be of combinatorial interest, with words in the language describing the shapes of rooted planar trees with arbitrary node degrees.
  Nevertheless, from now on when we refer to a ``tree contour language'' we will always mean the context-free tree contour language generated by a finite pointed species.
\end{rem}

\subsection{Representation theorem}
\label{section/represenation-theorem}

The achievement of the classical Chomsky-Schützenberger representation theorem is to separate any
context-free grammar $G=(\Sigma,\Nonterminals,\Sentence,\Productions)$
into two independent components:
\begin{enumerate}
\item
  a context-free grammar $G'$ with only one nonterminal over an alphabet
  \[\Sigma_{2k} = \set{ [_1, ]_1, \dots, [_k, ]_k}\]
  of size $2k$, for some $k$, which generates Dyck words of balanced brackets describing the shapes of parse trees with nodes labelled by
production rules of $G$; and
\item
  a finite-state automaton $M$ to check that the edges of these trees may be appropriately colored by
  the nonterminals of $G$ according to the labels of the nodes specifying the productions.
\end{enumerate}
The original context-free language generated by $G$ is then obtained as the image of the intersection of the
Dyck language generated by $G'$ with the regular language recognized by $M$, under a homomorphism $\Sigma_{2k}^* \to \Sigma^*$ that interprets
each bracket of the Dyck word by a word in the original alphabet, with some arbitrary choice involved in selecting whether to interpret opening or closing brackets as empty words.

In this section, we give a new proof of the representation theorem,
generalized to context-free grammars $G$ over any category~$\Ccategory$.
Since the category $\Ccategory$ may have more than one object, the appropriate statement of the representation theorem cannot require the grammar describing the shapes of parse trees to have only one nonterminal.
Nonetheless, we can construct one that is \emph{$\Ccategory$-chromatic} in the following sense.
\begin{defi}
A categorical context-free grammar $G = (\cC,\sS,S,p)$ is \defin{$\cC$-chromatic}
if the functor $p : \FreeOper{\sS} \to\Words{\cC}$ is injective on colors.
\end{defi}
\noindent
The nonterminals of a $\cC$-chromatic grammar may thus be considered as pairs $(A,B)$ of objects of~$\cC$.

Moreover, rather than using Dyck words to represent parse trees, we find it more natural to use tree contour words, based on the observation given above in \S\ref{section/unigram}
 that every context-free language may be \emph{canonically} represented as the image of a tree contour language generated by a context-free grammar with the same set of nonterminals.
(We will discuss the relationship between contour words and Dyck words as Remark~\ref{remark/contour-vs-dyck} below.)

As preparation to our proof of the representation theorem, we establish:
\begin{prop}\label{proposition/node-colors-factorization}
Any map of species $\phi:\Sspecies \to \Tspecies$ factors as 
\begin{center}
\begin{tikzcd}[column sep=2em]
\Sspecies
\arrow[rr,"{\phi}"]
&&
\Tspecies
\quad = \quad
\Sspecies
\arrow[rr,"{\phi_{\colors}}"] 
&& 
\phi_C\,\sS
\arrow[rr,"{\phi_\nodes}"] &&
\Tspecies
\end{tikzcd}
\end{center}
for some species $\phi_C\,\sS$ such that $\phi_\colors$ is the identity on nodes and surjective on colors and $\phi_\nodes$ is injective on colors.
\end{prop}
\begin{proof}
  A map of species is given by a commutative diagram as on the right below:
\[
  \begin{tikzcd}
    \Sspecies \arrow[d,"\phi"] \\
    \Tspecies
  \end{tikzcd}
  \qquad
\begin{tikzcd}
  C_\Sspecies^*\arrow[d,"\phi_C^*"] & \arrow[l,"i"']V_\Sspecies\arrow[r,"o"]\arrow[d,"\phi_V"] & C_\Sspecies \arrow[d,"\phi_C"]\\
  C_\Tspecies^* & \arrow[l,"i'"']V_\Tspecies\arrow[r,"o'"] & C_\Tspecies
\end{tikzcd}
\]
After taking the epi-mono factorization of $\phi_C$,
\[
\begin{tikzcd}
C_\Sspecies
\arrow[r,"\phi_C"]
&
C_\Tspecies
\end{tikzcd}
=
\begin{tikzcd}
C_\Sspecies
\arrow[two heads,r,"e"] 
& 
\phi_C\,C_\Sspecies
\arrow[tail,r,"m"] &
C_\Tspecies
\end{tikzcd}
\]
we obtain the desired factorization of $\phi$:
\[
  \begin{tikzcd}
    \Sspecies \arrow[d,"\phi_\colors"] \\
    \phi_C\,\Sspecies\arrow[d,"\phi_\nodes"] \\
    \Tspecies
  \end{tikzcd}
  \qquad
\begin{tikzcd}
  C_\Sspecies^*\arrow[d,two heads,"e^*"] & \arrow[l,"i"']V_\Sspecies\arrow[r,"o"]\arrow[d,equals] & C_\Sspecies \arrow[d,two heads,"e"]\\
  (\phi_C\, C_\Sspecies)^*\arrow[d,tail,"m^*"] & \arrow[l,"i e^*"']V_\Sspecies\arrow[r,"oe"]\arrow[d,"\phi_V"] & \phi_C\,C_\Sspecies \arrow[d,tail,"m"]\\
  C_\Tspecies^* & \arrow[l,"i'"']V_\Tspecies\arrow[r,"o'"] & C_\Tspecies
\end{tikzcd} \qedhere
\]
\end{proof}
\begin{lem}\label{lemma/pullback-sliced-contour}
Let $\psi:\Sspecies\to\Sspecies'$ be a map of species, let $p_\sS : \FreeOper{\sS} \to \Words{\Contour{\FreeOper{\sS}}}$ and $p_{\sS'} : \FreeOper{\sS'} \to \Words{\Contour{\FreeOper{\sS'}}}$ be the underlying functors of the universal grammars derived from the unit of the contour / splicing adjunction, and let
\begin{equation}\label{naturality-square-Oper}
\begin{tikzcd}[column sep=4em,row sep=2em]
{\FreeOper{\Sspecies}}
\arrow[d,"{p_{\Sspecies}}"{swap}]\arrow[r,"{\FreeOper{\psi}}"] 
& 
{\FreeOper{\Sspecies'}}
\arrow[d,"{p_{\Sspecies'}}"]
\\
{\Words{\Contour{\FreeOper{\Sspecies}}}}
\arrow[r,"{\Words{\Contour{\FreeOper{\psi}}}}"{swap}] &
{\Words{\Contour{\FreeOper{\Sspecies'}}}}
\end{tikzcd}
\end{equation}
be the associated naturality square.
If $\psi$ is injective on nodes then the canonical functor of operads from $\FreeOper{\Sspecies}$ to the pullback of $p_{\Sspecies'}$ along $\Words{\Contour{\FreeOper{\psi}}}$ is fully faithful.
\end{lem}
\begin{proof}
The naturality square \eqref{naturality-square-Oper} in $\Operads{}$ induces a commutative square in $\Species{}$:
\begin{equation}\label{naturality-square-Spec}
\begin{tikzcd}[column sep=4em,row sep=2em]
{\Sspecies}
\arrow[d,"{\phi_{\Sspecies}}"{swap}]\arrow[r,"{\psi}"]
&
{\Sspecies'}
\arrow[d,"{\phi_{\Sspecies'}}"]
\\
{\Words{\Contour{\FreeOper{\Sspecies}}}}
\arrow[r,"{\Words{\Contour{\FreeOper{\psi}}}}"{swap}] &
{\Words{\Contour{\FreeOper{\Sspecies'}}}}
\end{tikzcd}
\end{equation}
Neither \eqref{naturality-square-Oper} nor \eqref{naturality-square-Spec} will in general be a pullback square.
However, by Lemma~\ref{lemma/pullback-free-along-ULF}, we know that the pullback of $p_{\Sspecies'}$ along $\Words{\Contour{\FreeOper{\psi}}}$
is obtained from a corresponding pullback in the category of species
\begin{equation}\label{pullbackR}
\begin{tikzcd}[column sep=4em,row sep=2em]
{\Tspecies}
\arrow[d,"{\rho}"{swap}]\arrow[r,"{\pi}"]
\arrow[rd, phantom, "\lrcorner" very near start]
&
{\Sspecies'}
\arrow[d,"{\phi_{\Sspecies'}}"]
\\
{\Words{\Contour{\FreeOper{\Sspecies}}}}
\arrow[r,"{\Words{\Contour{\FreeOper{\psi}}}}"{swap}] &
{\Words{\Contour{\FreeOper{\Sspecies'}}}}
\end{tikzcd}
\qquad\vrule\ \vrule\qquad
\begin{tikzcd}[column sep=4em,row sep=2em]
  \FreeOper{\Tspecies}\arrow[d,"{r}"{swap}]\arrow[r,"\FreeOper{\pi}"]
  \arrow[rd, phantom, "\lrcorner" very near start]
  & \FreeOper{\Sspecies'}\arrow[d,"p_{\Sspecies'}"]
\\
{\Words{\Contour{\FreeOper{\Sspecies}}}}\arrow[r,"{\Words{\Contour{\FreeOper{\psi}}}}"{swap}] & {\Words{\Contour{\FreeOper{\Sspecies'}}}}
\end{tikzcd}
\end{equation}
where the pullback species $\Tspecies$ has colors defined as tuples $(R,(R_1^{\tagU},R_2^{\tagD}))$
of a color $R$ of $\sS'$ and colors $R_1$, $R_2$ of $\sS$ such that ${\psi(R_1)=\psi(R_2)=R}$, and
where the $n$-ary nodes of $\sR$ are defined as pairs $(x,f)$
of a $n$-ary node $x$ of $\Sspecies'$ and an $n$-ary operation $f$ of $\Words{\Contour{\FreeOper\Sspecies}}$,
necessarily of the form $f = (y,0){-}\dots{-}(y,n)$
for $y$ the unique $n$-ary node of $\Sspecies$ such that $\psi(y)=x$,
by assumption that $\psi:\Sspecies\to\Sspecies'$ is injective on nodes.
The canonical map of species ${\Sspecies}\to{\Tspecies}$ induced by \eqref{naturality-square-Spec} and the left-hand side of \eqref{pullbackR} transports every color $R$ of~$\Sspecies$ to the color~$((\psi(R)^{\tagU},\psi(R)^{\tagD}),R)$
and every $n$-ary node $y$ of~$\Sspecies$ to the $n$-ary node $(\psi(y),(y,0){-}\dots{-}(y,n))$ of $\Tspecies$.
From this follows that the canonical map of species ${\Sspecies}\to{\Tspecies}$ is injective on colors and bijective on nodes.
Moreover, there are no nodes in $\Tspecies$ whose colors are outside of the image of $\Sspecies$.
We conclude that the canonical functor of operads $\FreeOper{\Sspecies}\to\FreeOper{\Tspecies}$ induced by \eqref{naturality-square-Oper} and the right-hand side of \eqref{pullbackR} is fully faithful.
\end{proof}
\noindent
Now, let $G=(\Ccategory,\Sspecies,S,p)$ be any categorical context-free grammar, and using Prop.~\ref{proposition/node-colors-factorization} consider the corresponding $\Ccategory$-chromatic grammar $G_\nodes = (\Ccategory,\phi_C\,\Sspecies,(A,B),p_\nodes)$, where $p(S) = (A,B)$.
Applying the factorization \eqref{equation/universal-CFG-factorization}, we build the commutative diagram below
\begin{equation}\label{equation/pullback-of-C-chromatic}
\begin{tikzcd}[column sep=3em,row sep=1em]
{\FreeOper{\Sspecies}}
\arrow[dd,"{p_{\Sspecies}}"{swap}]\arrow[rr,"{\FreeOper{\phi_\colors}}"]
&& 
{\FreeOper{\phi_C\,\Sspecies}}
\arrow[dd,"{p_{\phi_C\,\Sspecies}}"]
\\
&  &
\\
{\Words{\Contour{\FreeOper{\Sspecies}}}}
\arrow[rr,"{\Words{\Contour{\FreeOper{\phi_\colors}}}}"{}] 
\arrow[rdd,"{\Words{q_{G}}}"{swap}] 
&&
{\Words{\Contour{\FreeOper{\phi_C\,\Sspecies}}}}
\arrow[ldd,"{\Words{q_{G_\nodes}}}"] 
\\\\
&
\Words{\Ccategory}
\end{tikzcd}
\end{equation}
combining a naturality square with a triangle whose commutativity follows from equation~$\phi = \phi_\colors\,\phi_\nodes$
and the contour / splicing adjunction.
Note also that $p_\colors \defeq \Contour{\FreeOper{\phi_\colors}}$ is a ULF functor of categories by Prop.~\ref{proposition/Contour-species-ULF}, and finitary
because $\phi_\colors$ is finitary (and even finite).
Therefore pairing the functor with the initial state $S^\tagU$ and accepting state $S^\tagD$ defines a finite-state automaton $M_\colors = (\Contour{\FreeOper{\phi_C\,\Sspecies}},\Contour{\FreeOper{\Sspecies}},p_\colors,S^\tagU,S^\tagD)$.
By the intersection construction (Cor.~\ref{corollary/intersection-of-languages}) together with Lemma~\ref{lemma/pullback-sliced-contour} and the strong translation principle (Prop.~\ref{translation-principle}),
we deduce that
%
$$
p_\colors \, \ContourLang{\dot\Sspecies} \, = \, \ContourLang{\phi_C\,\dot\Sspecies} \, \cap \, \Lang{M_\colors}\,.
$$
Finally, using that $G$ is the image of the universal grammar $\UnivGrammar{\dot\Sspecies}$ and considering the commutative diagram~\eqref{equation/pullback-of-C-chromatic}, we conclude that
\begin{equation}\label{equation/explicit-CS-rep}
  \Lang{G} = q_G\,\ContourLang{\dot\Sspecies} = q_{G_\nodes}\,p_\colors \, \ContourLang{\dot\Sspecies} = q_{G_\nodes}\, (\ContourLang{\phi_C\,\dot\Sspecies} \, \cap \, \Lang{M_\colors})
\end{equation}
and thus derive our generalized Chomsky-Schützenberger representation theorem.
\begin{thm}\label{theorem/CS}
A language $L \subseteq \cC(A,B)$ is a CFL of arrows iff it is a functorial image of the intersection of a $\Ccategory$-chromatic tree contour language with a regular language.
\end{thm}
\noindent
\begin{exa}
  For the original example of Figure~\ref{fig:example-cfg}, the colors-nodes factorization
\begin{center}
\begin{tikzcd}[column sep=2em]
\Sspecies
\arrow[rr,"{\phi}"]
&&
\Words{\Sigma}
\quad = \quad
\Sspecies
\arrow[rr,"{\phi_{\colors}}"] 
&& 
\phi_C\,\Sspecies
\arrow[rr,"{\phi_\nodes}"] &&
\Words{\Sigma}
\end{tikzcd}
\end{center}
  of Prop.~\ref{proposition/node-colors-factorization} may be visualized as below:
  \begin{center}\includegraphics[scale=0.25]{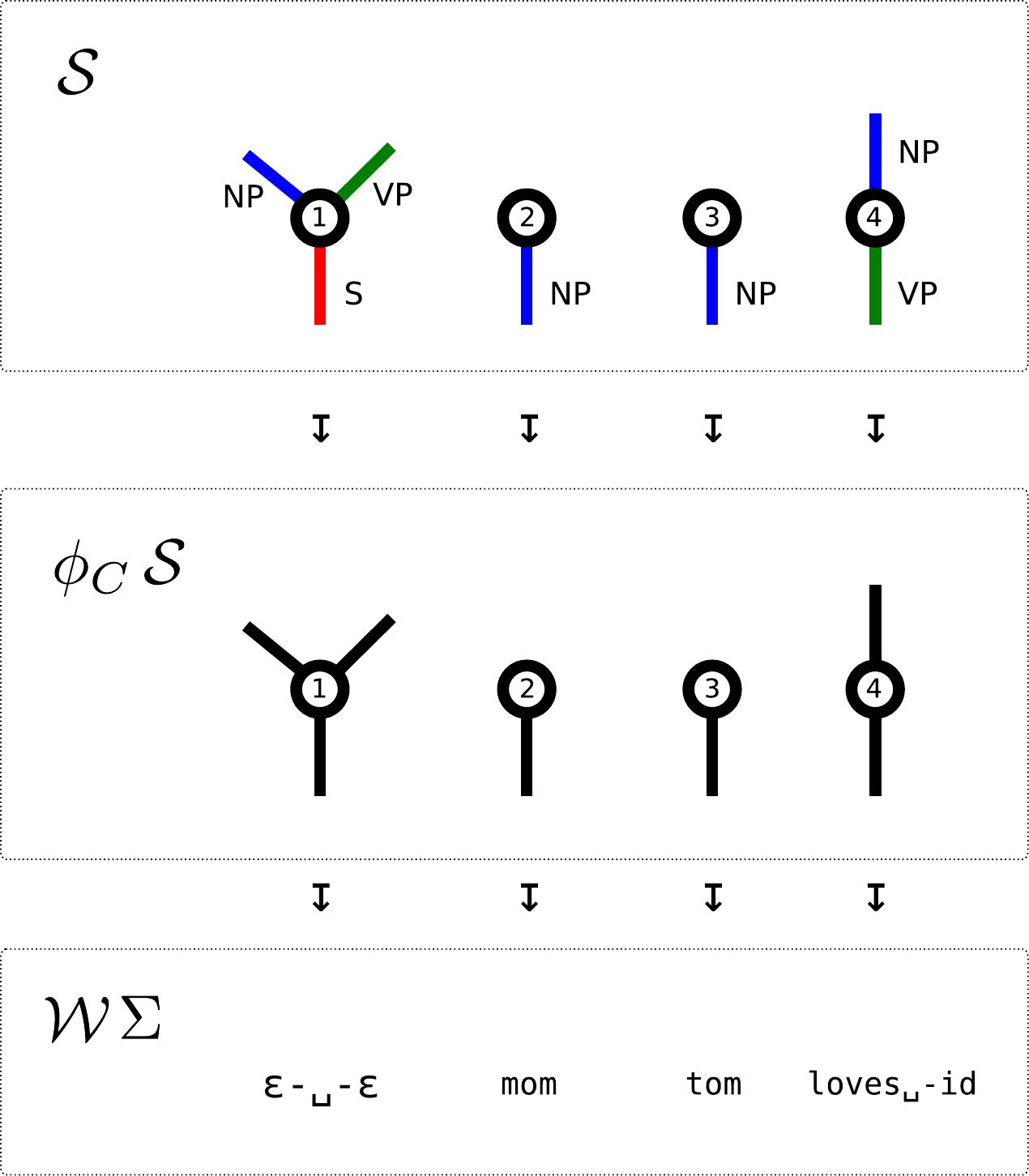}\end{center}
  Some examples of contour words in the tree contour language $\ContourLang{\phi_C\,\dot\Sspecies}$ and their corresponding images under the functor $q_{G_\nodes}$ include:
  \begin{center}
  \includegraphics[scale=0.25]{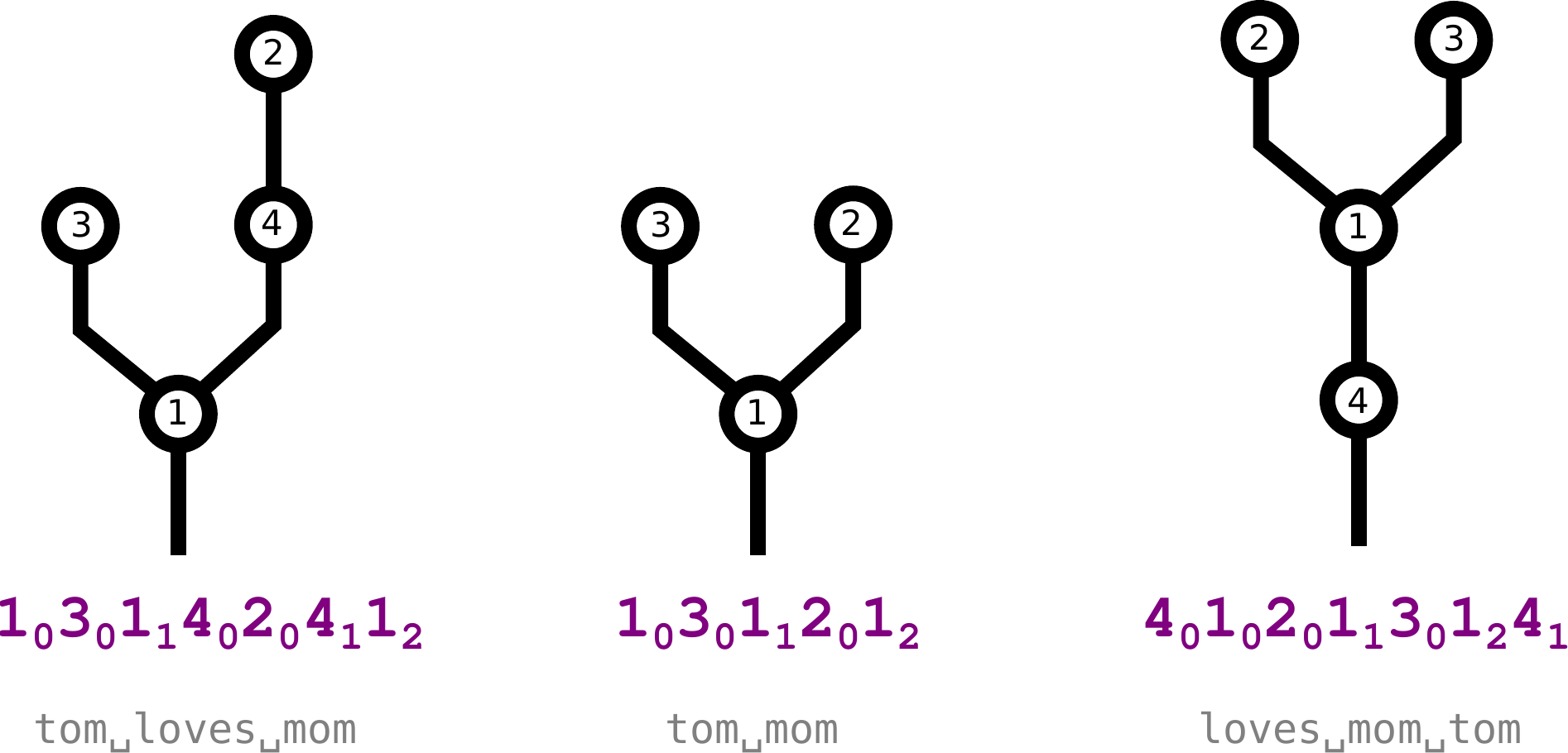}
  \end{center}
  Observe that only the first contour word encodes a valid parse tree mapping to a word of the language.
  Finally, the ``coloring automaton'' $M_\colors$ derived by the contour construction may be depicted by the following transition diagram:
  \[
    \begin{tikzpicture}[node distance=3cm,thick,scale=0.75, state/.style={circle, draw, minimum size=1cm}, every node/.style={transform shape}]
      \node[state,initial] (Su) at (0,0) {$\mathsf{S}^\tagU$};
      \node[state,right of=Su] (NPu) {$\mathsf{NP}^\tagU$};
      \node[state,below right of=NPu] (VPu) {$\mathsf{VP}^\tagU$};
      \node[state,above right of=VPu] (NPd) {$\mathsf{NP}^\tagD$};
      \node[state,right of=NPd] (VPd) {$\mathsf{VP}^\tagD$};
      \node[state,right of=VPd,accepting] (Sd) {$\mathsf{S}^\tagD$};
      \draw
      (Su) edge node[above]{$\mathsf{1_0}$} (NPu)
      (NPu) edge [bend left] node[above]{$\mathsf{2_0}$} (NPd)
      (NPu) edge [bend right] node[below]{$\mathsf{3_0}$} (NPd)
      (NPd) edge node[below right]{$\mathsf{1_1}$} (VPu)
      (VPu) edge node[below left]{$\mathsf{4_0}$} (NPu)
      (NPd) edge node[above]{$\mathsf{4_1}$} (VPd)
      (VPd) edge node[above]{$\mathsf{1_2}$} (Sd);
    \end{tikzpicture}
  \]
  The reader can verify that the automaton only accepts the first contour word above.
\end{exa}
\noindent
\begin{rem}\label{remark/contour-vs-dyck}
We view the use of contour words rather than Dyck words in Theorem~\ref{theorem/CS} as a mild improvement of the original theorem by Chomsky and Sch\"utzenberger, since it removes the need for choosing arbitrarily whether to erase left brackets or right brackets in defining the homomorphism.
However, there is a close relationship between tree contour words and Dyck words, and the original representation can be recovered from formula~\eqref{equation/explicit-CS-rep}.
We now sketch this for the classical case of context-free languages of words over an alphabet~$\Sigma$. 
As suggested by the right-hand side of Figure~\ref{fig:contour-word}, there is a simple mapping $s$ (for \emph{subdivide}) from contour words to Dyck words defined by translating each corner $(x,i)$ of an $n$-ary node $x$ into a pair of brackets, as follows:
\[
s\quad = \quad (x,i) \mapsto \begin{cases}[^x_{0}\ ]^x_{0} & \text{if }i = n = 0\\ [^x_{0}\ [^x_{1} &  \text{if }i = 0\text{ and } n > 0 \\ ]^x_{i}\ [^x_{i+1} & \text{if }0 < i < n \\ ]^x_{n}\ ]^x_{0} & \text{if }0 < i = n\end{cases}
\]
It is easy to check that the image of a tree contour word under this translation is a well-bracketed word over the alphabet containing a pair of brackets $[^x_i$ and $]^x_i$ for every $n$-ary node $x$ and index $0 \le i \le n$.
Indeed, every production rule
\[
R \to (x,0)\:R_1\:(x,1)\dots R_n\:(x,n)
\]
of the universal grammar $\UnivGrammar{\dot\sS}$ is translated to a production rule
\begin{equation}\label{S-dyck-production}
R \to [^x_0\:[^x_1\: R_1\:]^x_1\:[^x_2 \cdots R_n\:]^x_n\:]^x_0
\end{equation}
which manifestly preserves well-bracketing.
Let us refer to the language of well-bracketed words generated by the collection of productions \eqref{S-dyck-production}, indexed by nodes $x : R_1,\dots,R_n \to R$ of a species $\sS$ with start symbol $S \in \sS$, as the \emph{$\dot\sS$-Dyck language} $\DyckLang{\dot\sS}$.
By construction, it is immediate that
\begin{equation}\label{s-contour-is-S-dyck}
s(\ContourLang{\dot\sS}) = \DyckLang{\dot\sS}
\end{equation}
for every pointed species $\dot\sS$.
Conversely, there are two natural ways to define a translation from $\dot\sS$-Dyck words back to $\dot\sS$-contour words: the ``green'' translation
\[
g\quad=\quad [^x_i \ \mapsto
\begin{cases} (x,0) & \text{if } i=0\\ \epsilon& \text{if } i>0\end{cases}
\qquad
]^x_i \mapsto 
\begin{cases}\epsilon & \text{if }i=0 \\ (x,i) & \text{if }i>0\end{cases}\quad,
\]
which corresponds to interpreting the green arrows in the right hand side of Figure~\ref{fig:contour-word} as corners and erasing the red arrows,
or alternatively the ``red'' translation
\[
r\quad=\quad [^x_i \ \mapsto
\begin{cases} \epsilon & \text{if } i=0\\ (x,i)& \text{if } i>0\end{cases}
\qquad
]^x_i \mapsto 
\begin{cases}(x,0) & \text{if }i=0 \\ \epsilon & \text{if }i>0\end{cases} \quad,
\]
which does the opposite.
Both translations are left inverse to the $s$-translation in the sense that $g \circ s = r \circ s = \id$.
It follows that if we define the homomorphism $h$ from bracketed words to $\Sigma^*$ to be either $h = q_{G_\nodes} \circ g$ or $h = q_{G_\nodes} \circ r$, then we can derive that
\begin{align}
  \Lang{G} &= q_{G_\nodes}\, (\ContourLang{\phi_C\,\dot\Sspecies} \, \cap \, \Lang{M_\colors}) & \text{by }\eqref{equation/explicit-CS-rep} \notag \\
           &= h\, (s\, (\ContourLang{\phi_C\,\dot\Sspecies} \, \cap \, \Lang{M_\colors})) & \text{since }g\circ s = r\circ s = \id \notag\\
           &= h\, (\DyckLang{\phi_C\,\dot\Sspecies} \, \cap \, s(\Lang{M_\colors})) & \text{by }\eqref{s-contour-is-S-dyck}\text{ and injectivity of }s. \label{equation/s-dyck-cs-theorem}
\end{align}
Next, it is not hard to verify that the image $s(\Lang{M_\colors})$ of the regular language of the coloring automaton under the $s$-translation is itself a regular language (just use a pair of transitions to simulate each transition of $M_\colors$).
We have therefore derived that every context-free language is a homomorphic image of the intersection of a $\phi_C\,\dot\sS$-Dyck language (for some $\dot\sS$) with a regular language.

Words in $\DyckLang{\phi_C\,\dot\sS}$ are generated by productions of the form
\begin{equation}\label{S-dyck-production-uncolored}
S \to [^x_0\:[^x_1\: S\:]^x_1\:[^x_2 \cdots S\:]^x_n\:]^x_0
\end{equation}
for every node $x : R_1,\dots,R_n \to R$ in $\sS$ (i.e., replacing every non-terminal in \eqref{S-dyck-production} by the unique start symbol $S$).
Now let $k = \sum_{x\in\sS}1+|x|$ be the sum of the arity-plus-one of every node in $\sS$, and let us write $\DyckLang{k}$ for the language of unrestricted Dyck words over the family of $k$ different pairs of brackets labelled $[^x_i$ and $]^x_i$ for $0 \le i \le n = |x|$ as above.
One can prove by induction that an arbitrary Dyck word is a $(\phi_C\,\dot\sS)$-Dyck word just in case it satisfies the following conditions:
\begin{enumerate}
\item it starts with $[^x_0$ for some $x$;
\item every $[^x_0$ is immediately followed by either $]^x_0$ (if $|x|=0$) or $[^x_1$ (if $|x|>0$);
\item for $i>0$, every $[^x_i$ is followed by an opening bracket $[^y_0$ for some $y$;
\item for $i>0$, every $]^x_i$ is immediately followed by either $[^x_{i+1}$ (if $i < |x|$) or $]^x_0$ (if $i = |x|$);
\item $]^x_0$ is never followed by an opening bracket $[^y_i$.
\end{enumerate}
Since these conditions can be checked by a simple finite automaton $M_{\mathsf{index}}$, we have that
\begin{equation}\label{S-dyck-dyck}
\DyckLang{\phi_C\,\dot\sS} = \DyckLang{k} \cap \Lang{M_{\mathsf{index}}}
\end{equation}
from which the Chomsky-Schützenberger representation theorem follows by substituting \eqref{S-dyck-dyck} into \eqref{equation/s-dyck-cs-theorem} and using that regular languages are closed under intersection.

Although we have limited our attention here to the classical result, we believe that the connection between contour words and Dyck words deserves a deeper categorical analysis, which we leave to future work.
The $s$-translation has a clear geometric interpretation as decomposing each corner of the contour into alternating actions of walking along an edge and turning around a node (we again refer the reader to Figure~\ref{fig:contour-word}).
In the theory of maps on surfaces, an analogous decomposition is used to embed the \emph{oriented cartographic group} into the \emph{unoriented cartographic group} (cf.~\cite{ShabatVoevodsky1990,JonesSingerman1994}), which suggests that there may be a corresponding refinement of the contour category construction.
\end{rem}

\subsection{Related work}
\label{section/csrep/related}

As already mentioned, our proof that (generalized) context-free languages are closed under intersection with regular languages may be seen as a decomposition of the Bar-Hillel construction \cite{BarHillel+1961}.
In a similar spirit, Ludmann, Pogodalla, and De~Groote have recently analyzed Kanazawa's proof that languages generated by abstract categorial grammars are closed under intersection with regular languages, and shown that it defines a pullback in an appropriate category \cite{LudmannPogodallaDeGroote2022}.

The representation theorem is one of many results appearing in Chomsky and Schützenberger's joint article on ``The algebraic theory of context-free languges'' \cite{ChomskySchuetzenberger1963}.
The theorem was sketched in a technical report by Chomsky \cite{Chomsky1962contextfree} with attribution to their collaboration, while the first detailed proof seems to have been by Schützenberger \cite{Schutzenberger1963context}.
The original formulation used \emph{two-sided}\footnote{The two-sided Dyck language on an alphabet of $2n$ letters $x_{\pm i}$ $(1 \le i \le n)$ may be defined as the kernel of the canonical homomorphism $\phi : \mathrm{Mon}\langle x_{-n},\dots,x_{-1},x_{+1},\dots,x_{+n}\rangle \to \mathrm{Grp}\langle x_1,\dots,x_n\rangle$ from the free monoid on $2n$ generators to the free group on $n$ generators defined by $\phi(x_{+i}) = x_i$ and $\phi(x_{-i}) = x_i^{-1}$.} Dyck languages, but would have worked just as well with ordinary (one-sided) Dyck languages and most standard treatments of the theorem now adopt that formulation \cite[\S{G}]{Kozen1997book}.
Numerous variations of the representation theorem and extensions to other classes of languages have since been established, see \cite{Crespi2016talk} for an overview and bibliography.

The contour / splicing adjunction is fundamental to our analysis of the representation theorem, providing an unexpected geometric lens on context-free languages that is evocative of Girard's geometry of interaction. 
Although the adjunction is not identified, this geometric perspective is also apparent in Slavnov's recent work~\cite{Slavnov2020} (inspired both by abstract categorial grammars and by linear logic proof-nets), wherein he constructs a compact closed monoidal category of \emph{word cobordisms} reminiscent of the operad of spliced words.

After the original version of this article was published \cite{MZ2022mfps}, the contour / splicing adjunction was subsequently analyzed and extended by Earnshaw, Hefford, and Román, who described applications to process algebra \cite{EarnshawHR24,Roman2023phd}.
They noticed that the splicing functor $\Wordsonly : \Cat \to \Operads{}$ is closely related to a construction by Day and Street~\cite[Example~7.3]{DayStreet2004}, of a canonical promonoidal structure on $\cC^\op \otimes \cC$ for any $\ccat{V}$-enriched category~$\cC$.
We have also been further exploring the 2-categorical aspects of the adjunction between categories and operads in ongoing work with Peter Faul.


\section{Conclusion}
\label{section/conclusion}

We have seen that the classical notions of context-free grammar and of nondeterministic finite-state automaton can be reformulated and generalized in a categorical setting, in a way that seems to have significant explanatory power.
Clearly many natural questions arise in regards to integrating the deep body of work on automata theory and formal language theory within this setting, some of which we have already signalled.
Another obvious question is how to define pushdown automata over categories.
We find it remarkable that when one considers the nodes of a species as the production rules of a context-free grammar, the corners of its contour category correspond exactly to the \emph{items} used in LR parsing and Earley parsing \cite{Knuth1965,Earley1970}, which we think suggests contour categories have a role to play in analyzing left-to-right parsing algorithms from a categorical perspective.
One of our original motivations for studying context-free grammars as functors and looking at parsing as a lifting problem was to develop a better understanding of the universal properties of algorithms like these, having in mind that such an understanding could then be transferred to other settings including type checking and type inference algorithms.
This attempt at unification is supported by the emergence of a number of categorical and fibrational structures which seem to connect these different kinds of symbolic systems, revealing that they are founded on many of the same basic patterns.




\section*{Acknowledgment}
We thank Peter~Faul for many exciting discussions of the contour / splicing adjunction, Thea~Li for stimulating conversations about determinization and $\epsilon$-removal for categorical automata, and Farzad~Jafarrahmani for fruitful interactions related to the material in \cite{MZcflsmods}.
We also gratefully acknowledge Bryce~Clarke for many helpful insights into the contents of the paper, which improved some of the proofs, and for telling us about \cite{Pare2012mealy}.
We thank Daniela~Petri\c{s}an, Thomas~Colcombet, and the other members of the working group on ``Categories for Automata and Language Theory'' at IRIF for helping us to better understand the relationship with \cite{ColcombetPetrisan2020} and with other works.
We thank Benjamin~Steinberg for telling us about \cite{Jones1996profinite}, which led us to \cite{TherienSznajder1988,WeissTherien1986}.
We gratefully thank the LMCS referees for numerous comments that improved the presentation in the article.
Finally, special thanks to Hayo~Thielecke for giving a very lucid presentation of LL and LR parsing using abstract machines in 2017, which was one of our original sources of inspiration.

\bibliographystyle{alphaurl}
\bibliography{refs}


\end{document}